\documentstyle[amssymb,amsfonts]{amsart}

\newenvironment{proof}{\noindent {\bf Proof} }{\endprf\par}
\def \endprf{\hfill  {\vrule height6pt width6pt depth0pt}\medskip}

\newcommand{\bea}{\begin{eqnarray}}
\newcommand{\eea}{\end{eqnarray}}
\def\beaa{\begin{eqnarray*}}
\def\eeaa{\end{eqnarray*}}
\def\bseq{\begin{subeq}}
\def\eseq{\end{subeq}}
\def\ba{\begin{array}}
\def\ea{\end{array}}
\def\be#1{\begin{equation} \label{#1}}

\theoremstyle{plain}
  \newtheorem{theorem}[subsection]{Theorem}
  
  \newtheorem{proposition}[subsection]{Proposition}
  \newtheorem{lemma}[subsection]{Lemma}

\theoremstyle{remark}

\theoremstyle{definition}

\begin{document}

	  \title{Global regularity of Wave Maps from $\bf{R}^{3+1}$ to $\mathbf{H}^{2}$}
  \author{Joachim Krieger}

  \vspace{-0.3in}
  \begin{abstract}
  We consider Wave Maps with
  smooth compactly supported initial data of small $\dot{H}^{3/2}$-norm
  from $\bf{R}^{3+1}$ to hyperbolic 2-space
  and show that they stay smooth globally in time. Our methods are
  based on the introduction of a global Coulomb Gauge as in
  \cite{Str-Sh}, followed by a dynamic separation as in \cite{Kl-M
  3}. We then rely on an adaptation of T.Tao's methods used in his
  recent breakthrough result \cite{Tao 2}.
  \end{abstract}

\maketitle
\section{Introduction}

Let $M$ be a Riemannian manifold with metric $(g_{ij})=g$. A Wave
Map $u:\mathbf{R}^{n+1}\rightarrow$ $M,\,n\geq 2$ is by definition
a solution of the Euler-Lagrange equations associated with the
functional $u\rightarrow \int_{\mathbf{R}^{n+1}}<\partial_{\alpha
}u,
\partial^{\alpha} u>_{g}d\sigma$. Here the usual Einstein summation
convention is in force, while $d\sigma$ denotes the volume measure
on $\mathbf{R}^{n+1}$ with respect to the standard metric. In
local coordinates, $u$ is seen to satisfy the equation
\\
\begin{equation}\label{local formulation}
\Box
u_{i}+\Gamma^{i}_{jk}(u)\partial_{\alpha}u_{j}\partial^{\alpha}u_{k}
= 0
\end{equation}
\\
where $\Gamma^{i}_{jk}$ refer to the Riemann-Christoffel symbols
associated with the metric $g$. \\
The relevance of this model problem arises from its connections
with more complex nonlinear wave equations of mathematical
physics: for example, Einstein's vacuum equations under
$U(1)$-symmetry attain the form of a Wave Maps equation coupled
with additional elliptic equations.
\\

We are interested in the well-posedness of the Cauchy problem for
\eqref{local formulation} with initial data $u[0]\times
\partial_{t}u[0]$ at time $t=0$ in $H^{s}\times H^{s-1}$.
Classical theory relying on the energy inequality and Sobolev
inequalities allows one to
deduce local well-posedness in $H^{s}$ for $s>\frac{n}{2}+1$. \\
Ideally, one would like to prove local well-posedness in
$H^{\frac{n}{2}}$, as this would immediately imply global in time
well-posedness. The reason for this is that
$\dot{H}^{\frac{n}{2}}$ is the Sobolev space invariant under the
natural scaling associated with \eqref{local formulation}.
Unfortunately, it is known that "strong well-posedness" in the
sense of analytic or even $C^{2}$-dependence on the initial data
fails at the $H^{\frac{n}{2}}$-level, $n\geq 3$. Thus the best
result to be hoped for is global regularity of Wave
Maps with initial data of small $\dot{H}^{\frac{n}{2}}$-norm.\\
In two space dimensions, the scale invariant Sobolev space
coincides with the classical $\dot{H}^{1}$, and numerical data as
well as the conjectured non-concentration of energy suggest global
regularity for Wave Maps with arbitrary smooth initial data,
\textit{provided the target is negatively curved}\footnote{This
has been proved by Struwe for Wave Maps to the hyperbolic plane
which are corotationally symmetric} . This underlines
the importance of the hyperbolic plane as target manifold.\\

In the quest for reaching the critical $\frac{n}{2}$ regularity,
local well-posedness for \eqref{local formulation}with initial
data in $H^{\frac{n}{2}+\epsilon},\,\epsilon >0$ was proved for
$n\geq 3$ by Klainerman and Machedon in \cite{Kl-M 1}, and for
$n=2$ in \cite{Kl-S 2}. Later, Tataru established local
well-posedness in the Besov space $B^{\frac{n}{2},1}$, \cite{Tat},
\cite{Tatu}. Note that $\dot{B}^{\frac{n}{2},1}$ has the same
scaling as $\dot{H}^{\frac{n}{2}}$, but unlike the latter controls
$L^{\infty}$.
\\

An important breakthrough with respect to global regularity was
recently achieved by T.Tao in the case of Wave Maps to the sphere
\cite{Tao 1}, \cite{Tao 2}, proving global regularity for smooth
initial data small in $\dot{H}^{\frac{n}{2}}$: Tao's work
exemplifies the importance of taking the global geometry of the
target into account, an aspect largely ignored by the local
formulation \eqref{local formulation}. Embedding the target sphere
in an ambient Euclidean space, the Wave Maps equation considered
by Tao takes the form
\\
\begin{equation}\label{sphere}
\Box\phi = -\phi\partial_{\alpha}\phi^{t}\partial^{\alpha}\phi =
-(\phi\partial_{\alpha}\phi^{t}-\partial_{\alpha}\phi\phi^{t})\partial^{\alpha}\phi
\end{equation}
\\
$\alpha$ as usual runs over the space-time indices $0,1,...n$. The
nonlinearity encodes both geometric (skew-symmetry of
$\phi\partial_{\alpha}\phi^{t}-\partial_{\alpha}\phi^{t}$) as well
as algebraic information ('null-form' structure). Tao manages to
analyze all possible frequency interactions of the nonlinearity up
to the case in which the derivatives fall on high frequency terms
while the undifferentiated term has very low frequency. This bad
case is then gauged away, using the skew-symmetric structure. With
this method, which served as inspiration for the following
developments, as well as sophisticated methods from harmonic
analysis, Tao manages to go all the way to $n=2$(note that the
smaller the dimension, the more difficult the problem on account
of the increasing scarcity of available Strichartz estimates).
\\

After Tao, Klainerman and Rodnianski \cite{Kl-R}, extended this
result to Wave Maps form $\mathbf{R}^{n+1}, n\geq 5$ to more
general and in particular noncompact targets. More precisely,
Klainerman and Rodnianski consider parallelizable targets which
are well-behaved at infinity. Upon introducing a global
orthonormal frame $\{e_{i}\}$, they define the new variables
$\phi^{i}_{\alpha}$ defined by $u_{*}\partial_{\alpha} =
\phi^{i}_{\alpha}e_{i}$. It turns out that these satisfy the
system of equations
\\
\begin{equation}\label{div-curl 1}
\partial_{\beta}\phi^{i}_{\alpha} -
\partial_{\alpha}\phi^{i}_{\beta} =
C^{i}_{jk}\phi^{j}_{\alpha}\phi^{k}_{\beta}
\end{equation}
\\
\begin{equation}\label{div-curl 2}
\partial_{\alpha}\phi^{i\alpha} =
-\Gamma^{i}_{jk}\phi^{j}_{\beta}\phi^{k}_{\gamma}m^{\beta\gamma}
\end{equation}
\\
where $m_{\beta\gamma}$ is the standard Minkowski metric on
$\mathbf{R}^{n+1}$ and $C^{i}_{jk},\Gamma^{i}_{jk}$ are defined as
follows:
\\
\begin{equation}
[e_{j},e_{k}] = C^{i}_{jk}e_{i}
\end{equation}
\\
\begin{equation}
\nabla_{e_{j}}e_{k} = \Gamma^{i}_{jk}e_{i}
\end{equation}
\\

There is again a skew-symmetric structure present in this
formulation on account of $\Gamma^{i}_{jk} = -\Gamma^{k}_{ji}$.
Moreover, by contrast with Tao's formulation \eqref{sphere}, the
boundedness of $\phi$ is replaced here by the boundedness of the
$C^{i}_{jk}, \Gamma^{i}_{jk}$. Klainerman and Rodnianski impose in
addition the condition that all derivatives of these coefficients
be bounded, or in their terminology that $M$
be 'boundedly parallelizable.'\\
If one now passes to the wave equation satisfied by the vector
$\phi_{\alpha}:=\{\phi^{i}_{\alpha}\}$, one obtains
\\
\begin{equation}\label{KlR}
\Box\phi_{\alpha} = -R_{\mu}\partial^{\mu}\phi_{\alpha} +E
\end{equation}
\\
where $R_{\mu}$ is skew-symmetric and moreover depends linearly on
$\phi$, provided we assume the $C^{i}_{jk}, \Gamma^{i}_{jk}$ to be
constant for simplicity's sake. $E$ is a cubic polynomial in
$\phi$. By contrast with \eqref{sphere}, the leading term in the nonlinearity is
'quadratic in $\phi$'.\\
It is now possible to control all possible frequency interactions
on the right hand side ($n\geq 5$) except when $R_{\mu}$ is
localized to very low frequency while $\partial^{\mu}\phi$ is at
large frequency. However, as Klainerman and Rodnianski observed,
the curvature
\\
\begin{equation}
\partial_{\nu}R_{\mu}-\partial_{\mu}R_{\nu}+[R_{\mu},R_{\nu}]
\end{equation}
\\
when $R$ is reduced to low frequencies is 'very small', in the
sense that it is quadratic in $\phi$, hence amenable to good
Strichartz estimates. To take advantage of this, they introduce a
Coulomb Gauge $\sum_{j=1}^{3}\partial_{j}\tilde{R}_{j}=0$, which
allows one to replace the $R_{\mu}$ in \eqref{KlR} by
$\tilde{R}_{\mu}$ which is 'quadratic in $\phi$', effectively
replacing the nonlinearity by a term which is trilinear in $\phi$
and hence easily handled by Strichartz estimates. The general
philosophy here is that the higher the degree of the nonlinearity,
the more room is available to apply Strichartz estimates.
Klainerman and Rodnianski's method is thus similar to Tao's in
that it utilizes a microlocal Gauge Change to deal with specific
bad frequency interactions.
\\

The last result to be mentioned in this development is the
simplification and extension of the previous arguments to include
the case of $4+1$-dimensional Wave Maps to esssentially arbitrary
targets achieved by Shatah-Struwe \cite{Str-Sh} and (in more
restrictive formulation) Uhlenbeck-Stefanov-Nahmod \cite{N-St-Uh}.
The former observed that using a Coulomb Gauge, in a similar
fashion as above, at the beginning without carrying out a
frequency decomposition allows one to reduce the nonlinearity to a
form directly amenable to Strichartz estimates. This allows them
to avoid the microlocal Gauge Change of Tao and leads to a
remarkable simplification of the argument. In addition, they are
also able to treat the case of dimension $4+1$.
\\

The methods in \cite{Kl-R} and \cite{Str-Sh} run into serious
difficulties for $3+1$-dimensional Wave Maps, and even more so for
$2+1$-dimensional Wave Maps.  This can be seen intuitively as
follows:\\
The global Coulomb Gauge puts the leading terms of the
nonlinearity roughly into the form $D^{-1}(\phi^{2})D\phi$. In
dimensions $4$ and higher, Shatah and Struwe can estimate such
terms relying on the Strichartz type inequality for Lorentz spaces
\\
\begin{equation}
||\phi||_{L_{t}^{2}L_{x}^{2n,2}}\leq C||\Box
\phi||_{L_{t}^{1}H^{\sigma}}+C||\phi[0]||_{H^{\frac{n}{2}-1}}
\end{equation}
\\
where $\sigma=\frac{n}{2}-2$. This can be used to estimate the
$L_{t}^{1}L_{x}^{\infty}$-norm of $D^{-1}(\phi^{2})$.
\footnote{Alternatively, as pointed out by Klainerman and
Rodnianski, one can utilize an improved bilinear version of
Strichartz estimates in \cite{Kl-Tat}
to handle these cases.}\\
However, in three space dimensions, the above estimate fails. In
order to handle the case when $D^{-1}(\phi^{2})$ has much lower
frequency that $D\phi$, one would have to use an endpoint
$L_{t}^{2}L_{x}^{\infty}$-Strichartz estimate, which is false,
even replacing the $L_{x}^{\infty}$-norm by $BMO$, see \cite{Tao
3} .
\\

The present paper starts with the basic formulation
\eqref{div-curl 1}, \eqref{div-curl 2} of Klainerman and
Rodnianski applied to the simplified context of $\mathbf{H}^{2}$,
but utilizes the Coulomb Gauge right at the beginning as do Shatah
and Struwe. The main innovation over the preceding then is to
introduce a special null-structure into the nonlinearity by way of
what we term a \textit{dynamic separation}\footnote{This
terminology was suggested by S.Klainerman}, a method introduced
first in \cite{Kl-M 3}: in our context, we introduce twisted
variables $\psi_{\alpha}: = e^{-i\Phi}\phi_{\alpha}$ in complex
notation for some potential function $\Phi$, and utilize the
div-curl system satisfied by these to split them into a dynamic
part, which has the form of a gradient, and an elliptic part,
which satisfies an elliptic div-curl system. Substituting these
components into the nonlinearity results in a fairly complicated
trilinear null-structure\footnote{This is to be contrasted with
the null-structure in \cite{Kl-M 3}, which is bilinear}, as well
as error terms which are at least
quintilinear in $\phi$. \\
In order to estimate the trilinear null-structure, we have to
utilize the technical framework set forth in \cite{Tao 2}. Also,
as we are at the level of the derivative with respect to $u$ in
\eqref{local formulation} (and hence lose one degree of
differentiability), we have to modify several key lemmata of Tao
for our purposes in the case of high-high interactions. Moreover,
we have to prove a Gauge Change estimate (Theorem~\ref{Gauge
Change}) which is new for the spaces introduced in \cite{Tao 2}.
This Gauge Change result clearly breaks down for $2+1$-dimensional
Wave Maps, and so there is little hope to extend our result to
that case by following the same route.\\
However, our result does extend to more general and in particular
higher-dimensional targets. Details on this shall be contained in
a forthcoming note.
\\

Our main theorem is the following: Denote the $2$-dimensional
hyperbolic plane by $\mathbf{H}^{2} =
\mathbf{R}\times\mathbf{R}_{>0}$. We utilize standard coordinates
$(x,y)\in\mathbf{R}\times\mathbf{R}_{>0}$, with respect to which
the metric becomes
$dg = \frac{dx^{2}+dy^{2}}{y^{2}}$. \\
A Wave Map will be described in terms of a pair of functions
$(x,y)$ on $\mathbf{R}^{3+1}$ with range as specified in the
preceding paragraph. Then we have the theorem
\\
\begin{theorem}\label{theorem:Main}
Let $M$ be the hyperbolic 2-plane. There exists a number
$\epsilon>0$ such that for all compactly supported smooth initial
data $(x,y)[0], (\partial_{t}x,\partial_{t}y)[0]$ with\\
$||(\frac{\nabla(x,y)}{y}[0],\frac{\partial_{t}(x,y)}{y}[0])||_{\dot{H}^{\frac{1}{2}}}<
\epsilon $, the WM-problem \eqref{local formulation} with these
initial data has a globally smooth solution.
\end{theorem}

\textit{Acknowledgments}: The author would like to thank his Ph.D.
advisor Sergiu Klainerman as well as Igor Rodnianski and Terence
Tao for helpful suggestions and comments as well as reading the
manuscript.\\
The research for this paper was conducted in the fall 2001.\newline

\newpage
\section{Outline of the Argument}
\subsection{Basic formulation of the problem}
This section will serve as outline for the rest of the paper.\\
Let us restate the main theorem, using the notation introduced at
the end of last section:
\\
\begin{theorem}\label{theorem:Main}
Let $M$ be the hyperbolic 2-plane. There exists a number
$\epsilon>0$ such that for all compactly supported smooth initial
data $(x,y)[0], (\partial_{t}x,\partial_{t}y)[0]$ with\\
$||(\frac{\nabla(x,y)}{y}[0],\frac{\partial_{t}(x,y)}{y}[0])||_{\dot{H}^{\frac{1}{2}}}<
\epsilon $, the WM-problem \eqref{local formulation} with these
initial data has a globally smooth solution.
\end{theorem}

We translate this problem to the level of the derivative,
utilizing the formulation \eqref{div-curl 1}, \eqref{div-curl 2}
with respect to the global orthonormal frame $\{y\partial_{x},
y\partial_{y}\}$. More explicitly, we have
\begin{equation}
\phi^{1}_{\alpha} =
\frac{\partial_{\alpha}x}{y},\,\phi^{2}_{\alpha} =
\frac{\partial_{\alpha}y}{y}
\end{equation}
\\
The div-curl system satisfied by these quantities is then of the
following form:
\\
\begin{equation}\label{1}
\partial_{\beta}\phi^{1}_{\alpha}-\partial_{\alpha}\phi^{1}_{\beta}=\phi^{1}_{\alpha}\phi^{2}_{\beta}-\phi^{1}_{\beta}\phi^{2}_{\alpha}
\end{equation}
\\
\begin{equation}\label{2}
\partial_{\beta}\phi^{2}_{\alpha}-\partial_{\alpha}\phi^{2}_{\beta}=0
\end{equation}
\\
\begin{equation}\label{3}
\partial_{\alpha}\phi^{1\alpha}=-\phi^{1}_{\alpha}\phi^{2\alpha}
\end{equation}
\\
\begin{equation}\label{4}
\partial_{\alpha}\phi^{2\alpha}=\phi^{1}_{\alpha}\phi^{1\alpha}
\end{equation}
\\
$\alpha,\beta$ here vary over the space-time indices $0,1,2,3$,
and Einstein's summation convention is in force.
\\

Once we can show that the $\phi^{i}_{\alpha}$ stay smooth globally
in time, the actual Wave Map can be obtained by integration from
$(\frac{\partial_{t}x}{y},\frac{\partial_{t}y}{y}) =
(\phi^{1}_{0},\phi^{2}_{0})$.\\
Letting $\phi_{\alpha}$ denote the column vector with entries
$\phi^{1}_{\alpha},\phi^{2}_{\alpha}$, we obtain the following
wave equations:
\\
\begin{equation}\label{simple formulation}
\Box\phi_{\alpha} = M_{\nu}\partial^{\nu}\phi_{\alpha} +
"\phi^{3}",
\end{equation}
\\
where
\\
\begin{equation}M_{\nu} =
\begin{pmatrix}0&-2\phi_{\nu}\\2\phi_{\nu}&0
\end{pmatrix},
\end{equation}
\\
and $"\phi^{3}"$ refers to a vector with entries that are cubic
polynomials in the $\phi^{i}_{\alpha}$. The fine structure of
these entries will actually be relevant later on, but we leave it
out for the present discussion.
\\

As explained in the introduction, this formulation does not lend
itself to good estimates.

\subsection{Introducing the global Coulomb Gauge}
We now try to modify the matrix $M_{\nu}$ by adding a term of the
form $2\partial_{\nu}A$, in such a way that the resulting matrix
$\tilde{M}_{\nu} = M_{\nu}+2\partial_{\nu}A$ has better
properties. More precisely, we want this to depend 'quadratically'
on $\phi$. This can be achieved by utilizing the Coulomb Gauge
condition $\sum_{j=1}^{3}\partial_{j}\tilde{M}_{j}=0$, whence $A =
-\frac{1}{2}\triangle_{x}^{-1}\sum_{j=1}^{3}\partial_{j}M_{j}$.\\
Indeed, observe that the $\tilde{M}_{\nu}$ satisfy the following
div-curl system:
\begin{equation}\label{div-curl}
\sum_{i=1}^{3}\partial_{i}\tilde{M}_{i} = 0,\,
\partial_{\nu}\tilde{M}_{\mu}-\partial_{\mu}\tilde{M}_{\nu} =
\begin{pmatrix}0&2(\phi^{1}_{\nu}\phi^{2}_{\mu}-\phi^{1}_{\mu}\phi^{2}_{\nu})\\
-2(\phi^{1}_{\nu}\phi^{2}_{\mu}-\phi^{1}_{\mu}\phi^{2}_{\nu})&0
\end{pmatrix},
\end{equation}
whence
\\
\begin{equation}
\tilde{M}_{\nu} =
\begin{pmatrix}0&2\sum_{i=1}^{3}\triangle^{-1}\partial_{i}(\phi^{1}_{\nu}\phi^{2}_{i}-\phi^{1}_{i}\phi^{2}_{\nu})\\
-2\sum_{i=1}^{3}\triangle^{-1}\partial_{i}(\phi^{1}_{\nu}\phi^{2}_{i}-\phi^{1}_{i}\phi^{2}_{\nu})&0
\end{pmatrix},
\end{equation}
\\
or in a first approximation $\tilde{M}_{\nu} = "D^{-1}(\phi^{2})"$.\\
We can now set $U = e^{A}$ and obtain
\begin{equation}\label{better formulation}
U^{-1}\Box(U\phi_{\alpha}) =
U^{-1}\Box(U)\phi_{\alpha}+\tilde{M}_{\nu}\partial^{\nu}\phi_{\alpha}+"\phi^{3}"
\end{equation}
Of course, we use the commutativity of the Gauge group for
2-dimensional target. The difference between this wave equation
for $U\phi_{\alpha}$ and \eqref{simple formulation} is that the
nonlinearity here consists of trilinear expressions. In
particular, this modification suffices to handle the case of
$4+1$-dimensional Wave Maps. For this, observe for example that
one can easily estimate the $L_{t}^{1}L_{x}^{2}$-norm of
$\tilde{M}_{\nu}\partial^{\nu}\phi_{\alpha}$ since this is morally
$D^{-1}(\phi^{2})D\phi$ and
\begin{equation}
||D^{-1}(\phi^{2})D\phi||_{L_{t}^{1}L_{x}^{2}} \leq C
||\phi||_{L_{t}^{2}L_{x}^{8,2}}^{2}||\phi||_{L_{t}^{\infty}H_{x}^{1}}
\end{equation}
The right-hand terms are controlled by means of Strichartz'
inequalities. Similarly, one can estimate the remaining terms of the nonlinearity in the $L_{t}^{1}L_{x}^{2}$-norm.\\
This is Shatah and Struwe's method for $\mathbf{H}^{2}$ . One can
also estimate this term using the improved bilinear Strichartz
estimate for $D^{-1}(\phi^{2})$ in \cite{Kl-Tat}, as observed by
Klainerman and Rodnianski.
\\

For the 3-dimensional case, Strichartz' estimates alone don't seem
sufficient. This can be seen by analyzing the case when
$D^{-1}(\phi^{2})$ has very low frequency while $D\phi$ has large
frequency; in order to recoup the exponential loss caused by
$D^{-1}$, one seems to be forced to employ a
$L_{t}^{2}L_{x}^{\infty}$ Strichartz estimate, which unfortunately
doesn't exist. To proceed, we need to take into account more of
the special structure of the nonlinear terms.

\subsection{Implementing the dynamic separation.}

For convenience's sake, and as we are dealing with $2$-dimensional target,
introduce complex notation:\\
Let $\phi_{\alpha} = \phi^{1}_{\alpha}+i\phi^{2}_{\alpha}$, and
define the twisted variables best suited for dynamic separation as
follows :
\\
\begin{equation}
\psi_{\alpha}:=\psi^{1}_{\alpha}+i\psi^{2}_{\alpha}=e^{-i\Phi}\phi_{\alpha},
\end{equation}
\\
where $\Phi:=\triangle^{-1}\sum_{k=1}^{3}\partial_{k}\phi^{1}_{k}$
($\triangle$ stands for $\triangle_{x}$). The $\psi_{\alpha}$
satisfy the following wave equation:
\\
\begin{equation}\label{WAVE}
\Box\psi_{\alpha} = [-i\Box\Phi-\partial_{\nu}\Phi
\partial^{\nu}\Phi]\psi_{\alpha}+2ie^{-i\Phi}\triangle^{-1}\sum_{k=1}^{3}\partial_{k}[\phi^{1}_{k}\phi^{2}_{\nu}-\phi^{2}_{k}\phi^{1}_{\nu}]\partial^{\nu}\phi_{\alpha}+e^{-i\Phi}"\phi^{3}"
\end{equation}
\\
which of course is just \eqref{better formulation} spelled out
more completely. At first glance the nonlinear terms here do not
seem to exhibit any straightforward null structure. To reveal it,
we want to decompose the variables into \textit{dynamic} and
\textit{elliptic} components. This allows us to rewrite the
nonlinear terms as sums of trilinear terms with a well-defined
null structure, and terms of higher degree of linearity.\\
The middle term on the right-hand side of \eqref{WAVE} will turn
out to be the most difficult, and indeed both low- and
high-frequency interactions seem forbidding. Before implementing
and using the 'dynamic separation' for the $\psi_{\alpha}$, we
need to express this middle term purely in terms of the
$\psi_{\alpha}$(up to error terms). For this, observe that we have
$\phi^{1}_{k}\phi^{2}_{\nu}-\phi^{2}_{k}\phi^{1}_{\nu} =
\Im(\phi_{k}\bar{\phi_{\nu}}) = \Im(\psi_{k}\bar{\psi_{\nu}})$,
whence we can rewrite it as
\\
\begin{equation}\label{villain2}
\sum_{j=1}^{3}\triangle^{-1}\partial_{j}[\psi^{1}_{j}\psi^{2}_{\nu}-\psi^{2}_{j}\psi^{1}_{\nu}]\partial^{\nu}\psi_{\alpha}
+\sum_{j=1}^{3}i\triangle^{-1}\partial_{j}[\psi^{1}_{j}\psi^{2}_{\nu}-\psi^{2}_{j}\psi^{1}_{\nu}]\partial^{\nu}\Phi\psi_{\alpha}
\end{equation}
\\
The 2nd term in this equation, being quadrilinear, is to be
considered an error term, and turns out to be actually fairly
easily controllable once the first term is dealt with.
\\

Now consider the div-curl system satisfied by the $\psi_{\alpha}$:
\\
\begin{equation}\label{Div-Curl}\begin{split}
&\partial_{\alpha}\psi_{\beta}-\partial_{\beta}\psi_{\alpha}=
(\partial_{\alpha}\phi_{\beta}-\partial_{\beta}\phi_{\alpha})e^{-i\Phi}
-i\sum_{j=1}^{3}(\phi_{\beta}\triangle^{-1}\partial_{j}\partial_{\alpha}\phi^{1}_{j}
-\phi_{\alpha}\triangle^{-1}\partial_{j}\partial_{\beta}\phi^{1}_{j})e^{-i\Phi}\\
&=i\psi_{\beta}\triangle^{-1}\partial_{j}(\psi^{1}_{\alpha}\psi^{2}_{j}-\psi^{2}_{\alpha}\psi^{1}_{j})
-i\psi_{\alpha}\triangle^{-1}\partial_{j}(\psi^{1}_{\beta}\psi^{2}_{j}-\psi^{2}_{\beta}\psi^{1}_{j})\\
\end{split}\end{equation}
\\

\begin{equation}\begin{split}
&\partial_{\nu}\psi^{\nu}=\partial_{\nu}\phi^{\nu}e^{-i\Phi}
-i\phi^{\nu}\triangle^{-1}\partial^{\nu}\partial_{j}\phi^{1}_{j}e^{-i\Phi}
=i\psi^{\nu}\triangle^{-1}\partial_{j}(\psi^{1}_{\nu}\psi^{2}_{j}-\psi^{2}_{\nu}\psi^{1}_{\j})\\
\end{split}\end{equation}
\\
Thus the right-hand side consists of trilinear expressions in
$\psi$, by contrast with the bilinear expressions on the
right-hand side of \eqref{1}-\eqref{4}.
\\

The dynamic separation for $\psi_{\alpha}$ now consists in
decomposing
\\
\begin{equation}\label{dynamic separation}
\psi_{\nu} = -R_{\nu}\sum_{k=1}^{3}R_{k}\psi_{k} +\chi_{\nu}: =
-R_{\nu}\Psi +\chi_{\nu}
\end{equation}
\\
$R_{\nu}$ here stands for the Riesz operator
$(\sqrt{-\triangle_{x}})^{-1}\partial_{\nu}$, $\nu = 0,1,2,3$. The
$\chi_{\nu}$("elliptic part") in turn are determined via a simple
elliptic div-curl system, namely
\begin{equation}
\partial_{j}\chi_{j} = 0
\end{equation}
\begin{equation}
\partial_{i}\chi_{\nu}-\partial_{\nu}\chi_{i}
=\partial_{i}\psi_{\nu}-\partial_{\nu}\psi_{i}
\end{equation}
\\
This in addition to \eqref{Div-Curl} implies immediately that we
can write
\\
\begin{equation}\label{elliptic}
\chi_{\nu} =
i\sum_{i,j=1}^{3}\partial_{i}\triangle^{-1}(\psi_{\nu}\triangle^{-1}\partial_{j}(\psi^{1}_{i}\psi^{2}_{j}-\psi^{1}_{j}\psi^{2}_{i})-\psi_{i}\triangle^{-1}\partial_{j}(\psi^{1}_{\nu}\psi^{2}_{j}-\psi^{1}_{j}\psi^{2}_{\nu}))
\end{equation}
\\
Passing to real parts and imaginary parts, we can write
$\psi^{1}_{\nu} = -R_{\nu}\Psi^{1}+\Re\chi_{\nu}$, $\psi^{2}_{\nu}
= -R_{\nu}\Psi^{2} +\Im\chi_{\nu}$, where $\Psi^{a} =
\sum_{k=1}^{3}R_{k}\psi^{a}_{k}$.
\\

Having implemented the dynamic separation, we can now see the
null-structure, corresponding to substituting the dynamic parts
$-R_{\nu}\Psi^{a}$ for $\psi^{a}_{\nu}$, $a=1,2$ in $[,]$ of the
first term of \eqref{villain2}. There results
\\
\begin{equation}\label{null-form}
\sum_{j=1}^{3}\triangle^{-1}\partial_{j}[R_{j}\Psi^{1}R_{\nu}\Psi^{2}-R_{j}\Psi^{2}R_{\nu}\Psi^{1}]\partial^{\nu}\psi_{\alpha}
\end{equation}
\\
A reason for calling this term a null-form, in addition to the
fact that it appears to intertwine a
$Q_{0}$-structure(corresponding to
$\partial_{\nu}\phi\partial^{\nu}\psi$ with a $Q_{\nu
j}$-structure(corresponding to
$\partial_{i}\phi\partial_{\nu}\psi-\partial_{\nu}\phi\partial_{i}\psi$),
is most directly exemplified by the following elementary identity,
which we state as an easily verified lemma:
\\
\begin{lemma}
Let $f,g,h$ be Schwartz functions. Then we have
\\
\begin{equation}\label{full null-form}\begin{split}
&2\sum_{j=1}^{3}\triangle^{-1}\partial_{j}[R_{\nu}fR_{j}g-R_{j}f
R_{\nu}g]\partial^{\nu}h \\
&=\sum_{j=1}^{3}\Box[\triangle^{-1}\partial_{j}[\nabla^{-1}fR_{j}g]h]-\sum_{j=1}^{3}\Box\triangle^{-1}\partial_{j}[\nabla^{-1}fR_{j}g]h\\
&-\sum_{j=1}^{3}\triangle^{-1}\partial_{j}[\nabla^{-1}fR_{j}g]\Box
h-\nabla^{-1}f\Box((\nabla^{-1}g)
h)\\&+\nabla^{-1}f\Box(\nabla^{-1}g)
h+\nabla^{-1}f(\nabla^{-1}g)\Box h\\
\end{split}\end{equation}
\\
\end{lemma}

\textit{Remark}: The bilinear null form in \cite{Kl-M 3} exhibits
similar structure, though our formulation, which avoids the
Fourier transform, is more simple and explicit.

This identity will become useful once we have reduced all entries
to small modulation, i.e. to have Fourier support close to the
light cone, for then the $\Box$-operators kick in. However, when
the expression has relatively large modulation(statements such as
this one are always by comparison with frequencies and modulations
of the entries), one can usually estimate it in the
$\dot{X}_{0}^{-\frac{1}{2},-\frac{1}{2},1}$-norm. Moreover,
restricting entries to 'large' modulation allows one to estimate
them in $L_{t}^{2}L_{x}^{2}$ (using
$\dot{X}_{k}^{\frac{1}{2},\frac{1}{2},\infty}$-spaces) which is often strategically advantageous. \\
The embedded $Q_{\nu j}$ structure will become useful when both
entries in the expression
$[R_{j}\Psi^{1}R_{\nu}\Psi^{2}-R_{\nu}\Psi^{1}R_{j}\Psi^{2}]$ have
very large frequency.\\
Unfortunately, in spite of the conceptual simplicity of this
reasoning, filling in the details is still a formidable task, in
part as one has so many cases to consider.

Now consider the error terms arising upon substituting at least
one $\chi_{\nu}$ for $\psi_{\nu}$ in the first term of
\eqref{villain2}. One thereby obtains
\begin{equation}\label{easy1}
\nabla^{-1}(\nabla^{-1}(\nabla^{-1}(\psi^{2})\psi)R_{\nu}\Psi^{a})\nabla_{x,t}\psi
\end{equation}
\\
\begin{equation}\label{easy2}
\nabla^{-1}(\nabla^{-1}(\nabla^{-1}(\psi^{2})\psi)\nabla^{-1}(\nabla^{-1}(\psi^{2})\psi))\nabla_{x,t}\psi
\end{equation}
\\
the terms being written schematically, and we have referred to
\eqref{elliptic}. These terms will be totally straightforward,
invoking the $L_{t}^{4}L_{x}^{4}$-Strichartz
estimate. \\
Therefore, most of the work for this paper will go into proving an
estimate for the crucial term \eqref{null-form}, as well as
proving that the 'twisted variables' $\psi_{\alpha}$ are
controlled, and in turn allow one to control the $\phi_{\alpha}$,
a nontrivial task on account of the complex Banach spaces used for
our argument, viz. the following sections.
\\

Finally, having disposed of the middle term in \eqref{WAVE}, we
also have to control the remaining terms. Each of these has its
own flavor, but there are connections between all of them. For
example, the term $(\Box\Phi)\psi_{\alpha}$ calls for expanding
$\Box\Phi$ by means of the wave equation \eqref{simple
formulation}, and then invoking the dynamic separation associated
with the basic div-curl system satisfied by the $\phi_{\alpha}$ in
order to introduce a simple $Q_{0}$-type null-structure. \\
This null-structure is of course already present in the next term
$\partial_{\nu}\Phi\partial^{\nu}\Phi\psi_{\alpha}$ in
\eqref{WAVE}, where one has
to use it only to deal with high-frequency interactions.\\
Finally, as to the term $"\phi^{3}"e^{-i\Phi}"$, this term will be
trivial to estimate provided all entries are at low frequencies.
Otherwise, one again has to apply dynamic separation, this time to
the variables $\phi_{\alpha}$, in order to split it into a
$Q_{0}$-type null-form and quadrilinear error terms, which turn
out to be altogether elementary. Here the fine structure of
$"\phi^{3}"$ will of course play a crucial role. \\
This then summarizes the basic strategy for estimating
\eqref{WAVE}.

\subsection{The Bootstrapping argument} In order to prove the
global regularity of the $\phi^{i}_{\alpha}$, we utilize a
bootstrapping argument, quite similar to the one in \cite{Tao 2}.
More precisely, we introduce certain translation invariant Banach
spaces $S[k]([-T,T]\times\mathbf{R}^{3})$,
$N[k]([-T,T]\times\mathbf{R}^{3})$,\,$k\in\mathbf{Z}, T>0 $ which
enjoy a list of remarkable properties. The norms
$||.||_{S[k]([-T,T]\times\mathbf{R}^{3})}$ will be used to
estimate the components at frequency $\sim 2^{k}$ of the
$\phi^{i}_{\alpha}$ which are known to be smooth on the time
interval $[-T,T]$, while the norms
$||.||_{N[k]([-T,T]\times\mathbf{R}^{3})}$ will be used to
estimate the components at frequency $\sim 2^{k}$ of the
nonlinearity, again restricted to and smooth on the time interval
$[-T,T]$. Of course, $||.||_{S[k]}$ will have to majorize the
energy $||.||_{\dot{H}^{\frac{1}{2}}}$ as well as a certain range
of Strichartz norms, all applied to
functions microlocalized at frequency $\sim 2^{k}$.\\
Our goal will be to bootstrap each of the norms
$||P_{k}\phi^{i}_{\alpha}||_{S[k]([-T,T]\times\mathbf{R}^{3})}$.
As a matter of fact, we will only have to bootstrap
$||P_{0}\phi^{i}_{\alpha}||_{S[0]([-T,T]\times\mathbf{R}^{3})}$,
because the $S[k]$ will be constructed compatible with 'dilations'
compatible with the div-curl system \eqref{1}-\eqref{4} : denoting
$\phi_{\lambda}:=2^{\lambda}\phi(x2^{\lambda})$, we will have
$||P_{k+\lambda}\phi_{\lambda}||_{S[k+\lambda]([-T,T]\times\mathbf{R}^{3})}
= ||P_{k}\phi||_{S[k]([-T,T]\times\mathbf{R}^{3})},
k,\lambda\in\mathbf{Z}$. Here $P_{k}$ denotes the Littlewood-Paley
projector to frequency $\sim 2^{k}$. A similar identity holds for
$N[k]([-T,T]\times\mathbf{R}^{3})$.

The $S[k]$ and $N[k]$ (leaving out the time-parameter $T$ for
simplicity's sake) will be related by the fundamental energy
inequality:
\begin{equation}\label{energy}
||P_{k}\phi||_{S[k]([-T,T]\times\bf{R}^{3})}\leq C[||\Box
P_{k}\phi||_{N[k]([-T,T]\times\bf{R}^{3})}+||P_{k}\phi[0]||_{\dot{H}^{\frac{1}{2}}\times\dot{H}^{-\frac{1}{2}}}]
\end{equation}
where $C$ is independent of $T$. In order to use this inequality,
we need to estimate the $N[k]$-norm of the nonlinearity. For this,
it will be important to us amongst other things that there are
\begin{enumerate}
\item null-form estimates of the form
\begin{equation}
||P_{0}[R_{\nu}P_{k_{1}}\phi \partial^{\nu}P_{k_{2}}\psi]||_{N[0]}
\leq
C2^{-\delta\max\{k_{1},0\}}||P_{k_{1}}\phi||_{S[k_{1}]}||P_{k_{2}}\psi||_{S[k_{2}]},\,\delta>0
\end{equation}
\item Bilinear estimates that make up for the missing
$L_{t}^{2}L_{x}^{\infty}$-estimates. These come about by using
null frame spaces, and have roughly the form
\begin{equation}
||P_{k_{1}}\phi P_{k_{2}}\psi||_{L_{t}^{2}L_{x}^{2}}\leq
C2^{\frac{k_{1}-k_{2}}{2}}||P_{k_{1}}\phi||_{S[k_{1}]}||P_{k_{2}}\psi||_{S[k_{2}]}
\end{equation}
provided $\phi,\psi$ are microlocalized on small caps whose
distance is at least comparable to their radius, and provided
their Fourier support lives fairly closely to the cone.
\item Trilinear estimates: the preceding observations will
play a role in establishing the genuine trilinear estimate
\begin{equation}\begin{split}
&||P_{0}\sum_{j=1}^{3}\triangle^{-1}\partial_{j}[R_{\nu}P_{k_{1}}\psi_{1}R_{j}P_{k_{2}}\psi_{2}-R_{j}P_{k_{1}}\psi_{1}R_{\nu}P_{k_{2}}\psi_{2}]\partial^{\nu}P_{k_{3}}\psi_{3}||_{N[0]}\\&\leq
C2^{-\delta_{1}|k_{1}-k_{2}|}2^{-\delta_{2}|k_{3}|}\prod||P_{k_{j}}\phi_{j}||_{S[k_{j}]},\,\delta_{1},\delta_{2}>0
\end{split}\end{equation}
This will be the crux of the paper.
\item The $S[k]$ have to be well-behaved under the Gauge Change. In
particular, we need an assertion of the form that provided
$||P_{k}\phi||_{S[k]}$ are small in a suitable sense, then so are
$||P_{k}[f(\nabla^{-1}\phi)\phi]||_{S[k]}$, where $\nabla^{-1}$
stands for a linear combination of operators of the form
$\triangle^{-1}\partial_{j}$, and $f(x)$ is a smooth function all
of whose derivatives are bounded.
\end{enumerate}

\section{Technical preparations}
The spaces $S[k], N[k]$ and many of their properties were
considered in Tao's seminal paper \cite{Tao 2}, although their
origins can be traced back to Tataru's \cite{Tatu}. Most of this
section(except the trilinear inequality and the Gauge Change
result) is due to these 2 authors; we will therefore be rather
brief with the
definitions.\\
First, we introduce Tao's concept of frequency envelope, as in
\cite{Tao 1},\cite{Tao 2}: for any Schwartz function $\psi$ on
$\mathbf{R}^{3}$, we consider the quantities
\begin{equation}
c_{a}:=(\sum_{k\in\bf{Z}}2^{-\sigma|a-k|}||P_{k}\psi||_{\dot{H}^{\frac{1}{2}}}^{2})^{\frac{1}{2}}
\end{equation}
Here $P_{k}, k\in\mathbf{Z}$ are the standard Littlewood-Paley
operators that localize to frequency $\sim 2^{k}$, i.e. they are
given by Fourier multipliers $m_{k}(|\xi|) =
m_{0}(\frac{|\xi|}{2^{k}})$, where $m_{0}(\lambda)$ is a smooth
function compactly supported within $\frac{1}{2}\leq \lambda\leq
2$
with\\ $\sum_{k\in\bf{Z}}m_{0}(\frac{\lambda}{2^{k}})=1,\,\lambda>0 $.\\
The $\sigma>0$ is chosen to be smaller than any of the exponential
decays occuring later in the paper. E.g. $\frac{1}{1000}$ would
suffice. We note that all of the generic constants $C$ occuring in
the sequel depend at most on this envelope exponent $\sigma$. \\
Note that
\begin{equation}
c_{k}2^{-\sigma|a-k|}\leq c_{a}\leq 2^{\sigma|a-k|}c_{k}
\end{equation}
as well as $\sum_{k\in\bf{Z}}c_{k}^{2}\leq C
||\psi||_{\dot{H}^{\frac{1}{2}}}^{2}$. \\
One reason for why this concept is useful is that provided we know
that the frequency localized components $P_{k}\rho$ for some other
Schwartz function $\rho$ on $\mathbf{R}^{3}$(think: the
time-evolved Wave Map) have $\dot{H}^{\frac{1}{2}}$-norm bounded
by a multiple $Cc_{k}$, we can immediately bound the
$\dot{H}^{\frac{1}{2}+\epsilon}$-norm of $\rho$ for $\epsilon>0$
small enough. This will allow us later to continue the Wave Map,
by referring to local well-posedness of the div-curl system
\eqref{1}-\eqref{4} in $H^{\frac{1}{2}+\epsilon}$, and finite
speed of propagation.
\\

We introduce the following norms on frequency localized Schwartz
functions on $\mathbf{R}^{3+1}$ for our bootstrapping argument:
for every $l>10$, choose a covering $K_{l}$ of $S^{2}$ by finitely
overlapping caps $\kappa$ of radius $2^{-l}$. This is to be chosen
such that
the set of concentric caps with half the radius still covers the sphere. Now let\\

\begin{equation}\begin{split}
&||\psi ||_{S[k]}:=\\&
||\nabla_{x,t}\psi||_{L_{t}^{\infty}\dot{H}_{x}^{-\frac{1}{2}}} +
||\nabla_{x,t}\psi||_{\dot{X}_{k}^{-\frac{1}{2},\frac{1}{2},\infty}}
+\sup_{\pm}\sup_{l>10}(\sum_{\kappa\in
K_{l}}||\tilde{P}_{k,\pm\kappa}Q^{\pm}_{<k-2l}\psi||_{S[k,\kappa]}^{2})^{\frac{1}{2}}\\
\end{split}\end{equation}
\\

where it is understood that $\psi$ lives at frequency $\sim 2^{k},
k\in\mathbf{Z}$. The operators $\tilde{P}_{k,\kappa}$ are given by
symbols $\tilde{m}_{k}(|\xi|)a_{\kappa}(\frac{\xi}{|\xi|})$, where
$a:S^{2}\rightarrow \bf{R}$ is a smooth function with support
contained in the concentric cap inside $\kappa$ with half the
radius of $\kappa$, and $\tilde{m}_{k}$ localizes frequency to
size $\sim 2^{k}$ and satisfies $\tilde{m}_{k}m_{k} = m_{k}$,
where $m_{k}$ is the multiplier chosen above.  We also require
that $\sum_{\kappa\in K_{l}}\tilde{P}_{k,\kappa}=\tilde{P}_{k}$,
the latter being defined in the obvious way.\\
$Q^{\pm}_{<k-2l}$ localizes the modulation, i.e. $||\tau|-|\xi||$,
to size $<2^{k-2l}$ and also restricts the Fourier support to
$\tau><0$, i.e. to the upper or lower half-space. More precisely,
it is given by the multiplier
$\sum_{i<k-2l}m_{i}(||\tau|-|\xi||)\chi_{>0}(\pm\tau)$.\\
The norm $||\phi||_{\dot{X}_{k}^{-\frac{1}{2},\frac{1}{2},1}}$
refers to
$2^{-\frac{k}{2}}\sum_{j\in\mathbf{Z}}2^{\frac{j}{2}}||Q_{j}\phi||_{L_{t}^{2}L_{x}^{2}}$.
\\

As we are in $3$ space dimensions, we strengthen $S[k,\kappa]$
with respect to Tao's original  definition to also contain
Strichartz norms such as $||.||_{L_{t}^{4}L_{x}^{4}},\,
2^{\frac{k}{6}}||.||_{L_{t}^{6}L_{x}^{3}},\,
2^{-\frac{k}{6}}||.||_{L_{t}^{3}L_{x}^{6}}$. Hence we define
\\
\begin{equation}\begin{split}
||\psi||_{S[k,\kappa]}:=&2^{\frac{k}{2}}||\psi||_{NFA^{*}[\kappa]}+|\kappa|^{-\frac{1}{2}}2^{-\frac{k}{2}}||\psi||_{PW[\kappa]}
+2^{\frac{k}{2}}||\psi||_{L_{t}^{\infty}L_{x}^{2}}\\&+
\sup_{\frac{1}{p}+\frac{1}{q}\leq
\frac{1}{2},p>2+\mu}2^{k(\frac{1}{p}+\frac{3}{q}-1)}||\psi||_{L_{t}^{p}L_{x}^{q}}\\
\end{split}\end{equation}
\\

for some very small $\mu>0$, e.g. $\mu<\frac{1}{1000}$ will do.
The definitions of the individual ingredients are as follows:
\begin{equation}
||\psi||_{NFA*[\kappa]}:=\sup_{\omega\notin
2\kappa}dist(\omega,\kappa)||\phi||_{L_{t_{\omega}}^{\infty}L_{x_{\omega}}^{2}}
\end{equation}
Here $(t_{\omega},x_{\omega})$ refer to null-frame coordinates,
i.e. $t_{\omega} = (t,x)\cdot \frac{1}{\sqrt{2}}(1,\omega),
x_{\omega}=(t,x)-t_{\omega}\frac{1}{\sqrt{2}}(1,\omega)$.
\\

Also, define $||,||_{PW[\kappa]}$ to be the norm associated with
the atomic Banach space whose atoms are the set $A$ of all
Schwartz functions $\psi$ with
$||\psi||_{L_{t_{\omega}^{2}L_{x_{\omega}}^{\infty}}}\leq 1$ for
some $\omega\in\kappa$. In other words,
\begin{equation}\begin{split}
||\psi||_{PW[\kappa]} = &\inf\{|\lambda|
|\exists\{0\leq\lambda_{i}\leq 1\}, \{\psi_{i}\}\subset A,1\leq
i\leq N\, s.t.\,
\sum_{i}\lambda_{i}=1,\\&\lambda\sum_{i}\lambda_{i}\psi_{i}=\psi\}\\
\end{split}\end{equation}
\\
Of course, the Banach space $S[k]$ is obtained by completing the
Schwartz functions on $\mathbf{R}^{3+1}$ with respect to
$||.||_{S[k]}$.

Next, we will place frequency localized pieces of the nonlinearity
into the following spaces $N[k]$, again introduced by Tao: they
are the atomic Banach spaces whose atoms are
\begin{enumerate}
\item Schwartz functions $F$ at frequency between $2^{k-4}$ and $2^{k+4}$ with
$||F||_{L_{t}^{1}L_{x}^{2}}\leq 2^{\frac{k}{2}}$.
\item Schwartz functions $F$ with frequency between $2^{k-4}$ and $2^{k+4}$
and modulation between $2^{j-5}$ and $2^{j+5}$ such that
$||F||_{L_{t}^{2}L_{x}^{2}}\leq 2^{\frac{j}{2}}2^{\frac{k}{2}}$.
\item Schwartz functions $F$ for which there exists a number
$l>10$ and Schwartz functions $F_{\kappa}$ with Fourier support in
the region $\{(\tau,\xi)|\pm\tau>0,\,||\tau|-|\xi||\leq
2^{k-2l-100}, 2^{k-4}\leq |\xi|\leq 2^{k+4},
\Theta\in\frac{1}{2}\kappa\}$ such that $F = \sum_{\kappa\in
K_{l}}F_{\kappa}$ and $(\sum_{\kappa\in
K_{l}}||F_{\kappa}||_{NFA[\kappa]}^{2})^{\frac{1}{2}}\leq
2^{\frac{k}{2}}$ . Here $\Theta = \frac{\tau\xi}{|\tau||\xi|}$ and
$NFA[\kappa]$ is the dual space of $NFA[\kappa]^{*}$, i.e.
\begin{equation}
||\psi||_{NFA[\kappa]}= \inf_{\omega\notin 2\kappa}
\frac{1}{dist(\omega,\kappa)}||\psi||_{L_{t_{\omega}}^{1}L_{x_{\omega}}^{2}}
\end{equation}
\end{enumerate}
The reason for introducing these spaces is that they allow one to
establish essential bilinear estimates, in particular the
following:
\begin{equation}\label{crux}
||\phi\psi||_{L_{t}^{2}L_{x}^{2}}\leq C
\frac{2^{\frac{k'}{2}}|\kappa'|^{\frac{1}{2}}}{dist(\kappa,\kappa')2^{\frac{k}{2}}}||\phi||_{S[k,\kappa]}||\psi||_{S[k',\kappa']}
\end{equation}
\begin{equation}
||\phi\psi||_{NFA[\kappa]}\leq
C\frac{2^{\frac{k'}{2}}|\kappa'|^{\frac{1}{2}}}{dist(\kappa,\kappa')}||\phi||_{L_{t}^{2}L_{x}^{2}|}||\psi||_{S[k',\kappa']}
\end{equation}
These follow immediately from the definitions, or as in \cite{Tao
2} .
\\
Another reason for including the null-frame version
$L^{2}_{t_{\omega}}L^{\infty}_{x_{\omega}}$ in $S[k,\kappa]$
instead of the forbidden $L_{t}^{2}L_{x}^{\infty}$ is that we have
the following inequality valid for all Schwartz functions with
Fourier support localized along $\kappa$ and at frequency$\sim
2^{k}$:
\\
\begin{equation}\label{cap localized}
||\phi||_{S[k,\kappa]}\leq
C||\phi||_{\dot{X}_{k}^{\frac{1}{2},\frac{1}{2},1}}
\end{equation}
\\
Of course, in light of the ingredients of $S[k]$,  this is a
substitute for the missing $L_{t}^{2}L_{x}^{\infty}$-Strichartz
estimate.
\\

Moreover, it is easily seen that the $S[k]$-norm of
$P_{k}\phi$ as defined above majorizes\\
$\sup_{\frac{1}{p}+\frac{1}{q}=
\frac{1}{2},K>p>2+\mu}2^{k(\frac{1}{p}+\frac{3}{q}-2)}||\nabla_{x,t}
P_{k}\phi||_{L_{t}^{p}L_{x}^{q}}$, $K$ an arbitrarily large
constant, up to multiplication with a constant(depending on $K$).
Indeed, splitting $P_{0}\phi = P_{0}Q_{<-100}\phi + P_{0}Q_{\geq
-100}\phi$, we can control the 2nd high-modulation summand by
invoking the fact that
\begin{equation}
||P_{0}Q_{j}\phi||_{L_{t}^{p}L_{x}^{q}}\leq C
2^{j(\frac{1}{2}-\frac{1}{p})}||P_{0}Q_{j}\phi||_{L_{t}^{2}L_{x}^{2}}\leq
2^{-\frac{j}{p}}C||P_{0}Q_{j}\phi||_{\dot{X}_{0}^{\frac{1}{2},\frac{1}{2},\infty}}
\end{equation}
for $p\leq q$ by Bernstein's inequality, which gives a uniform
exponential gain for $j\geq -100$ provided $p\leq q,\,K>p>2+\mu$,
and interpolating between the inequality thus gotten for
$L_{t}^{4}L_{x}^{4}$ and the inequality
$||P_{0}Q_{j}\phi||_{L_{t}^{\infty}L_{x}^{2}}\leq C
2^{\frac{j}{2}}||P_{0}Q_{j}\phi||_{L_{t}^{2}L_{x}^{2}}$. For
$P_{0}Q_{<-100}\phi$, we split the unit sphere into caps of size
$2^{-50}$, say, observe that
\\
\begin{equation}\begin{split}
||P_{0}Q_{<100}\phi||_{L_{t}^{p}L_{x}^{q}}^{2}=&||(P_{0}Q_{<-100}\phi)^{2}||_{L_{t}^{\frac{p}{2}}L_{x}^{\frac{q}{2}}}
=||(\sum_{\kappa\in
K_{-50}}P_{0,\kappa}Q_{<-100}\phi)^{2}||_{L_{t}^{\frac{p}{2}}L_{x}^{\frac{q}{2}}}\\
&\leq C (\sum_{\kappa\in
K_{-50}}||P_{0,\kappa}Q_{<-100}\phi||_{L_{t}^{p}L_{x}^{q}})^{2}\\
\end{split}\end{equation}
\\

Since we will be implementing a bootstrapping argument, we can
only assume the a priori existence of a solution on a finite time
interval $[-T,T]$. We therefore need to localize the above
(frequency-localized) norms to this interval. To wit
\begin{equation}
||P_{k}\phi ||_{S[k]([-T,T]\times\mathbf{R}^{3})}:=
\inf_{\psi\in\mathcal{S}(\mathbf{R}^{3+1}),\,\psi|_{[-T,T]} =
\phi}||P_{k}\psi||_{S[k](\mathbf{R}^{3+1})}
\end{equation}
\begin{equation}
||P_{k}\phi ||_{N[k]([-T,T]\times\mathbf{R}^{3})}:=
\inf_{\psi\in\mathcal{S}(\mathbf{R}^{3+1}),\,\psi|_{[-T,T]} =
\phi}||P_{k}\psi||_{N[k](\mathbf{R}^{3+1})}
\end{equation}
\\

We can now formulate the following energy inequality, which is the
essential link between the $N[k]$ and $S[k]$-norm that will allow
us to finish the bootstrapping argument:
\begin{equation}\label{energy}
||P_{k}\phi||_{S[k]([-T,T]\times\bf{R}^{3})}\leq C[||\Box
P_{k}\phi||_{N[k]([-T,T]\times\bf{R}^{3})}+||\phi[0]||_{\dot{H}^{\frac{1}{2}}\times\dot{H}^{-\frac{1}{2}}}]
\end{equation}
where $C$ is independent of $T$. This is proved as in \cite{Tao
2}; the only difference between our $S[k,\kappa]$ norm and Tao's
$S[k, \kappa]$-norm (other than the different scaling, which
doesn't affect the proof) is the addition of
$\sup_{\frac{1}{p}+\frac{1}{q}\leq
\frac{1}{2},p>2+\mu}2^{k(\frac{1}{p}+\frac{3}{q}-1)}||P_{k}\phi||_{L_{t}^{p}L_{x}^{q}}$.
However, this is majorized by
$||\phi||_{\dot{X}_{k}^{\frac{1}{2},\frac{1}{2},1}}$ and moreover
well behaved under multiplication with $L^{\infty}$ functions,
which is all that is required to prove the above inequality, viz.
Tao's proof.
\\

It is important that the $S[k]([-T,T]\times\mathbf{R}^{3})$-norms
of the frequency localized components of a Schwartz function are
in a sense uniformly lower semicontinuous with respect to $T$, as
demonstrated in \cite{Tao 2}. In particular, we may assume that
$T>0$ has been chosen such that the component functions $\phi$ of
our Wave Map satisfy
\begin{equation}
||P_{k}\phi||_{S[k]([-T,T]\times\mathbf{R}^{3})}\leq Cc_{k}
\end{equation}
where $c_{k}$ is a frequency envelope associated with the initial
conditions $\phi[0]\times\partial_{t}\phi[0]$ as above, i.e.
\begin{equation}
c_{k}:=
(\sum_{k'}2^{-\delta|k'-k|}(||P_{k'}\phi||_{\dot{H}^{\frac{1}{2}}}+||P_{k'}\partial_{t}\phi||_{\dot{H}^{-\frac{1}{2}}})^{2})^{\frac{1}{2}}
\end{equation}
Moreover, since we assume that $\phi$ is rapidly decaying in space
directions,  we can construct a Schwartz function $\tilde{\phi}$
with $\tilde{\phi}|_{[-T,T]}=\phi$ and such that
$||P_{k}\tilde{\phi} ||_{S[k]}\leq 2Cc_{k}$. This is achieved by
using a partition of unity. We will always substitute
$\tilde{\phi}$ for $\phi$ when making actual estimates.
\\

\begin{bf}Notation\end{bf}
The Riesz operators $R_{\nu}$, $\nu\in\{0,1,2,3\}$, refer to
operators $\partial_{\nu}{(\sqrt{-\triangle_{x}})}^{-1}$. We
usually omit the subscript for operators like
$\nabla_{x},\triangle_{x}$, understanding that they refer only to
space variables.\\
The symbol $\nabla^{-1}$ is either a shorthand for an operator
$\triangle^{-1}\partial_{i}$, or else refers to
$(\sqrt{-\triangle})^{-1}$, depending on the context.\\
We use the notation $P_{k+O(1)} = \sum_{k_{1}=k+O(1)}P_{k_{1}},
Q_{j+O(1)} = \sum_{j_{1}=j+O(1)}Q_{j_{1}}$. Also,
$||\phi||_{S[k+O(1)]} =
\sum_{k_{1}=k+O(1)}||P_{k_{1}}\phi||_{S[k_{1}]}$ etc.
\\
The following terminology, introduced by T.Tao in \cite{Tao 2},
shall be useful in the future: we call a Fourier multiplier
disposable if it is given by convolution with a translation
invariant measure of mass $\leq O(1)$. In particular, operators
such as $P_{k}, P_{k}Q_{<>j}$ where $j\geq k+O(1)$ are disposable,
see above reference. By contrast, $Q_{j}$ is not disposable.
However, it acts boundedly on Lebesgue spaces of the form
$L_{t}^{p}L_{x}^{2}$.\\
Whenever we consider an expression of the form $P_{0}(AB[CD])$,
for example, we shall refer to $A,B,C,D$ as inputs and the whole
expression as output. Also, when referring to $[,]$, we mean
$[CD]$, while $(,)$ would refer to $P_{0}(AB[CD])$; thus the shape
of brackets matters in the discussion. When considering a part of
the whole expression such as $[CD]$, we may also refer to this as
output, and $C,D$ as inputs, depending on the context.
\\

\newpage

\begin{bf}Two Important Inequalities\end{bf}
\\

Bernstein's inequality in the form $||P_{k}\phi||_{L_{x}^{p}}\leq
C 2^{\frac{3k}{2}-\frac{3k}{p}}||P_{k}\phi||_{L_{x}^{2}}$ or
variations thereof will be frequently used in the sequel.
\\

Moreover, the following improvement of Bernstein's inequality
which is a consequence of Strichartz' inequality shall be used
later on:
\begin{equation}\label{improved Bernstein}
||P_{k}Q_{j}\phi||_{L_{t}^{2}L_{x}^{\infty}}\leq
C2^{\frac{1}{4}\min\{0,j-k\}}2^{\frac{3k}{2}}||P_{k}Q_{j}\phi||_{L_{t}^{2}L_{x}^{2}}
\end{equation}
For this see \cite{Tao 2}. The intuition here is that Strichartz
estimates in the form
$2^{k(\frac{1}{p}+\frac{3}{q}-1)}||P_{k}Q_{j}\phi||_{L_{t}^{p}L_{x}^{q}}\leq
C||P_{k}Q_{j}\phi||_{\dot{X}_{k}^{\frac{1}{2},\frac{1}{2},1}}$ for
a Strichartz pair $(p,q), p\leq q$ can be interpreted as an
improvement over Bernstein's inequality; indeed, one gains
exponentially in  $j-k$, if this is small.
\\

\begin{bf}Summary of the key properties satisfied by these
spaces\end{bf}
\\

One difficulty in the present approach consists in showing that
the "renormalized $\phi$", i.e. $\psi: =
e^{-i\triangle^{-1}\sum_{k=1}^{3}\partial_{k}\tilde{\phi}^{1}_{k}}\tilde{\phi}$
is under approximately the same frequency envelope as $\phi$,
where $\tilde{\phi}, \tilde{\phi}^{1}_{k}$ are
(real-valued)\footnote{Note that the $S[k]$, $N[k]$ are
conjugation invariant. Thus we can always find real-valued
extensions of our component functions with the required
properties.} Schwartz functions agreeing with $\phi, \phi^{1}_{k}$
on $[-T,T]$ and for which the $S[k]$-norms of the frequency
localized pieces sit under approximately the same frequency
envelope, and later on translate the improved estimates for $\psi$
back to $\phi$. This step is trivial in $4$ dimensions and higher,
as one doesn't have to invoke the Fourier transform. However, it
becomes rather delicate in our present setting. For this, we have
the following theorem:
\\

\begin{theorem}\label{Gauge Change}
Let $\tilde{\phi}^{1}_{j}$ be as above, and
$\psi\in\mathcal{S}(\mathbf{R}^{3+1})$ satisfy
$||P_{k}\psi||_{S[k]}\leq \mu c_{k}$. Then we have
\\
\begin{equation}
||P_{k}(e^{i\sum_{j=1}^{3}\triangle^{-1}\partial_{j}\tilde{\phi}^{1}_{j}}\psi)||_{S[k]}\leq
C\mu c_{k}
\end{equation}
\end{theorem}

This theorem will follow immediately from the following natural
generalization:
\begin{proposition}
Let $f(x)$ be a smooth function with all derivatives up to and
including fourth order bounded. Then, under the same assumptions
as in the theorem,  we have
\begin{equation}
||P_{k}(f(\triangle^{-1}\sum_{j=1}^{3}\partial_{j}\tilde{\phi}^{1}_{j})\psi)||_{S[k]}\leq
C\mu c_{k}
\end{equation}
\end{proposition}
The proof of this Proposition, as well as the following technical
lemmata and propositions, shall be relegated to a later section.

In addition to the preceding "Gauge-Change"-estimate, we need to
utilize slightly modified versions of bilinear estimates proved in
\cite{Tao 2}, and most importantly the following version of Tao's
lemma(14.1):
\begin{lemma}\label{bilinear1}(T.Tao):
Let $j\leq\min(k_{1},k_{2})+O(1)$. Then
\begin{equation}
||P_{k}(F\psi)||_{N[k]}\leq
C2^{-\delta_{1}|k-\max(k_{1},k_{2})|}2^{-\delta_{2}|j-\min(k_{1},k_{2})|}||F||_{\dot{X}_{k_{1}}^{\frac{1}{2},-\frac{1}{2},\infty}}
||\psi||_{S[k_{2}]}
\end{equation}
for all Schwartz functions $F$ with Fourier support at frequency
$2^{k_{1}}$ and modulation $2^{j}$ while $\psi$ is at frequency
$2^{k_{2}}$, $\delta_{1},\delta_{2}>0$. Moreover, we also have
\begin{equation}
||\nabla P_{k}(F\psi)||_{N[k]}\leq
C2^{-\delta_{1}|k-\max(k_{1},k_{2})|}2^{-\delta_{2}|j-\min(k_{1},k_{2})|}||F||_{\dot{X}_{k_{1}}^{\frac{1}{2},-\frac{1}{2},\infty}}
||\nabla\psi||_{S[k_{2}]}
\end{equation}
\end{lemma}
The only difference here is that our spaces $S[k], N[k]$ are
scaled-down versions of Tao's spaces. In particular, for high-high
interactions, the above version as expressed in the first
inequality is slightly stronger than the original version.
Nevertheless, the proof follows
Tao's original proof almost identically, and only deviates in one small detail.

The intuition behind this rather technical estimate is that it
allows one for example to play the (small)size of the modulation
of $F$ against the frequency of $\psi$. This shall be important
for estimating the trilinear null-form \eqref{null-form}, in
particular in the delicate case when all inputs have relatively
low modulation but still too large to employ successfully the full
algebraic null-structure, see e.g. the proof of
lemma~\ref{low-high interactions} in the last section.

This lemma entails the following fundamental null-form estimate,
due to Tao, which we again modify for our scaled-down spaces
$S[k]$.
\begin{lemma}\label{bilinear2}(T.Tao)
Let $\phi$,$\psi$ be Schwarz functions on $\mathbf{R}^{3+1}$. Then
letting $k_{1}=k_{2}+O(1)\geq O(1)$, we have
\begin{equation}
||P_{0}[R_{\nu}P_{k_{1}}\phi\partial^{\nu}P_{k_{2}}\psi]||_{N[0]}\leq
C2^{-\delta
k_{1}}||P_{k_{1}}\phi||_{S[k_{1}]}||P_{k_{2}}\psi||_{S[k_{2}]}
\end{equation}
for some $\delta>0$.
\end{lemma}

In addition to the above "high-high" version, we will also need
the following "low-high" version, proved by T.Tao:
\begin{lemma}\label{bilinear3}
Let $k_{1}=O(1), k_{2}<O(1)$. Then
\begin{equation}
||P_{0}\nabla_{x}[R_{\nu}P_{k_{1}}\phi
R^{\nu}P_{k_{2}}\psi]||_{N[0]}\leq
C||P_{k_{1}}\phi||_{S[k_{1}]}||P_{k_{2}}\psi||_{S[k_{2}]}
\end{equation}
\end{lemma}

We will need the following fundamental bilinear inequality, which
is again a slight modification of an inequality of Tao:
\begin{lemma}\label{bilinear4}
Let $\phi,F$ be Schwartz functions, and $k_{1}=k_{2}+O(1)$. Then
we have
\begin{equation}
||P_{0}(P_{k_{1}}\phi P_{k_{2}}F)||_{N[0]}\leq C2^{-\delta k_{1}}
||P_{k_{1}}\phi||_{S[k_{1}]}||\nabla_{x}(P_{k_{2}}F)||_{N[k_{2}]}
\end{equation}
for some $\delta >0$.\\
Moreover, we have the estimate
\begin{equation}
||P_{0}\nabla_{x}(\phi P_{k_{2}}F)||_{N[0]}\leq
C(||\phi||_{L_{t}^{\infty}L_{x}^{\infty}}+\sup_{k}||P_{k}\nabla_{x}\phi||_{S[k]})||\nabla_{x}(P_{k_{2}}F)||_{N[k_{2}]}
\end{equation}
\end{lemma}

\begin{proof}
The proof of this is again almost identical to Tao's proof; in
particular, for the first inequality, one simply has to substitute
the above version of Tao's lemma(14.1) in the appropriate places.
\end{proof}

Finally, we have the important trilinear inequality
\begin{proposition}\label{Null-Form}
Let $\psi_{l},\, l=1,2,3$ be Schwartz functions on
$\mathbf{R}^{3+1}$. We then have the estimate
\begin{equation}\begin{split}
&||P_{0}(\sum_{j=1}^{3}\triangle^{-1}\partial_{j}[R_{\nu}P_{k_{1}}\psi_{1}R_{j}P_{k_{2}}\psi_{2}-R_{j}P_{k_{1}}\psi_{1}R_{\nu}P_{k_{2}}\psi_{2}]\partial^{\nu}P_{k_{3}}\psi_{3})||_{N[0]}\\&\leq
C 2^{-\delta_{1} |k_{1}-k_{2}|}2^{-\delta_{2}
|k_{3}|}\prod_{l=1}^{3}||P_{k_{l}}\psi_{l}||_{S[k_{l}]}\\
\end{split}\end{equation}
for appropriate constants $\delta_{1},\delta_{2}>0$.
\end{proposition}

\newpage
\section{Proof of theorem~\ref{theorem:Main}}

As explained in section 2, the proof of the Theorem can be reduced
to the proof of global regularity of the $\phi^{i}_{\alpha}$, or
$\phi$ for short. In order to formulate a 'bootstrapping
Proposition' for $\phi$, we first need to verify that smallness of
the initial conditions as stated in Theorem~\ref{theorem:Main}
implies a corresponding statement for the initial conditions of
$\phi$.

\begin{lemma}
Under the conditions of Theorem~\ref{theorem:Main}, we have
$||\phi[0]\times\partial_{t}\phi[0]||_{\dot{H}^{\frac{1}{2}}\times
\dot{H}^{-\frac{1}{2}}}< C \epsilon$.
\end{lemma}

\begin{proof}: Note that the wave equation satisfied by $(x,y)$ is

\begin{equation}
\Box x = \frac{2}{y}\partial_{\nu}x\partial^{\nu}y,\, \Box y =
\frac{-1}{y}(\partial_{\nu}x\partial^{\nu}x-\partial_{\nu}y\partial^{\nu}y)
\end{equation}

We need to show that
$||\partial_{t}(\frac{\partial_{t}x}{y})||_{H^{-\frac{1}{2}}}$ is
small, the other cases being similar or simpler. Now
\begin{equation}
\partial_{t}(\frac{\partial_{t}x}{y}) =
\sum_{i=1}^{3}\partial_{i}(\frac{\partial_{i}x}{y})-\sum_{i=1}^{3}\frac{\partial_{\nu}x\partial^{\nu}y}{y^{2}}
\end{equation}
Clearly, the first term has small $H^{-\frac{1}{2}}$-norm by our
assumptions. As to the 2nd, this is a simple consequence of
$||fg||_{\dot{H}^{-\frac{1}{2}}}\leq
C||f||_{\dot{H}^{\frac{1}{2}}}||g||_{\dot{H}^{\frac{1}{2}}}$.
\end{proof}

We now formulate the bootstrapping argument implying global
regularity of the $\phi$ for suitably small initial conditions:

\begin{proposition}\label{bootstrap}
Let the  smooth functions $\phi^{a}_{\alpha}$ satisfy the system
\eqref{1}-\eqref{4} on $[-T,T]\times\mathbf{R}^{3}$, with initial
conditions as in the preceding lemma. Assume that
\begin{equation}
||P_{k}\phi^{a}_{\alpha}||_{S[k]([-T,T]\times\mathbf{R}^{3})}\leq
M c_{k}
\end{equation}
for a frequency envelope $c_{k}$ covering the initial data . Then
we have
\begin{equation}
||P_{k}\phi^{a}_{\alpha}||_{S[k]([-T,T]\times\mathbf{R}^{3})}\leq
\frac{M}{2}c_{k}
\end{equation}
provided that $\epsilon$ is small enough, independently of $T$,
and $M\gg 1$.
\end{proposition}

\begin{bf}Remark\end{bf}: Global regularity of $\phi$ follows now
via the local result in \cite{Kl-M 3}. See also \cite{Kl-S 2}.

\newpage
\begin{proof}
Recall the quantities $\psi_{\alpha}:=
e^{-i\Phi}(\phi^{1}_{\alpha}+i\phi^{2}_{\alpha})$ introduced in
section 2, where $\Phi =
\triangle^{-1}\sum_{j=1}^{3}\partial_{j}\phi^{1}_{j}$. They
satisfy the equation
\begin{equation}
\Box\psi_{\alpha} = I+II+III \nonumber
\end{equation}
where
\begin{equation}
I = [-i\Box\Phi-\partial_{\nu}\Phi\partial^{\nu}\Phi]\psi_{\alpha}
\end{equation}
\begin{equation}
II =
2i(\triangle^{-1}\sum_{j=1}^{3}\partial_{j}[\psi^{1}_{j}\psi^{2}_{\nu}-\psi^{2}_{j}\psi^{1}_{\nu}]\partial^{\nu}\psi_{\alpha}+i\sum_{j=1}^{3}\triangle^{-1}\partial_{j}[\psi^{1}_{j}\psi^{2}_{\nu}-\psi^{2}_{j}\psi^{1}_{\nu}]\partial^{\nu}\Phi\psi_{\alpha})
\end{equation}
\begin{equation}
III =
e^{-i\Phi}[\phi^{1}_{\alpha}(\phi^{1}_{\nu}\phi^{1\nu}+\phi^{2}_{\nu}\phi^{2\nu})+2i\phi^{1\nu}(\phi^{1}_{\nu}\phi^{2}_{\alpha}-\phi^{1}_{\alpha}\phi^{2}_{\nu})]
\end{equation}
\\

We now consider each of these in turn, starting with the most
difficult and pivotal term II. Recall from last section that we
have to localize to frequency $\sim 2^{0}$ and evaluate the
$N[0]$-norm. The atomic nature of this space implies that we can
estimate in any of the norms $L_{t}^{1}L_{x}^{2}$,
$\dot{X}_{0}^{-\frac{1}{2}, -\frac{1}{2},1}$, as well as the
complicated norm involving null-frame spaces, i.e. the third class
of atoms
defining $N[0]$. \\
In order to make estimates, we have to substitute Schwartz
functions for the inputs of the nonlinearity coinciding with them
on the time interval $[-T,T]$, as well as having frequency
localized components whose $S[k]$-norms are bounded by a multiple
of the frequency envelope $Mc_{k}$. In the sequel, we shall not
make a distinction between these extensions and the local Wave
Map, implicitly carrying out the above substitution immediately
before estimating and after having reformulated the individual
ingredients of the nonlinearity to become amenable to estimates.

\subsection{Treating II: First term of II}

Recall form section 2 that using dynamic separation, we have
reduced this to estimating \eqref{easy1}, \eqref{easy2} as well as
the difficult \eqref{null-form}. \\
Note that on account of the $L_{t}^{4}L_{x}^{4}$-boundedness of
$\nabla^{-1}(\nabla^{-1}(\psi^{2})\psi)$ (which of course follows
from Sobolev's inequality as well as passing from frequency
localized pieces to the full function via the fundamental theorem
of Littlewood-Paley theory), we can immediately reduce control of
\eqref{easy2} to estimating \eqref{easy1}.\\
As to the latter, we distinguish between low-high, high-high and
high-low interactions, i.e. reducing $\nabla_{x,t}\psi$ to medium,
high and low frequency:
\\

Low-High interactions: we use here the fact that $L_{x}^{4} =
L_{x}^{4,4}$, the latter a Lorentz space. Also, we utilize the
improved Hardy-Littlewood-Sobolev inequality for these spaces:

\begin{equation}\begin{split}
&||P_{0}[P_{<O(1)}\nabla^{-1}(\nabla^{-1}(\nabla^{-1}(\psi^{2})\psi)R_{\nu}\Psi^{a})P_{O(1)}\nabla_{x,t}\psi]||_{L_{t}^{1}L_{x}^{2}}\\
&\leq
C||P_{<O(1)}\nabla^{-1}(\nabla^{-1}(\nabla^{-1}(\psi^{2})\psi)R_{\nu}\Psi^{a})||_{L_{t}^{1}L_{x}^{3,1}}||P_{O(1)}\nabla_{x,t}\psi||_{L_{t}^{\infty}L_{x}^{2}}\\
&\leq
C||\nabla^{-1}(\nabla^{-1}(\psi^{2})\psi)||_{L_{t}^{\frac{4}{3}}L_{x}^{12,\frac{4}{3}}}||R_{\nu}\Psi^{a}||_{L_{t}^{4}L_{x}^{4}}||P_{O(1)}\nabla_{x,t}\psi||_{L_{t}^{\infty}L_{x}^{2}}\\
&\leq
C||\nabla^{-1}(\psi^{2})||_{L_{t}^{2}L_{x}^{6,2}}||\psi||_{L_{t}^{4}L_{x}^{4}}||R_{\nu}\Psi^{a}||_{L_{t}^{4}L_{x}^{4}}||P_{O(1)}\nabla_{x,t}\psi||_{L_{t}^{\infty}L_{x}^{2}}\\
&\leq C||\psi||_{L_{t}^{4}L_{x}^{4}}^{3}||R_{\nu}\Psi^{a}||_{L_{t}^{4}L_{x}^{4}}||P_{O(1)}\nabla_{x,t}\psi||_{L_{t}^{\infty}L_{x}^{2}}\\
&\leq CM^{5}c_{0}\\
\end{split}\end{equation}

High-High interactions: This is more elementary. Note that we
utilize the exponential decay provided by
$||P_{k}\phi||_{L_{x}^{2}}\leq
C2^{-\frac{k}{2}}||P_{k}\phi||_{H^{\frac{1}{2}}}$, provided
$k>O(1)$, as well as Bernstein's inequality to majorize
$L_{t}^{1}L_{x}^{2}$ of the output by  a multiple of
$L_{t}^{1}L_{x}^{\frac{6}{5}}$ of the output, which lives at
frequency $\sim 1$:

\begin{equation}\begin{split}
&||\sum_{k\geq O(1)}P_{0}[\nabla^{-1}P_{k+O(1)}(\nabla^{-1}(\nabla^{-1}(\psi^{2})\psi)R_{\nu}\Psi^{a})\nabla_{x,t}P_{k}\psi||_{L_{t}^{1}L_{x}^{2}}\\
&\leq C\sum_{k\geq
O(1)}2^{-k}||P_{k+O(1)}(\nabla^{-1}(\nabla^{-1}(\psi^{2})\psi)R_{\nu}\Psi^{a})||_{L_{t}^{1}L_{x}^{3}}||\nabla_{x,t}\psi||_{L_{t}^{\infty}L_{x}^{2}}\\
&\leq C\sum_{k\geq O(1)}M^{5}2^{-\frac{k}{2}}c_{k}\leq
CM^{5}c_{0}\\
\end{split}\end{equation}

High-Low interactions: this is even more elementary as an extra
derivative ($\nabla_{x,t}$) falls on a low frequency function,
while the first $\nabla^{-1}$ is harmless. We therefore skip this
case.
\\

Now we estimate the null-form, which is straightforward in light
of Proposition~\ref{Null-Form}: the exponential gains in the
statement of this Proposition allow one to counteract the
'exponential loss' inherent in the frequency envelope, provided
its exponent $\sigma$ is chosen small enough:

\begin{equation}\begin{split}
&||P_{0}[\sum_{j=1}^{3}\triangle^{-1}\partial_{j}(R_{\nu}\Psi^{1}R_{j}\Psi^{2}-R_{j}\Psi^{1}R_{\nu}\Psi^{2})\partial^{\nu}\psi_{\alpha}]||_{N[0]}\\
&\leq \sum_{k_{1},k_{2},k_{3}}||P_{0}[\sum_{j=1}^{3}\triangle^{-1}\partial_{j}(R_{\nu}P_{k_{1}}\Psi^{1}R_{j}P_{k_{2}}\Psi^{2}-R_{j}P_{k_{1}}\Psi^{1}R_{\nu}P_{k_{2}}\Psi^{2})\partial^{\nu}P_{k_{3}}\psi_{\alpha}]||_{N[0]}\\
&\leq
\sum_{k_{1},k_{2},k_{3}}C2^{-\delta_{1}|k_{1}-k_{2}|}2^{-\delta_{2}|k_{3}|}M^{3}\prod_{a=1}^{3}c_{k_{a}}\\
&\leq CM^{3}\sum_{k}c_{k}^{2}c_{0}\leq CM^{3}\epsilon c_{0}\\
\end{split}\end{equation}

provided we choose the $\sigma$ in the definition of the frequency
envelope (viz. section 3) much smaller than
$\min\{\delta_{1},\delta_{2}\}$.
\\

\begin{bf}Second term of II\end{bf}
\\

Applying dynamic separation  to $[,]$ as explained in section 2.4,
we decompose this term into a term with null-structure, as well as
error terms of high linearity. Of course one is tempted here to
simply use the $L_{t}^{4}L_{x}^{4}$-Strichartz estimate; however,
in a sense our strategy of controlling each frequency localized
piece is backfiring here, as we 'have to recover the initial
frequency envelope' in our estimates, which requires in particular
obtaining exponential gains for high-frequency interactions. The
$L_{t}^{4}L_{x}^{4}$-estimate is too weak for this, and we have to
invoke the hidden null-structure. It goes without saying that
being able to prove a genuine trilinear null-form estimate on the
level of Lebesgue spaces for \eqref{null-form} would render these
convoluted considerations unnecessary, although the overall
complexity of the argument would barely improve.
\\

Hence we reduce the 2nd term of II to the following:
\\

Null-form:
\begin{equation}\label{quadr. null-form}
\sum_{j=1}^{3}(\triangle^{-1}\partial_{j}[R_{j}\Psi^{1}R_{\nu}\Psi^{2}-R_{j}\Psi^{2}R_{\nu}\Psi^{2}]\partial^{\nu}\Phi)\psi_{\alpha}
\end{equation}

Error terms: these are of the rough form

\begin{equation}
\nabla^{-1}(R_{\mu}\Psi^{a}\nabla^{-1}(\nabla^{-1}(\psi^{2})\psi))\partial^{\nu}\Phi\psi
\end{equation}

\begin{equation}
\nabla^{-1}(\nabla^{-1}(\nabla^{-1}(\psi^{2})\psi)\nabla^{-1}(\nabla^{-1}(\psi^{2})\psi))\partial^{\nu}\Phi\psi,
\end{equation}

where $\Psi^{a},\,a=1,2$ is as in section 2.4.\\
These error terms can be handled as before, utilizing e.g. Lorentz
spaces for low-high frequency interactions. As this doesn't offer
anything new, it is left out.\\

As to the null-form \eqref{quadr. null-form}, we decompose this
into various frequency interactions of $\psi_{\alpha}$ and $(,)$.
Most cases are elementary, with the exception of the high-high
case. However, we can then profit from
Proposition~\ref{Null-Form}.
\\

Estimation of the null-form \eqref{quadr. null-form}:
\\

1): Both $(,)$ and $\psi_{\alpha}$ at frequency $<O(1)$:

\begin{equation}\begin{split}
&||P_{0}(P_{<O(1)}(\triangle^{-1}\partial_{j}[R_{j}\Psi^{1}R_{\nu}\Psi^{2}-R_{\nu}\Psi^{1}R_{j}\Psi^{2}]\partial^{\nu}\Phi)P_{<O(1)}\psi_{\alpha})||_{L_{t}^{1}L_{x}^{2}}\\
&\leq
C||\nabla^{-1}[R_{j}\Psi^{1}R_{\nu}\Psi^{2}-R_{\nu}\Psi^{1}R_{j}\Psi^{2}]||_{L_{t}^{2}L_{x}^{6}}||\partial^{\nu}\Phi||_{L_{t}^{4}L_{x}^{4}}||P_{<O(1)}\psi_{\alpha}||_{L_{t}^{4}L_{x}^{12}}\\
&\leq CM^{4}\epsilon c_{0},\\
\end{split}\end{equation}

where we have used the fact that
$||P_{<O(1)}\psi_{\alpha}||_{L_{t}^{4}L_{x}^{12}}\leq Cc_{0}$, as
is easily verified.
\\

2.1): High-High interactions of $(,)$ and $\psi_{\alpha}$:
$\psi_{\alpha}$ at frequency $2^{k_{2}}$, $k_{2}>>1$, $\Phi$ at
frequency $<2^{k_{2}-100}$. This implies that $[,]$ is at
frequency $2^{k_{2}+O(1)}$. Hence we can estimate this by

\begin{equation}\begin{split}
&||P_{0}(\triangle^{-1}\partial_{j}P_{k_{2}+O(1)}[R_{j}\Psi^{1}R_{\nu}\Psi^{2}-R_{\nu}\Psi^{1}R_{j}\Psi^{2}]P_{<k_{2}-100}\partial^{\nu}\Phi
P_{k_{2}}\psi_{\alpha})||_{L_{t}^{1}L_{x}^{2}}\\
&\leq
C2^{-k_{2}}||\nabla_{x,t}\nabla^{-1}\Psi^{1}||_{L_{t}^{4}L_{x}^{4}}||\nabla_{x,t}\nabla^{-1}\Psi^{2}||_{L_{t}^{4}L_{x}^{4}}||\partial^{\nu}\Phi||_{L_{t}^{4}L_{x}^{4}}||P_{k_{2}}\psi_{\alpha}||_{L_{t}^{4}L_{x}^{4}}\\
&\leq CM^{4}\epsilon 2^{-k_{2}}c_{k_{2}}\\
\end{split}\end{equation}

This can be summed over $k_{2}>>1$, to give a bound of the form
$CM^{4}\epsilon c_{0}$.
\\

2.2): High-High interactions of $(,)$ and $\psi_{\alpha}$:
$\psi_{\alpha}$ at frequency $2^{k_{2}},\,k_{2}>>1$, $\Phi$ at
frequency $\geq 2^{k_{2}-100}$: As mentioned before, this case
does not appear covered by the Strichartz estimates, on account of
the possibility that all inputs could have large frequencies while
there is 'no exponential gain for free', as before. The easiest
way around this difficulty is to refer to the already proved
estimate for the null-form \eqref{null-form}, in addition to the
estimate lemma~\ref{bilinear4}. This gives the estimate

\begin{equation}\begin{split}
&||P_{0}(\sum_{j=1}^{3}\triangle^{-1}\partial_{j}[R_{j}\Psi^{1}R_{\nu}\Psi^{2}-R_{\nu}\Psi^{1}R_{j}\Psi^{2}]P_{\geq
k_{2}-100}\partial^{\nu}\Phi P_{k_{2}}\psi_{\alpha})||_{N[0]}\\
&\leq C2^{-\delta
k_{2}}||\nabla(\sum_{j=1}^{3}\triangle^{-1}\partial_{j}[R_{j}\Psi^{1}R_{\nu}\Psi^{2}-R_{\nu}\Psi^{1}R_{j}\Psi^{2}]P_{\geq
k_{2}-100}\partial^{\nu}\Phi P_{k_{2}})||_{N[k_{2}+O(1)]}\\
&||P_{k_{2}}\psi_{\alpha}||_{S[k_{2}]}\\
&\leq CM^{4}\epsilon 2^{-\delta k_{2}}c_{k_{2}}\\
\end{split}\end{equation}

This can be summed over $k_{2}>>1$ to yield a bound of the
required form, provided the $\sigma$ in the definition of the
frequency envelope is much smaller than $\delta_{2}$. We are done
with the estimation of II.
\\

\subsection{Treating I: First term of I}

Recall the definition of $\Phi$: $\Phi =
\sum_{k=1}^{3}\triangle^{-1}\partial_{k}\phi^{1}_{k}$. Thus upon
using the wave equation satisfied by $\phi$, the first term of I
is seen to reduce to

\begin{equation}
\sum_{j=1}^{3}\triangle^{-1}\partial_{j}(-\phi^{1}_{\nu}\partial^{\nu}\phi^{2}_{j}+"\phi^{3}")\psi_{\alpha}
\end{equation}

 In spite of the trilinear character of this expression, it is not yet amenable to good estimates, as can be seen when
 for example $\partial^{\nu}\phi^{2}_{j}$ is at very large frequency. We therefore employ the div-curl system and
 corresponding dynamic separation for the $\phi^{1}_{\alpha}$ in
 order to introduce a null-structure plus error terms. There
 results the expression
 
 \begin{equation}
 \sum_{j,l}^{3}\triangle^{-1}\partial_{j}[(R_{\nu}R_{l}\phi^{1}_{l}+"\nabla^{-1}(\phi^{2})")\partial^{\nu}\phi^{2}_{j}+"\phi^{3}"]\psi_{\alpha}
 \end{equation}

 Here "$\nabla^{-1}(\phi^{2})$" refers to
 $\sum_{r=1}^{3}\triangle^{-1}\partial_{r}(\phi^{1}_{r}\phi^{2}_{\nu}-\phi^{1}_{\nu}\phi^{2}_{r})$.
We estimate each of these terms in turn:
\\

\begin{bf}A\end{bf}: $\triangle^{-1}\partial_{j}(\phi^{3})\psi_{\alpha}$: this
term is entirely elementary, for the simple reason that high-high
interactions of $(,)$ and $\psi_{\alpha}$ imply an automatic
exponential gain in the frequency on account of the
$\triangle^{-1}\partial_{j}$. For completeness' sake, we provide
the simple details:
\\

A.1): $\psi_{\alpha}$ at frequency $O(1)$, $(,)$ at frequency
$<O(1)$:

\begin{equation}\begin{split}
&||P_{0}[\sum_{k<O(1)}\triangle^{-1}\partial_{j}P_{k}(\phi^{3})P_{O(1)}\psi_{\alpha}]||_{L_{t}^{1}L_{x}^{2}}\\
&\leq
C\sum_{k<O(1)}C2^{\frac{k}{2}}||\phi^{3}||_{L_{t}^{\frac{4}{3}}L_{x}^{\frac{4}{3}}}||P_{O(1)}\psi_{\alpha}||_{L_{t}^{4}L_{x}^{4}}\\
&\leq CM^{4}\epsilon c_{0},\\
\end{split}\end{equation}

where we have employed Bernstein's inequality estimating
$L_{x}^{4}$ in terms of $L_{x}^{\frac{4}{3}}$ for the term
$\phi^{3}$.
\\

A.2): $\psi_{\alpha}$ at frequency $<O(1)$, $(,)$ at frequency
$O(1)$:

\begin{equation}\begin{split}
&||P_{0}[P_{O(1)}\triangle^{-1}\partial_{j}(\phi^{3})\sum_{k<O(1)}P_{k}\psi_{\alpha}]||_{L_{t}^{1}L_{x}^{2}}\\
&\leq
C\sum_{k<O(1)}||\phi^{3}||_{L_{t}^{\frac{4}{3}}L_{x}^{\frac{4}{3}}}||P_{k}\psi_{\alpha}||_{L_{t}^{4}L_{x}^{\infty}}\\
&\leq \sum_{k<O(1)}CM^{4}2^{\frac{3k}{4}}c_{k}\leq CM^{4}\epsilon
c_{0},\\
\end{split}\end{equation}

where we have applied Bernstein's inequality estimating
$L_{x}^{2}$ in terms of $L_{x}^{\frac{4}{3}}$ to the output.
\\

A.3): $\psi_{\alpha}$ at frequency $2^{k}$, $k>>1$, $(,)$ at
frequency $2^{k+O(1)}$: we have

\begin{equation}\begin{split}
&\sum_{k>>1}||P_{0}[\triangle^{-1}\partial_{j}P_{k+O(1)}(\phi^{3})P_{k}\psi_{\alpha}]||_{L_{t}^{1}L_{x}^{2}}\\
&\leq
\sum_{k>>1}C2^{-k}||\phi^{3}||_{L_{t}^{\frac{4}{3}}L_{x}^{\frac{4}{3}}}||P_{k}\psi_{\alpha}||_{L_{t}^{4}L_{x}^{4}}\\
&\leq C\sum_{k>>1}M^{4}2^{-k}\epsilon c_{k}\leq CM^{4}\epsilon
c_{0}\\
\end{split}\end{equation}

This then disposes of this term.
\\

\begin{bf}B\end{bf}:
$\triangle^{-1}\partial_{j}(R_{\nu}R_{l}\phi^{1}_{l}\partial^{\nu}\phi^{2}_{j})\psi_{\alpha}$:
this term of course has an obvious null-structure inherent in it.
We will split the term into various frequency interactions between
$(,)$ and $\psi_{\alpha}$:
\\

B.1): Low -high interactions: $(,)$ at frequency $<O(1)$,
$\psi_{\alpha}$ at frequency $O(1)$: we can't simply put all
inputs into $L_{t}^{3}L_{x}^{6}$, as the inputs of $(,)$ could
have very large frequency. In this bad case, we will revert to
lemma~\ref{bilinear2}, estimating $(,)$ in $N[k]$ and  gaining
exponentially in the difference between the logarithm of the
frequency of $(,)$ and its (large frequency) inputs. Then we can
apply lemma~\ref{bilinear4} to conclude. Of course, when the
inputs of $(,)$ have frequency at most comparable to the frequency
of $(,)$, we put all inputs into $L_{t}^{3}L_{x}^{6}$. Writing the
details:

\begin{equation}\begin{split}
&||P_{0}(\sum_{k<O(1)}\triangle^{-1}\partial_{j}P_{k}(R_{\nu}R_{l}\phi^{1}_{l}\partial^{\nu}\phi^{2}_{j})P_{O(1)}\psi_{\alpha})||_{N[0]}\\
&\leq \sum_{k<O(1)}\sum_{a\geq
k+O(1)}C||P_{k}\triangle^{-1}\nabla\partial_{j}(P_{a}R_{\nu}R_{l}\phi^{1}_{l}\partial^{\nu}P_{a+O(1)}\phi^{2}_{j})||_{N[k]}\\&||P_{O(1)}\psi_{\alpha}||_{S[O(1)]}\\
&+\sum_{k<O(1)}C||P_{k}\triangle^{-1}\nabla\partial_{j}(R_{\nu}R_{l}P_{<k+O(1)}\phi^{1}_{l}\partial^{\nu}P_{<k+O(1)}\phi^{2}_{j})P_{O(1)}\psi_{\alpha}||_{L_{t}^{1}L_{x}^{2}}\\
&\leq C\sum_{k<O(1)}\sum_{a\geq k+O(1)}2^{\delta(k-a)}||P_{a}\phi^{1}_{l}||_{S[a]}||P_{a+O(1)}\phi^{2}_{j}||_{S[a+O(1)]}||P_{O(1)}\psi_{\alpha}||_{S[O(1)]}\\
&+C\sum_{k<O(1)}\sum_{l_{1},l_{2}<k+O(1)}2^{\frac{l_{1}+l_{2}}{6}}||P_{l_{1}}\phi^{1}_{l}||_{S[l_{1}]}||P_{l_{2}}\phi^{2}_{j}||_{S[l_{2}]}||P_{O(1)}\psi_{\alpha}||_{S[O(1)]}\\
&\leq CM^{3}\epsilon c_{0}\\
\end{split}\end{equation}

where the $\delta$ in the preceding comes from
lemma~\ref{bilinear2}.
\\

B.2): High-low interactions: these are dealt with in identical
fashion, using lemma~\ref{bilinear2}, and lemma~\ref{bilinear4}
for the high-low interaction case(note that we have
$||P_{<O(1)}\psi_{\alpha}||_{L_{t}^{\infty}L_{x}^{\infty}}\leq
Cc_{0}$), hence left out.
\\

B.3): High-high interactions: $(,)$ at frequency $2^{k_{1}}$,
$k_{1}>>1$, $\psi_{\alpha}$ at frequency $2^{k_{2}}$ with
$k_{2}=k_{1}+O(1)$: we distinguish between the case when both
inputs of $(,)$ have frequency at least comparable to $2^{k_{1}}$,
and its complement,i.e. one input of $(,)$ at much lower frequency
than $(,)$. The former case is dealt with as in the preceding
cases, using lemmata~\ref{bilinear2}, ~\ref{bilinear4}, hence not
further elaborated. The treatment of the latter case, i.e. both
inputs of $(,)$ have at most frequency comparable with
$2^{k_{1}}$, is subsumed under the following lemma:

\begin{lemma}\label{almost endpoint}: Let $1<< k_{3}$,
$k_{2}\leq k_{1}-10\leq k_{3}$. Also, let $\psi_{j},\,j=1,2,3$ be
Schwartz functions. Then the  following estimates hold:

\begin{equation}
||P_{0}[\nabla^{-1}(R_{\nu}P_{k_{2}}\psi_{2}P_{k_{1}}\partial^{\nu}\psi_{1})P_{k_{3}}\psi_{3}]||_{N[0]}\leq
C2^{-\delta_{1}k_{3}}2^{\delta_{2}(k_{2}-k_{1})}\prod_{a=1}^{3}||P_{k_{a}}\psi_{a}||_{S[k_{a}]}
\end{equation}
\\
\begin{equation}
||P_{0}[R_{\nu}P_{k_{2}}\psi_{2}\nabla^{-1}P_{k_{1}}\psi_{1}\partial^{\nu}P_{k_{3}}\psi_{3}]||_{N[0]}\leq
C2^{-\delta_{1}k_{3}}2^{\delta_{2}(k_{2}-k_{1})}\prod_{a=1}^{3}||P_{k_{a}}\psi_{a}||_{S[k_{a}]}
\end{equation}

In the first inequality, $\nabla^{-1}$ refers to any operator of
the form $\triangle^{-1}\partial_{j}$, $j=1,2,3$.  Also, if
$k_{1}>>1$, $k_{2}\leq k_{1}$, then

\begin{equation}
||P_{0}[R_{\nu}P_{k_{2}}\psi_{2}P_{k_{1}}R^{\nu}\psi_{1}\psi_{3}]||_{N[0]}\leq
C2^{-\delta_{1}k_{1}}2^{\delta_{2}\min\{0,k_{2}\}}\prod_{a=1}^{2}||P_{k_{a}}\psi_{a}||_{S[k_{a}]}\sup_{k}\{||P_{k}\psi_{3}||_{S[k]}\}
\end{equation}

\end{lemma}

\begin{proof}: The proofs of these three assertions are
essentially identical, hence we will only deal with the first of
them: the idea of the proof consists in splitting into the case of
moderately small $k_{2}$, namely $k_{2}\geq -k_{1}$, and very
small $k_{2}$, or $k_{2}<-k_{1}$. Of course, $k_{1}=k_{3}+O(1)$,
hence is large. In the first case, we utilize
lemma~\ref{bilinear3} and lemma~\ref{bilinear4}, the latter
furnishing an exponential decay accounting for both factors on the
right-hand side of the inequality.\\
For the second case, we place the low-frequency term into a
Strichartz space slightly worse than $L_{t}^{2}L_{x}^{\infty}$,
while the high-frequency terms are placed into small perturbations
of the neutral(as far as exponential gains or losses are
concerned) $L_{t}^{4}L_{x}^{4}$. The assumption on the smallness
of the low frequency allows to play this successfully against the
high frequencies, resulting in the desired exponential gain.
\\

Case 1): $k_{1}-10\geq k_{2}\geq -k_{1}$: use
lemma~\ref{bilinear2}, lemma~\ref{bilinear4} to conclude that

\begin{equation}\begin{split}
&||P_{0}[\nabla^{-1}(R_{\nu}P_{k_{2}}\psi_{2}P_{k_{1}}\partial^{\nu}\psi_{1})P_{k_{3}}\psi_{3}]||_{N[0]}\\&\leq
C2^{-\delta k_{3}}
||(R_{\nu}P_{k_{2}}\psi_{2}P_{k_{1}}\partial^{\nu}\psi_{1})||_{N[k_{3}+O(1)
]}||P_{k_{3}}\psi_{3}||_{S[k_{3}]}\\
&\leq
C2^{-\frac{\delta}{3}k_{3}}2^{\frac{\delta}{3}(k_{2}-k_{1})}\prod_{a=1}^{3}||P_{k_{a}}\psi_{a}||_{S[k_{a}]}\\
\end{split}\end{equation}

Case 2): $k_{2}<-k_{1}$:

\begin{equation}\begin{split}
&||P_{0}[\nabla^{-1}(\partial^{\nu}P_{k_{1}}\psi_{1}P_{k_{2}}R_{\nu}\psi_{2})P_{k_{3}}\psi_{3}||_{N[0]}\\
&\leq
C2^{-k_{3}}||P_{k_{1}}\partial^{\nu}\psi_{1}||_{L_{t}^{4-}L_{x}^{4+}}||P_{k_{2}}R_{\nu}\psi_{2}||_{L_{t}^{2+}L_{x}^{M}}
||P_{k_{3}}\psi||_{L_{t}^{4-}L_{x}^{4+}}\\
\end{split}\end{equation}

where we set $\frac{2}{4+}+\frac{1}{M}=1$,
$\frac{2}{4-}+\frac{1}{2+} = 1$, $\frac{1}{2+}+\frac{1}{M} =
\frac{1}{2}$, $\frac{1}{4}-\frac{1}{4+} = \lambda >0$. Using this
we can estimate the above by

\begin{equation}\begin{split}
&C2^{4\lambda k_{1}}2^{(\frac{1}{2}-4\lambda)k_{2}}\prod_{a=1}^{3}||P_{k_{a}}\psi_{a}||_{S[k_{a}]}\\
&\leq C
2^{(\frac{1}{8}-\lambda)(k_{2}-k_{1})}2^{(4\lambda-(\frac{1}{4}-2\lambda))k_{1}}\prod_{a=1}^{3}||P_{k_{a}}\psi_{a}||_{S[k_{a}]}\\
\end{split}\end{equation}

We can now choose $\lambda = \frac{1}{100}$ in order to obtain the
required estimate. Note that the Strichartz estimates invoked here
are still controlled by our choice of $S[k]$.
\end{proof}

\begin{bf}C\end{bf}:
$P_{0}(\sum_{j=1}^{3}\triangle^{-1}\partial_{j}[\sum_{r=1}^{3}\triangle^{-1}\partial_{r}(\phi^{1}_{r}\phi^{2}_{\nu}-\phi^{1}_{\nu}\phi^{2}_{r})\partial^{\nu}\phi^{2}_{j}]\psi_{\alpha})$
This is very similar to the 2nd term of II: we have simply traded
an extra $\nabla^{-1}$ for an extra $\nabla$ (in
$\partial^{\nu}\phi^{2}_{j}$), recalling that $\partial^{\nu}\Phi$
is morally equivalent (as far as scaling is concerned) to $\phi$.
As before, the quadrilinear structure appears insufficient to
render this term amenable to a simple energy estimate. The
difficulty arises as usual for certain high-high interactions. By
now however, enough estimates are available to dispose of this
case rather quickly, invoking in particular the estimates for II:
\\

C.1): Low-high interactions: $[,]$ at frequency $<O(1)$,
$\psi_{\alpha}$ at frequency $O(1)$:

\begin{equation}\begin{split}
&||P_{0}(\sum_{j=1}^{3}\triangle^{-1}\partial_{j}\sum_{k<O(1)}P_{k}[\sum_{r=1}^{3}\triangle^{-1}\partial_{r}(\phi^{1}_{r}\phi^{2}_{\nu}-\phi^{1}_{\nu}\phi^{2}_{r})\partial^{\nu}\phi^{2}_{j}]P_{O(1)}\psi_{\alpha})||_{L_{t}^{1}L_{x}^{2}}\\
&\leq
C\sum_{j=1}^{3}\sum_{k<O(1)}2^{-k}||P_{k}(P_{<k+O(1)}"\nabla^{-1}(\phi^{2})"P_{<k+O(1)}\partial^{\nu}\phi^{2}_{j})||_{L_{t}^{\frac{4}{3}}L_{x}^{4}}||P_{O(1)}\psi_{\alpha}||_{L_{t}^{4}L_{x}^{4}}\\
&+C\sum_{j=1}^{3}\sum_{k<O(1)}\sum_{a>k+O(1)}2^{-k}||P_{k}(P_{a}"\nabla^{-1}(\phi^{2})"P_{a+O(1)}\partial^{\nu}\phi^{2}_{j})||_{L_{t}^{\frac{4}{3}}L_{x}^{4}}||P_{O(1)}\psi_{\alpha}||_{L_{t}^{4}L_{x}^{4}}\\
&\leq
C\sum_{j=1}^{3}\sum_{k<O(1)}2^{-k}2^{\frac{k}{2}}||P_{<k+O(1)}"\nabla^{-1}(\phi^{2})"||_{L_{t}^{2}L_{x}^{6}}||P_{<k+O(1)}\partial^{\nu}\phi^{2}_{j}||_{L_{t}^{4}L_{x}^{4}}\\&||P_{O(1)}\psi_{\alpha}||_{L_{t}^{4}L_{x}^{4}}\\
&+C\sum_{j=1}^{3}\sum_{k<O(1)}C2^{\frac{k}{2}}2^{-a}||P_{\geq
a+O(1)}\phi||_{L_{t}^{6}L_{x}^{3}}||P_{a+O(1)}\partial^{\nu}\phi^{2}_{j}||_{L_{t}^{3}L_{x}^{6}}||\phi||_{L_{t}^{4}L_{x}^{4}}\\&||P_{O(1)}\psi_{\alpha}||_{L_{t}^{4}L_{x}^{4}}\\
&\leq \sum_{k<O(1)}CM^{4}\epsilon 2^{\frac{k}{2}}c_{0}+
C\sum_{k<O(1)}\sum_{a\geq k+O(1)}2^{\frac{k}{2}}c_{a}^{2}c_{0}\\
&\leq CM^{4}\epsilon c_{0}\\
\end{split}\nonumber\end{equation}

We have repeatedly used Bernstein's inequality, as well as
Sobolev's inequality and the fact that $P_{a}(\phi^{2})=
P_{a}(P_{\geq a-10}\phi\phi)+P_{a}(P_{<a-10}\phi P_{\geq
a-10}\phi)$.
\\

C.2): High-high interactions: $[,]$ and $\psi_{\alpha}$ at high
frequencies. Here we utilize the estimates proved for II:

\begin{equation}\begin{split}
&||P_{0}(\sum_{j=1}^{3}\triangle^{-1}\partial_{j}P_{a}[\sum_{r=1}^{3}\triangle^{-1}\partial_{r}(\phi^{1}_{r}\phi^{2}_{\nu}-\phi^{1}_{\nu}\phi^{2}_{r})\partial^{\nu}\phi^{2}_{j}]P_{a+O(1)}\psi_{\alpha})||_{N[0]}\\
&\leq \sum_{a\geq O(1)}C2^{-\delta
a}||\nabla\sum_{j=1}^{3}\triangle^{-1}\partial_{j}P_{a}[\sum_{r=1}^{3}\triangle^{-1}\partial_{r}(\phi^{1}_{r}\phi^{2}_{\nu}-\phi^{1}_{\nu}\phi^{2}_{r})\partial^{\nu}\phi^{2}_{j}]||_{N[a]}\\
&||P_{a+O(1)}\psi_{\alpha}||_{S[a+O(1)]}\\
&\leq\sum_{a\geq O(1)}CM^{4}2^{-\delta a}c_{a}^{2}\leq
CM^{4}\epsilon
c_{0}\\
\end{split}\end{equation}

We have used lemma~\ref{bilinear4} for the second step.
\\

C.3): High-low interactions: $[,]$ at frequency $O(1)$,
$\psi_{\alpha}$ at frequency $<O(1)$: this case is handled exactly
as the preceding one, utilizing lemma~\ref{bilinear4}, hence left
out.
\\

This finishes the treatment of the first term
$(\Box\Phi)\psi_{\alpha}$ of I.
\\

\begin{bf} Second term of I\end{bf}
\\

We only need to look at the cases when at least two inputs of
$P_{0}[\partial_{\nu}\Phi\partial^{\nu}\Phi\psi_{\alpha}]$ have
large frequency, since otherwise, one can simply place all inputs
into $L_{t}^{3}L_{x}^{6}$. However, one can then simply use the
estimate provided by lemma~\ref{almost endpoint}, in particular
the third inequality, to handle this case: in detail

\begin{equation}\begin{split}
&||P_{0}\sum_{k_{1}>10, k_{2}\leq
k_{1}}[\partial_{\nu}P_{k_{1}}\Phi\partial^{\nu}P_{k_{2}}
\Phi\psi_{\alpha}]||_{N[0]}\\& +||P_{0}\sum_{k_{2}>10, k_{1}<k_{2}
}[\partial_{\nu}P_{k_{1}}\Phi\partial^{\nu}P_{k_{2}}\Phi\psi_{\alpha}]||_{N[0]}\\
&\leq \sum_{k_{1}>10,k_{2}\leq
k_{1}}C2^{-\delta_{1}k_{1}}2^{\delta_{2}\min\{0,k_{2}\}}M^{3}\epsilon
c_{k_{1}}c_{k_{2}}\\&+\sum_{k_{2}>10,k_{1}<k_{2}}2^{-\delta_{1}
k_{2}}2^{\delta_{2}\min\{0,k_{1}\}}M^{3}\epsilon c_{k_{1}}c_{k_{2}}\\
&\leq CM^{3}c_{0}\\
\end{split}\end{equation}

\subsection{Treating III}

Utilizing the Proposition~\ref{Gauge Change}, we see that we need
to estimate $||P_{0}(\phi^{1}_{\nu}\phi^{1\nu}\psi)||_{N[0]}$,
where $\psi$ is a Schwartz function whose frequency localized
components have $S[k]$-norms majorized by the frequency envelope
$CM c_{k}$ for appropriate $C$. The other terms in III are entirely analogous.\\
Of course, as before, this only requires consideration if at least
one (and then of course at least two) inputs have very large
frequency. Our strategy shall again be to employ dynamic
separation to render a hidden $Q_{0}$-structure visible. We need
to estimate the following expression:

\begin{equation}
||P_{0}[\sum_{k_{1}>
10}P_{k_{1}}\phi^{1}_{\nu}\phi^{1\nu}P_{<k_{1}+O(1)}\psi]||_{N[0]}
\end{equation}

We apply dynamic separation to $\phi^{1}_{\nu}$, referring to
\eqref{div-curl 1}, \eqref{div-curl 2}, and decompose the term
within $||.||$ into the following two terms:

\begin{equation}
\sum_{k=1}^{3}\sum_{k_{1}>10}P_{0}[P_{k_{1}}R_{\nu}R_{k}\phi^{1}_{k}\phi^{1\nu}P_{<k_{1}+O(1)}\psi]
\end{equation}

\begin{equation}
\sum_{k_{1}>10}P_{0}[P_{k_{1}}"\nabla^{-1}(\phi^{2})"\phi^{1\nu}P_{<k_{1}+O(1)}\psi]
\end{equation}

The 2nd term in the immediately preceding is easy to estimate,
placing all inputs into $L_{t}^{4}L_{x}^{4}$ and profiting from
the exponential gain $2^{-k_{1}}$ coming from $\nabla^{-1}$. \\
In order to estimate the first term, we need to apply another
dynamic separation. This results in the terms
\\

1):

\begin{equation}
\sum_{k,l=1}^{3}P_{0}[P_{k_{1}}R_{\nu}R_{k}\phi^{1}_{k}R^{\nu}R_{l}\phi^{1}_{l}P_{<k_{1}+O(1)}\psi]
\end{equation}

2):

\begin{equation}
\sum_{k,r=1}^{3}P_{0}[P_{k_{1}}R_{\nu}\phi^{1}_{k}\triangle^{-1}\partial_{r}(\phi^{1}_{r}\phi^{2}_{\nu}-\phi^{2}_{r}\phi^{1}_{\nu})P_{<k_{1}+O(1)}\psi]
\end{equation}

The first term can be immediately estimated upon referring to the
last inequality in lemma~\ref{almost endpoint}, as before. The 2nd
term is easily estimated by applying lemma~\ref{bilinear4} and the
estimates in II: we split into the case when $(,)$ is at frequency
$>2^{k_{1}-10}$ and its complement, i.e. $(,)$ at frequency $\leq
2^{k_{1}-10}$.
\\

First case: $(,)$ at frequency $>2^{k_{1}-10}$: Simply place all
entries into $L_{t}^{4}L_{x}^{4}$, using Bernstein's inequality
estimating $L_{t}^{1}L_{x}^{2}$  in terms of $L_{t}^{1}L_{x}^{1}$
for the output, as well as using the exponential decay
$2^{-k_{1}}$ coming from
$\partial_{r}\triangle^{-1}P_{>k_{1}-10}$.
\\

Second case: $(,)$ at frequency $\leq 2^{k_{1}-10}$: this can be
rewritten as (leaving out $\sum_{k,r=1}^{3}$ for simplicity's
sake):

\begin{equation}\begin{split}
&P_{0}[P_{k_{1}}R_{\nu}\phi^{1}_{k}\triangle^{-1}\partial_{r}P_{\leq
k_{1}-10}(\phi^{1}_{r}\phi^{2}_{\nu}-\phi^{2}_{r}\phi^{1}_{\nu})P_{<k_{1}+O(1)}\psi]\\
&=P_{0}[P_{k_{1}+O(1)}(P_{k_{1}}R_{\nu}\phi^{1}_{k}\triangle^{-1}\partial_{r}P_{\leq
k_{1}-10}(\phi^{1}_{r}\phi^{2}_{\nu}-\phi^{2}_{r}\phi^{1}_{\nu}))P_{<k_{1}+O(1)}\psi]\\
&=P_{0}[P_{k_{1}+O(1)}(P_{k_{1}}R_{\nu}\phi^{1}_{k}\triangle^{-1}\partial_{r}
(\phi^{1}_{r}\phi^{2}_{\nu}-\phi^{2}_{r}\phi^{1}_{\nu}))P_{k_{1}+O(1)}\psi]\\
&-P_{0}[P_{k_{1}+O(1)}(P_{k_{1}}R_{\nu}\phi^{1}_{k}\triangle^{-1}\partial_{r}P_{>
k_{1}-10}(\phi^{1}_{r}\phi^{2}_{\nu}-\phi^{2}_{r}\phi^{1}_{\nu}))P_{<k_{1}+O(1)}\psi]\\
\end{split}\end{equation}

The 2nd term in the last equality can be treated exactly as in the
first case, of course. As to the first term, we use
lemma~\ref{bilinear4} and the estimates proved in II to conclude
that

\begin{equation}\begin{split}
&||P_{0}[P_{k_{1}+O(1)}(P_{k_{1}}R_{\nu}\phi^{1}_{k}\triangle^{-1}\partial_{r}
(\phi^{1}_{r}\phi^{2}_{\nu}-\phi^{2}_{r}\phi^{1}_{\nu}))P_{<k_{1}+O(1)}\psi]||_{N[0]}\\
&\leq C2^{-\delta k_{1}}
||P_{k_{1}+O(1)}\nabla(P_{k_{1}}R_{\nu}\phi^{1}_{k}\triangle^{-1}\partial_{r}
(\phi^{1}_{r}\phi^{2}_{\nu}-\phi^{2}_{r}\phi^{1}_{\nu}))||_{N[k_{1}+O(1)]}||P_{k_{1}+O(1)}\psi||_{S[k_{1}+O(1)]}\\
&\leq CM^{4}2^{-\delta k_{1}}c_{k_{1}}^{2}\\
\end{split}\end{equation}

This can be easily summed over $k_{1}>10$ to yield the required
bound $CM^{4}\epsilon c_{0}$.
\\

The proof of Proposition~\ref{bootstrap} is now completed by
fixing $M>>1$, then choosing $\epsilon>0$ small enough such that
using \eqref{energy} and Theorem~\ref{Gauge Change}, as well as
the estimates proved in this section, we can conclude
$||P_{0}\phi||_{S[0]}\leq \frac{M}{2}c_{0}$.
\end{proof}

\newpage

\section{Proof of the technical estimates}

\subsection{The Trilinear Estimate}

\begin{proposition}
Let $\psi_{l},\, l=1,2,3$ be Schwartz functions on
$\mathbf{R}^{3+1}$. We then have the estimate
\begin{equation}\begin{split}
&||P_{0}(\sum_{j=1}^{3}\triangle^{-1}\partial_{j}[R_{\nu}P_{k_{1}}\psi_{1}R_{j}P_{k_{2}}\psi_{2}-R_{j}P_{k_{1}}\psi_{1}R_{\nu}P_{k_{2}}\psi_{2}]\partial^{\nu}P_{k_{3}}\psi_{3})||_{N[0]}\\&\leq
C 2^{-\delta_{1} |k_{1}-k_{2}|}2^{-\delta_{2}
|k_{3}|}\prod_{l=1}^{3}||P_{k_{l}}\psi_{l}||_{S[k_{l}]}\\
\end{split}\end{equation}
for appropriate constants $\delta_{1},\delta_{2}>0$.
\end{proposition}

\begin{proof}
The proof is split into estimating the various frequency
interactions of $[,]$ and $\partial^{\nu}\psi_{3}$. The author
apologizes in advance for the extremely technical and mechanical
nature of these estimates. Many of the following cases are similar
and almost write themselves, yet we have opted for as complete an
account as possible. This section may be treated as an appendix to
the paper.

\begin{bf}1): Low-High Interactions: $[,]$ at frequency $<2^{-10}$, $\partial^{\nu}\psi_{3}$ at frequency between $2^{-10}$ and $2^{10}$, i.e. $k_{3}=O(1)$.\end{bf}
\\
We distinguish between the cases
\\

\begin{bf}A\end{bf}: $k_{1}<k_{2}-10$\\
\begin{bf}B\end{bf}: $k_{1}>k_{2}+10$\\
\begin{bf}C\end{bf}: $|k_{1}-k_{2}|\leq 10$
\\

On account of the apparent symmetry between $A$ and $B$, only one
of them has to be treated. We estimate the corresponding terms by
restricting the Fourier supports of the inputs of the trilinear
form further. I.e., in addition to frequency localization, we also
localize in modulation, corresponding to the quantity
$||\tau|-|\xi||$. The general strategy in the estimates to follow
is to reduce output and inputs to small modulation, as it is only
then that the algebraic structure as exemplified by \eqref{full
null-form} becomes useful. However, when the output has "large
modulation", it can usually be easily estimated in the
$\dot{X}_{0}^{-\frac{1}{2},-\frac{1}{2},1}$-norm. Similarly,
whenever an input has "large modulation", it can be placed into
$L_{t}^{2}L_{x}^{2}$, which allows one more flexibility for the
other inputs.

\newpage
\begin{bf}1.A): Estimation of A\end{bf}
\\

\begin{bf}Either output or third input at large modulation:\end{bf}
\\

1.A.1):Output reduced to modulation $>2^{k_{2}-100}$: we have
(denoting $(\sqrt{-\triangle})^{-1}$ as $\nabla^{-1}$)

\begin{equation}\begin{split}
&||P_{0}Q_{>k_{2}-100}(\sum_{j=1}^{3}\triangle^{-1}\partial_{j}P_{<-10}[R_{\nu}P_{k_{1}}\psi_{1}R_{j}P_{k_{2}}\psi_{2}-R_{j}P_{k_{1}}\psi_{1}R_{\nu}P_{k_{2}}\psi_{2}]\\&P_{k_{3}}\partial^{\nu}\psi_{3})||_{N[0]}\\
&\leq \sum_{j=1}^{3}||P_{0}Q_{>k_{2}-100}(P_{k_{2}+O(1)}\triangle^{-1}\partial_{j}[R_{\nu}P_{k_{1}}\psi_{1}R_{j}P_{k_{2}}\psi_{2}-R_{j}P_{k_{1}}\psi_{1}R_{\nu}P_{k_{2}}\psi_{2}]\\&P_{k_{3}}\partial^{\nu}\psi_{3})||_{X^{-\frac{1}{2},-\frac{1}{2},1}}\\
&\leq C2^{-\frac{k_{2}}{2}}2^{-k_{2}}\sum_{j=1}^{3}||[R_{\nu}P_{k_{1}}\psi_{1}R_{j}P_{k_{2}}\psi_{2}-R_{j}P_{k_{1}}\psi_{1}R_{\nu}P_{k_{2}}\psi_{2}]||_{L_{t}^{2}L_{x}^{\infty}}\\&||P_{k_{3}}\partial^{\nu}\psi_{3}||_{L_{t}^{\infty}L_{x}^{2}}\\
&\leq
C||\nabla_{x,t}\nabla^{-1}P_{k_{1}}\psi_{1}||_{L_{t}^{3}L_{x}^{6}}||\nabla_{x,t}\nabla^{-1}\psi_{2}||_{L_{t}^{6}L_{x}^{3}}||\nabla_{x,t}P_{k_{3}}\psi_{3}||_{L_{t}^{\infty}L_{x}^{2}}\\
&\leq C
2^{\frac{k_{1}-k_{2}}{6}}\prod_{l=1}^{3}||P_{l}\psi_{l}||_{S[k_{l}]}\\
\end{split}\end{equation}

1.A.2): $\partial^{\nu}P_{k_{3}}\psi_{3}$ reduced to modulation
$>2^{k_{2}-100}$, output reduced to modulation $\leq
2^{k_{2}-100}$:

\begin{equation}\begin{split}
&||P_{0}Q_{\leq
k_{2}-100}[\sum_{j=1}^{3}\triangle^{-1}\partial_{j}P_{<-10}(R_{\nu}P_{k_{1}}\psi_{1}R_{j}P_{k_{2}}\psi_{2}-R_{j}P_{k_{1}}\psi_{1}R_{\nu}P_{k_{2}}\psi_{2})\\&P_{k_{3}}Q_{>
k_{2}-100}\partial^{\nu}\psi_{3}]||_{N[0]}\\
&\leq ||P_{0}Q_{\leq
k_{2}-100}[\sum_{j=1}^{3}P_{k_{2}+O(1)}\triangle^{-1}\partial_{j}(R_{\nu}P_{k_{1}}\psi_{1}R_{j}P_{k_{2}}\psi_{2}-R_{j}P_{k_{1}}\psi_{1}R_{\nu}P_{k_{2}}\psi_{2})\\&P_{k_{3}}Q_{>
k_{2}-100}\partial^{\nu}\psi_{3}]||_{L_{t}^{1}L_{x}^{2}}\\
&\leq\sum_{j=1}^{3}
C2^{-k_{2}}||R_{\nu}P_{k_{1}}\psi_{1}R_{j}P_{k_{2}}\psi_{2}-R_{j}P_{k_{1}}\psi_{1}R_{\nu}P_{k_{2}}\psi_{2}||_{L_{t}^{2}L_{x}^{\infty}}
\\&||P_{k_{3}}Q_{>
k_{2}-100}\partial^{\nu}\psi_{3}||_{L_{t}^{2}L_{x}^{2}}\\
&\leq C
2^{\frac{k_{1}-k_{2}}{6}}\prod_{l=1}^{3}||P_{k_{l}}\psi_{l}||_{S[k_{l}]}\\
\end{split}\end{equation}

\begin{bf}Output and third input at small modulation. Part of
null-structure becomes useful: \end{bf}
\\

1.A.3): Both output and $\partial^{\nu}P_{k_{3}}\psi_{3}$ are
reduced to modulation $\leq 2^{k_{2}-100}$: notice that this
implies in particular that $[,]$ can be restricted to modulation
$\leq 2^{k_{2}+O(1)}$.\\
We use the elementary identity

\begin{equation}\label{special identity}\begin{split}
&2\sum_{j=1}^{3}\triangle^{-1}\partial_{j}(R_{j}fR_{\nu}g-R_{\nu}fR_{j}g)\partial^{\nu}h\\
&=\Box[\triangle^{-1}\partial_{j}(R_{j}f\nabla^{-1}g)h]-\triangle^{-1}\partial_{j}(R_{j}f\nabla^{-1}g)\Box
h\\&-\Box(\partial_{j}\triangle^{-1}(R_{j}f\nabla^{-1}g))h -
2R_{\nu}f\nabla^{-1}g\partial^{\nu}h\\
\end{split}\end{equation}

valid for all Schwartz functions, say. We consecutively estimate
each of the terms, with $f$ etc. replaced by the appropriately
microlocalized components of $P_{k_{1}}\psi$ etc.:
\\

1.A.3.a):

\begin{equation}\begin{split}
&||\Box P_{0}Q_{\leq
k_{2}-100}[\triangle^{-1}\partial_{j}(P_{k_{1}}R_{j}\psi_{1}\nabla^{-1}P_{k_{2}}\psi_{2})P_{k_{3}}Q_{\leq
k_{2}-100}\psi_{3}]||_{N[0]}\\
&\leq ||\Box P_{0}Q_{\leq
k_{2}-100}[\triangle^{-1}\partial_{j}P_{k_{2}+O(1)}(P_{k_{1}}R_{j}\psi_{1}\nabla^{-1}P_{k_{2}}\psi_{2})P_{k_{3}}Q_{\leq
k_{2}-100}\psi_{3}]||_{\dot{X}_{0}^{-\frac{1}{2},-\frac{1}{2},1}}\\
&\leq C
2^{-\frac{k_{2}}{2}}||P_{k_{2}+O(1)}(P_{k_{1}}R_{j}\psi_{1}\nabla^{-1}P_{k_{2}}\psi_{2})||_{L_{t}^{2}L_{x}^{\infty}}||P_{k_{3}}Q_{\leq
k_{2}-100}\psi_{3}||_{L_{t}^{\infty}L_{x}^{2}}\\
&\leq C
2^{\frac{k_{1}-k_{2}}{6}}\prod_{l=1}^{3}||P_{k_{l}}\psi_{l}||_{S[k_{l}]}\\
\end{split}\end{equation}

1.A.3.b):

\begin{equation}\begin{split}
&||P_{0}Q_{\leq k_{2}-100}[\Box
\triangle^{-1}\partial_{j}(R_{j}P_{k_{1}}\psi_{k_{1}}\nabla^{-1}P_{k_{2}}\psi_{2})P_{k_{3}}Q_{\leq
k_{2}-100}\psi_{3}||_{N[0]}\\
&\leq C||P_{0}Q_{\leq k_{2}-100}[\Box
\triangle^{-1}\partial_{j}P_{k_{2}+O(1)}Q_{\leq k_{2}+O(1)}
(R_{j}P_{k_{1}}\psi_{k_{1}}\nabla^{-1}P_{k_{2}}\psi_{2})\\&P_{k_{3}}Q_{\leq
k_{2}-100}\psi_{3}]||_{L_{t}^{1}L_{x}^{2}}\\
&\leq C2^{\frac{k_{1}+k_{2}}{6}}\prod_{l=1}^{3}||P_{k_{l}}\psi_{l}||_{S[k_{l}]}\\
\end{split}\end{equation}

as is seen by placing all entries into $L_{t}^{3}L_{x}^{6}$. Of
course this is even better than the estimate required for the
lemma, since $k_{1},k_{2}<O(1)$.
\\

1.A.3.c):

\begin{equation}\begin{split}
&||P_{0}Q_{\leq k_{2}-100}[
\triangle^{-1}\partial_{j}(R_{j}P_{k_{1}}\psi_{k_{1}}\nabla^{-1}P_{k_{2}}\psi_{2})\Box
P_{k_{3}}Q_{\leq
k_{2}-100}\psi_{3}]||_{N[0]}\\
&\leq
C||P_{k_{2}+O(1)}\triangle^{-1}\partial_{j}(R_{j}P_{k_{1}}\psi_{k_{1}}\nabla^{-1}P_{k_{2}}\psi_{2})||_{L_{t}^{2}L_{x}^{\infty}}
||\Box P_{k_{3}}Q_{\leq
k_{2}-100}\psi_{3}||_{L_{t}^{2}L_{x}^{2}}\\
&\leq C 2^{\frac{k_{1}-k_{2}}{6}}\prod_{l=1}^{3}||P_{k_{l}}\psi_{l}||_{S[k_{l}]}\\
\end{split}\end{equation}

again by placing $P_{k_{1}}\psi_{1}, P_{k_{2}}\psi_{2}$ into
$L_{t}^{3}L_{x}^{6}, L_{t}^{6}L_{x}^{3}$ respectively.
\\

1.A.3.d):
\\

$||P_{0}Q_{\leq
k_{2}-100}[R_{\nu}P_{k_{1}}\psi_{1}\nabla^{-1}P_{k_{2}}\psi_{2}\partial^{\nu}P_{k_{3}}Q_{\leq
k_{2}-100} \psi_{3}||_{N[0]}$.
\\

This term is more difficult and will be estimated with the aid of
the '$Q_{0}$-structure' inherent in it (viz. the complete
expansion \eqref{full null-form}). It really is a consequence of a
deep trilinear inequality proved by T.Tao, but we prove it here
for completeness' sake, and also since the proof in 3 dimensions
is somewhat simpler: first, we get rid of the multiplier $Q_{\leq
k_{2}-100}$ in front of it as well as in front of the third input;
for example:

\begin{equation}\begin{split}
&||P_{0}Q_{>k_{2}-100}[R_{\nu}P_{k_{1}}\psi_{1}\nabla^{-1}P_{k_{2}}\psi_{2}\partial^{\nu}P_{k_{3}}Q_{\leq k_{2}-100}\psi_{3}]||_{\dot{X}_{0}^{-\frac{1}{2},-\frac{1}{2},1}}\\
&\leq
C 2^{-\frac{k_{2}}{2}}||\nabla_{x,t}\nabla^{-1}P_{k_{1}}\psi_{1}||_{L_{t}^{4}L_{x}^{\infty}}||\nabla^{-1}P_{k_{2}}\psi_{2}||_{L_{t}^{4}L_{x}^{\infty}}||P_{k_{3}}Q_{\leq k_{2}-100}\partial^{\nu}\psi_{3}||_{L_{t}^{\infty}L_{x}^{2}}\\
&\leq C2^{\frac{3(k_{1}-k_{2})}{4}}\prod_{l=1}^{3}||P_{k_{l}}\psi_{l}||_{S[k_{l}]}\\
\end{split}\end{equation}

Let us restate what we are trying to prove

\begin{lemma}\label{low-high interactions}
Under the hypotheses of case 1.A), we have the estimate

\begin{equation}
||P_{0}[R_{\nu}P_{k_{1}}\psi_{1}\nabla^{-1}P_{k_{2}}\psi_{2}\partial^{\nu}P_{k_{3}}\psi_{3}||_{N[0]}
\leq C
2^{\delta(k_{1}-k_{2})}\prod_{l=1}^{3}||P_{k_{l}}\psi_{l}||_{S[k_{l}]}
\end{equation}

for some $\delta >0$.
\end{lemma}

\begin{proof}: The proof is again by considering many different
cases. However, the cases are somewhat more involved than in the
preceding; in particular, we will have to invoke
lemma~\ref{bilinear1} for some "intermediate cases" where inputs
have relatively small modulation but still too large to invoke the
algebraic identity. We first consider the case when
$P_{k_{1}}\psi_{1}$ has relatively large modulation, i.e.
modulation $\geq 2^{k_{1}+O(1)}$. We cannot yet invoke the
null-structure in this case. However, lemma~\ref{bilinear1} allows
us to play the modulation of $P_{k_{1}}\psi_{1}$ against the
frequency of $P_{k_{2}}\psi_{2}$:
\\

\begin{bf}1): $P_{k_{1}}\psi_{1}$ restricted to modulation $\geq
2^{k_{1}+100}$.\end{bf}
\\

1.1.1) $P_{k_{1}}\psi_{1}$ at modulation $2^{r}$ with
$k_{1}+100\leq r\leq k_{2}+100$. $P_{k_{3}}\psi_{3}$ at modulation
$<2^{r-100}$: From elementary geometrical considerations, we see
that we in this case
$R_{\nu}P_{k_{1}}Q_{r}\psi_{1}P_{k_{3}}Q_{<r-100}\psi_{3}$ has
Fourier support at distance approximately $2^{r}$ from the light
cone. Hence using lemma~\ref{bilinear1}, we can play this
modulation against the frequency $2^{k_{2}}$ of the third term
$P_{k_{2}}\psi_{2}$. Of course, we then have to also play $r$
against $k_{1}$, which is possible by placing
$P_{k_{1}}Q_{r}\psi_{1}$ into $L_{t}^{2}L_{x}^{\infty}$. The
upshot is an exponential gain in $k_{1}-k_{2}$. In detail:

\begin{equation}\begin{split}
&||\sum_{k_{2}+100\geq r\geq
k_{1}+100}P_{0}[R_{\nu}P_{k_{1}}Q_{r}\psi_{1}\nabla^{-1}P_{k_{2}}\psi_{2}\partial^{\nu}P_{k_{3}}Q_{<r-100}\psi_{3}||_{N[0]}\\
&\leq\sum_{k_{2}+100\geq r\geq k_{1}+100}C
2^{\delta(r-k_{2})}||P_{O(1)}Q_{r+O(1)}(R_{\nu}P_{k_{1}}Q_{r}\psi_{1}P_{k_{3}}Q_{<r-100}\psi_{3})||_{\dot{X}_{O(1)}^{\frac{1}{2},-\frac{1}{2},\infty}}
\\&||P_{k_{2}}\psi_{2}||_{S[k_{2}]}\\
&\leq\sum_{k_{2}+100\geq r\geq
k_{1}+100}C2^{-\frac{r}{2}}2^{\delta(r-k_{2})}||R_{\nu}P_{k_{1}}Q_{r}\psi_{1}||_{L_{t}^{2}L_{x}^{\infty}}||P_{k_{3}}Q_{<r-100}\psi_{3}||_{L_{t}^{\infty}L_{x}^{2}}
||P_{k_{2}}\psi_{2}||_{S[k_{2}]}\\
&\leq\sum_{k_{2}+100\geq r\geq
k_{1}+100}C2^{\delta(r-k_{2})}2^{k_{1}-r}\prod_{l=1}^{3}||P_{k_{l}}\psi_{l}||_{S[k_{l}]}\\
&\leq C2^{\delta(k_{1}-k_{2})}\prod_{l=1}^{3}||P_{k_{l}}\psi_{l}||_{S[k_{l}]}\\
\end{split}\nonumber\end{equation}

1.1.2.a) $P_{k_{1}}\psi_{1}$ at modulation $2^{r}$ with
$k_{1}+100\leq r\leq k_{2}+100$, $P_{k_{3}}\psi_{3}$ at modulation
$\geq 2^{r-100}$ but $\leq 2^{k_{2}+100}$, $P_{k_{2}}\psi_{2}$ at
modulation $<2^{r-100}$:  we need to group the terms differently,
again using lemma~\ref{bilinear1}. Note that the
$\dot{X}_{k}^{\frac{1}{2},-\frac{1}{2},\infty}$-norm has to be
applied to a term at large frequency in order to avoid an
exponential loss. In the present situation, the third term
$P_{k_{3}}Q_{\geq r-100}\psi_{3}$ is such a candidate. However, in
order to apply lemma~\ref{bilinear1}, we then need to place
$R_{\nu}P_{k_{1}}Q_{r}\psi_{1}\nabla^{-1}P_{k_{2}}Q_{<r-100}\psi_{2}$
into $S[k_{2}+O(1)]$. This is easily feasible on account of the
majorization $||P_{k}\phi||_{S[k]}\leq
C||P_{k}\phi||_{\dot{X}_{k}^{\frac{1}{2},\frac{1}{2},1}}$, which
follows for example from \eqref{cap localized}:

\begin{equation}\begin{split}
&||\sum_{k_{2}+100\geq r\geq k_{1}+100}\sum_{k_{2}+100\geq a\geq
r-100}\\&P_{0}[R_{\nu}P_{k_{1}}Q_{r}\psi_{1}\nabla^{-1}P_{k_{2}}Q_{<r-100}\psi_{2}\partial^{\nu}P_{k_{3}}Q_{a}\psi_{3}]||_{N[0]}\\
&\leq\sum_{k_{2}+100\geq r\geq k_{1}+100}\sum_{k_{2}+100\geq a\geq
r-100}\\&C
2^{\delta(a-k_{2})}||Q_{r+O(1)}\nabla[R_{\nu}P_{k_{1}}Q_{r}\psi_{1}\nabla^{-1}P_{k_{2}}Q_{<r-100}\psi_{2}]||_{S[k_{2}+O(1)]}||P_{k_{3}}Q_{a}\partial^{\nu}\psi_{3}||_{\dot{X}_{k_{3}}^{\frac{1}{2},-\frac{1}{2},\infty}}\\
&\leq\sum_{k_{2}+100\geq r\geq k_{1}+100}\sum_{k_{2}+100\geq a\geq
r-100}\\&C2^{\delta(a-k_{2})}||Q_{r+O(1)}\nabla[R_{\nu}P_{k_{1}}Q_{r}\psi_{1}\nabla^{-1}P_{k_{2}}Q_{<r-100}\psi_{2}]||_{\dot{X}_{k_{2}+O(1)}^{\frac{1}{2},\frac{1}{2},1}}||P_{k_{3}}Q_{a}\partial^{\nu}\psi_{3}||_{\dot{X}^{\frac{1}{2},-\frac{1}{2},\infty}}\\
&\leq\sum_{k_{2}+100\geq r\geq k_{1}+100}\sum_{k_{2}+100\geq a\geq
r-100}\\&C2^{\delta(a-k_{2})}2^{k_{2}}2^{-\frac{r}{2}}2^{-\frac{k_{1}}{2}}2^{\frac{3k_{1}}{2}}2^{-\frac{3k_{2}}{2}}2^{-a}2^{\frac{k_{2}}{2}}2^{\frac{r}{2}}\prod_{l=1}^{3}||P_{k_{l}}\psi_{l}||_{S[k_{l}]}\\
&\leq\sum_{k_{2}+100\geq r\geq k_{1}+100}\sum_{k_{2}+100\geq a\geq
r-100}
C2^{k_{1}-a}2^{\delta(a-k_{2})}\prod_{l=1}^{3}||P_{k_{l}}\psi_{l}||_{S[k_{l}]}\\
&\leq C2^{\delta'(k_{1}-k_{2})}\prod_{l=1}^{3}||P_{k_{l}}\psi_{l}||_{S[k_{l}]}\\
\end{split}\nonumber\end{equation}

for any $0<\delta'<\delta$.
\\

1.1.2.b) $P_{k_{1}}\psi_{1}$ at modulation $2^{r}$ with
$k_{1}+100\leq r\leq k_{2}+100$, $P_{k_{3}}\psi_{3}$ at modulation
$\geq 2^{r-100}$ but $\leq 2^{k_{2}+100}$, $P_{k_{2}}\psi_{2}$ at
modulation $\geq 2^{r-100}$: this and the following 4 large
modulation cases are routine and can be handled by the
$\dot{X}_{k}^{\frac{1}{2},\frac{1}{2},\infty}$ and
$L_{t}^{\infty}L_{x}^{2}$-components of $S[k]$. The only slight
difficulty here is that one has to play the modulation of the term
$P_{k_{2}}Q_{a}\psi_{2}$ against its own frequency, i.e. one has
to obtain an exponential gain in the difference $a-k_{2}$. This
requires the improvement of Bernstein's inequality mentioned in
the third section, i.e. \eqref{improved Bernstein}.

\begin{equation}\begin{split}
&||\sum_{k_{2}+100\geq r\geq k_{1}+100}\sum_{a\geq
r-100}P_{0}[R_{\nu}P_{k_{1}}Q_{r}\psi_{1}\nabla^{-1}P_{k_{2}}Q_{a}\psi_{2}\partial^{\nu}P_{k_{3}}Q_{k_{2}+100\geq
.\geq r-100}\psi_{3}]||_{L_{t}^{1}L_{x}^{2}}\\
&\leq\sum_{k_{2}+100\geq r\geq k_{1}+100}\sum_{a\geq
r-100}C||R_{\nu}P_{k_{1}}Q_{r}\psi_{1}||_{L_{t}^{2}L_{x}^{\infty}}||\nabla^{-1}P_{k_{2}}Q_{a}\psi_{2}||_{L_{t}^{2}L_{x}^{\infty}}
\\&||\partial^{\nu}P_{k_{3}}Q_{k_{2}+100\geq .\geq
r-100}\psi_{3}||_{L_{t}^{\infty}L_{x}^{2}}\\
&\leq\sum_{k_{2}+100\geq r\geq k_{1}+100}\sum_{a\geq
r-100}C2^{-\frac{r}{2}}2^{-\frac{k_{1}}{2}}2^{\frac{3k_{1}}{2}}2^{\min\{0,\delta(a-k_{2})\}}2^{-\frac{a}{2}}2^{-\frac{k_{2}}{2}}2^{\frac{3k_{2}}{2}}2^{-k_{2}}
\\&\prod_{l=1}^{3}||P_{k_{l}}\psi_{l}||_{S[k_{l}]}\\
&\leq C2^{\delta(k_{1}-k_{2})}\prod_{l=1}^{3}||P_{k_{l}}\psi_{l}||_{S[k_{l}]}\\
\end{split}\nonumber\end{equation}

We have used \eqref{improved Bernstein} in the 2nd step.
\\

1.1.2.c)  $P_{k_{1}}\psi_{1}$ at modulation $2^{r}$ with
$k_{1}+100\leq r\leq k_{2}+100$, $P_{k_{3}}\psi_{3}$ at modulation
$>2^{k_{2}+100}$. We have

\begin{equation}\begin{split}
&||\sum_{k_{2}+100\geq r\geq
k_{1}+100}P_{0}[R_{\nu}P_{k_{1}}Q_{r}\psi_{1}\nabla^{-1}P_{k_{2}}\psi_{2}\partial^{\nu}P_{k_{3}}Q_{>k_{2}+100}\psi_{3}||_{L_{t}^{1}L_{x}^{2}}\\
&\leq\sum_{k_{2}+100\geq r\geq
k_{1}+100}C||R_{\nu}P_{k_{1}}Q_{r}\psi_{1}||_{L_{t}^{2}L_{x}^{\infty}}||\nabla^{-1}P_{k_{2}}\psi_{2}||_{L_{t}^{\infty}L_{x}^{\infty}}||\partial^{\nu}P_{k_{3}}Q_{>k_{2}+100}\psi_{3}||_{L_{t}^{2}L_{x}^{2}}\\
&\leq\sum_{k_{2}+100\geq r\geq
k_{1}+100}C2^{-\frac{r}{2}}2^{-\frac{k_{1}}{2}}2^{\frac{3k_{1}}{2}}2^{-\frac{k_{2}}{2}}\prod_{l=1}^{3}||P_{k_{l}}\psi_{l}||_{S[k_{l}]}\\
&\leq\sum_{k_{2}+100\geq r\geq
k_{1}+100}C2^{k_{1}-\frac{r}{2}-\frac{k_{2}}{2}}\prod_{l=1}^{3}||P_{k_{l}}\psi_{l}||_{S[k_{l}]}\\
&\leq C2^{\frac{k_{1}-k_{2}}{2}}\prod_{l=1}^{3}||P_{k_{l}}\psi_{l}||_{S[k_{l}]}\\
\end{split}\end{equation}

1.2.1) $P_{k_{1}}\psi_{1}$ at modulation $>2^{k_{2}+100}$,
$P_{k_{2}}\psi_{2}$ at modulation $\geq 2^{k_{2}-100}$:

\begin{equation}\begin{split}
&||P_{0}[R_{\nu}P_{k_{1}}Q_{>k_{2}+100}\psi_{1}\nabla^{-1}P_{k_{2}}Q_{\geq
k_{2}-100}\psi_{2}\partial^{\nu}P_{k_{3}}\psi_{3}]||_{L_{t}^{1}L_{x}^{2}}\\
&\leq C
||R_{\nu}P_{k_{1}}Q_{>k_{2}+100}\psi_{1}||_{L_{t}^{2}L_{x}^{\infty}}||\nabla^{-1}P_{k_{2}}Q_{\geq
k_{2}-100}\psi_{2}||_{L_{t}^{2}L_{x}^{\infty}}||\partial^{\nu}P_{k_{3}}\psi_{3}||_{L_{t}^{\infty}L_{x}^{2}}\\
&\leq C2^{k_{1}-k_{2}}\prod_{l=1}^{3}||P_{k_{l}}\psi_{l}||_{S[k_{l}]}\\
\end{split}\end{equation}

1.2.2) $P_{k_{1}}\psi_{1}$ at modulation $>2^{k_{2}+100}$,
$P_{k_{2}}\psi_{2}$ at modulation $<2^{k_{2}-100}$,
$P_{k_{3}}\psi_{3}$ at modulation $\geq 2^{k_{2}-100}$:

\begin{equation}\begin{split}
&||P_{0}[R_{\nu}P_{k_{1}}Q_{>k_{2}+100}\psi_{1}\nabla^{-1}P_{k_{2}}Q_{<
k_{2}-100}\psi_{2}\partial^{\nu}P_{k_{3}}Q_{\geq k_{2}-100}\psi_{3}]||_{L_{t}^{1}L_{x}^{2}}\\
&\leq C
||R_{\nu}P_{k_{1}}Q_{>k_{2}+100}\psi_{1}||_{L_{t}^{2}L_{x}^{\infty}}||\nabla^{-1}P_{k_{2}}Q_{<
k_{2}-100}\psi_{2}||_{L_{t}^{\infty}L_{x}^{\infty}}||\partial^{\nu}P_{k_{3}}Q_{\geq
k_{2}-100}\psi_{3}||_{L_{t}^{2}L_{x}^{2}}\\
&\leq C2^{k_{1}-k_{2}}\prod_{l=1}^{3}||P_{k_{l}}\psi_{l}||_{S[k_{l}]}\\
\end{split}\end{equation}

1.2.3) $P_{k_{1}}\psi_{1}$ at modulation $>2^{k_{2}+100}$,
$P_{k_{2}}\psi_{2}$ at modulation $<2^{k_{2}-100}$,
$P_{k_{3}}\psi_{3}$ at modulation $<2^{k_{2}-100}$: note that
provided $P_{k_{1}}\psi_{1}$ is at modulation $2^{r}$,
$r>k_{2}+100$, then the output is at modulation $\sim 2^{r}$.
Hence

\begin{equation}\begin{split}
&||P_{0}[\sum_{r>k_{2}+100}R_{\nu}P_{k_{1}}Q_{r}\psi_{1}\nabla^{-1}P_{k_{2}}Q_{<k_{2}-100}\psi_{2}\partial^{\nu}P_{k_{3}}Q_{<k_{2}-100}\psi_{3}]||_{N[0]}\\
&\leq\sum_{r>k_{2}+100}||P_{0}Q_{r+O(1)}[R_{\nu}P_{k_{1}}Q_{r}\psi_{1}[\nabla^{-1}P_{k_{2}}Q_{<k_{2}-100}\psi_{2}\partial^{\nu}P_{k_{3}}Q_{<k_{2}-100}\psi_{3}]||_{\dot{X}_{k_{3}}^{-\frac{1}{2},-\frac{1}{2},1}}\\
&\leq\sum_{r>k_{2}+100}C2^{-\frac{r}{2}}||R_{\nu}P_{k_{1}}Q_{r}\psi_{1}||_{L_{t}^{2}L_{x}^{\infty}}||\nabla^{-1}P_{k_{2}}Q_{<k_{2}-100}\psi_{2}||_{L_{t}^{\infty}L_{x}^{\infty}}\\&
||\partial^{\nu}P_{k_{3}}Q_{<k_{2}-100}\psi_{3}||_{L_{t}^{\infty}L_{x}^{2}}\\
&\leq\sum_{r>k_{2}+100}C2^{-\frac{r}{2}}2^{-\frac{r}{2}}2^{-\frac{k_{1}}{2}}2^{\frac{3k_{1}}{2}}\prod_{l=1}^{3}||P_{k_{l}}\psi_{l}||_{S[k_{l}]}\\
&\leq C2^{k_{1}-k_{2}}\prod_{l=1}^{3}||P_{k_{l}}\psi_{l}||_{S[k_{l}]}\\
\end{split}\end{equation}

\begin{bf}2): $P_{k_{1}}\psi_{1}$ at modulation
$<2^{k_{1}+100}$.\end{bf}
\\

First, reduce to the case that $P_{k_{3}}\psi_{3}$ has modulation
$<2^{k_{1}-100}$: note that

\begin{equation}\begin{split}
&||P_{0}[R_{\nu}P_{k_{1}}Q_{<k_{1}+100}\psi_{1}\nabla^{-1}P_{k_{2}}\psi_{2}\partial^{\nu}P_{k_{3}}Q_{\geq
k_{1}-100} \psi_{3}||_{L_{t}^{1}L_{x}^{2}}\\
&\leq C
||R_{\nu}P_{k_{1}}Q_{<k_{1}+100}\psi_{1}||_{L_{t}^{4}L_{x}^{\infty}}||\nabla^{-1}P_{k_{2}}\psi_{2}||_{L_{t}^{4}L_{x}^{\infty}}||P_{k_{3}}Q_{\geq
k_{1}-100}\partial^{\nu}\psi_{3}||_{L_{t}^{2}L_{x}^{2}}\\
&\leq C2^{\frac{k_{1}-k_{2}}{4}}\prod_{l=1}^{3}||P_{k_{l}}\psi_{l}||_{S[k_{l}]}\\
\end{split}\end{equation}

Finally, we have manoeuvred ourselves into a position to exploit
the algebraic structure inherent in the trilinear form of the
lemma, i.e. the '2nd half in \eqref{full null-form}'. We thus have
to estimate the following terms:
\\

2.A):

\begin{equation}\label{almost done}
\Box(\nabla^{-1}P_{k_{1}}Q_{<k_{1}+100}\psi_{1}P_{k_{3}}Q_{<k_{1}-100}\psi_{3})\nabla^{-1}P_{k_{2}}\psi_{2}
\end{equation}

We first reduce $\nabla^{-1}P_{k_{1}}\psi$ to modulation
$<2^{k_{1}-100}$:

\begin{equation}\begin{split}
&||\Box(\nabla^{-1}P_{k_{1}}Q_{k_{1}-100\leq .
<k_{1}+100}\psi_{1}P_{k_{3}}Q_{<k_{1}-100}\psi_{3})\nabla^{-1}P_{k_{2}}\psi_{2}||_{N[0]}\\
&\leq
C2^{\delta(k_{1}-k_{2})}||Q_{<k_{1}+O(1)}\Box(\nabla^{-1}P_{k_{1}}Q_{k_{1}-100\leq
.
<k_{1}+100}\psi_{1}P_{k_{3}}Q_{<k_{1}-100}\psi_{3}||_{\dot{X}_{k_{3}+O(1)}^{\frac{1}{2},-\frac{1}{2},\infty}}\\&||P_{k_{2}}\psi_{2}||_{S[k_{2}]}\\
&\leq
2^{\delta(k_{1}-k_{2})}2^{\frac{k_{1}}{2}}||\nabla^{-1}P_{k_{1}}Q_{k_{1}-100\leq
.
<k_{1}+100}\psi_{1}||_{L_{t}^{2}L_{x}^{\infty}}||P_{k_{3}}\psi_{3}||_{L_{t}^{\infty}L_{x}^{2}}||P_{k_{2}}\psi_{2}||_{S[k_{2}]}\\
&\leq C2^{\delta(k_{1}-k_{2})}\prod_{l=1}^{3}||P_{k_{l}}\psi_{l}||_{S[k_{l}]}\\
\end{split}\end{equation}

We want to place $(,)$ into $L_{t}^{2}L_{x}^{2}$, then apply
lemma~\ref{bilinear1}. However, in order to avoid an exponential
loss in $k_{1}$, we would have to place
$\nabla^{-1}P_{k_{1}}\psi_{1}$ into $L_{t}^{2}L_{x}^{\infty}$,
which is impossible. Hence we need to invoke null-frame spaces,
and in particular the identity \eqref{crux}:
\\

2.A.1): Both inputs of $(,)$ have modulation $<<$ than the
modulation of $(,)$: assume that $(,)$ is localized to modulation
$2^{r}$. Hence by our assumptions $ r\leq k_{1}+O(1)$. Now assume
that the inputs of $(,)$ have modulation $\leq 2^{r-100}$. Then
using lemma(13.2) in \cite{Tao 2}, microlocalizing the inputs to
the upper or lower half-space, we can restrict the projections to
$S^{2}$ of the Fourier supports of the inputs to spherical caps
$\pm\kappa,\pm\kappa'$ of size $2^{\frac{r-k_{1}}{2}-50}$ at
distance $C2^{\frac{r-k_{1}}{2}}$, where the $\pm$-signs are
assigned corresponding to whether the function is microlocalized
in the upper or lower half-space $\tau><0$. Utilizing \eqref{crux}
yields

\begin{equation}\begin{split}
&||Q_{r}(\nabla^{-1}P_{k_{1}}Q_{<r-100}\psi_{1}P_{k_{3}}Q_{<r-100}\psi_{3}||_{L_{t}^{2}L_{x}^{2}}\\
&\leq 4\sup_{\pm,\pm}\sum_{\kappa,\kappa'\in K_{\frac{r-k}{2}-50},
\mbox{dist}(\pm\kappa,\pm\kappa')=C2^{\frac{r-k_{1}}{2}}}\\&||Q_{r}(\nabla^{-1}P_{k_{1},\kappa}Q^{\pm}_{<r-100}\psi_{1}P_{k_{3},\kappa'}Q^{\pm}_{<r-100}\psi_{3}||_{L_{t}^{2}L_{x}^{2}}\\
&\leq C\sup_{\pm,\pm}\sum_{\kappa,\kappa'\in K_{\frac{r-k}{2}-50},
\mbox{dist}(\pm\kappa,\pm\kappa')=C2^{\frac{r-k_{1}}{2}}}\\&||\nabla^{-1}P_{k_{1},\kappa}Q^{\pm}_{<r-100}\psi_{1}P_{k_{3},\kappa'}Q^{\pm}_{r-100\leq
.< r-k_{1}+ k_{3}-100}\psi_{3}||_{L_{t}^{2}L_{x}^{2}}\\
&+C\sup_{\pm,\pm}\sum_{\kappa,\kappa'\in K_{\frac{r-k}{2}-50},
\mbox{dist}(\pm\kappa,\pm\kappa')=C2^{\frac{r-k_{1}}{2}}}\\&||\nabla^{-1}P_{k_{1}\kappa}Q^{\pm}_{<r-100}\psi_{1}P_{k_{3},\kappa'}Q^{\pm}_{<r-k_{1}+k_{3}-100}\psi_{3}||_{L_{t}^{2}L_{x}^{2}}\\
\end{split}\nonumber\end{equation}

Now

\begin{equation}\begin{split}
&\sup_{\pm,\pm}\sum_{\kappa,\kappa'\in K_{\frac{r-k}{2}-50},
\mbox{dist}(\pm\kappa,\pm\kappa')=C2^{\frac{r-k_{1}}{2}}}||\nabla^{-1}P_{k_{1},\kappa}Q^{\pm}_{<r-100}\psi_{1}P_{k_{3},\kappa'}Q^{\pm}_{<
r-k_{1}+ k_{3}-100}\psi_{3}||_{L_{t}^{2}L_{x}^{2}}\\
&\leq C\sup_{\pm,\pm}\sum_{\kappa,\kappa'\in K_{\frac{r-k}{2}-50},
\mbox{dist}(\pm\kappa,\pm\kappa')=C2^{\frac{r-k_{1}}{2}}}2^{-\frac{k_{1}}{2}}||P_{k_{1},\kappa}Q^{\pm}_{<r-100}\psi_{1}||_{S[k_{1},\pm\kappa]}\\
&||P_{k_{3},\kappa'}Q^{\pm}_{\leq r-k_{1}+
k_{3}-100}\psi_{3}||_{S[k_{3},\pm\kappa']}\\
&\leq C2^{-\frac{k_{1}}{2}}\sup_{\pm,\pm}(\sum_{\kappa\in
K_{\frac{r-k_{1}}{2}-50}}||P_{k_{1},\pm\kappa}Q^{\pm}_{<r-100}\psi_{1}||_{S[k_{1},\kappa]}^{2})^{\frac{1}{2}}\\
&(\sum_{\kappa\in
K_{\frac{r-k_{1}}{2}-50}}||P_{k_{3},\pm\kappa}Q^{\pm}_{<r-k_{1}+k_{3}-100}\psi_{3}||_{S[k_{3},\kappa]}^{2})^{\frac{1}{2}}\\
&\leq C
2^{-\frac{k_{1}}{2}}||P_{k_{1}}\psi_{1}||_{S[k_{1}]}||P_{k_{3}}\psi_{3}||_{S[k_{3}]}\\
\end{split}\nonumber\end{equation}

Moreover

\begin{equation}\begin{split}
&\sup_{\pm,\pm}\sum_{\kappa,\kappa'\in K_{\frac{r-k}{2}-50},
\mbox{dist}(\pm\kappa,\pm\kappa')=C2^{\frac{r-k_{1}}{2}}}\\&
||\nabla^{-1}P_{k_{1},\kappa}Q^{\pm}_{<r-100}\psi_{1}P_{k_{3},\kappa'}Q^{\pm}_{r-100\leq
.< r-k_{1}+ k_{3}-100}\psi_{3}||_{L_{t}^{2}L_{x}^{2}}\\
&\leq C\sup_{\pm,\pm}\sum_{\kappa,\kappa'\in
K_{\frac{r-k_{1}}{2}-50},\mbox{dist}(\kappa,\kappa') = C
2^{\frac{r-k_{1}}{2}}}||\nabla^{-1}P_{k_{1},\kappa}Q^{\pm}_{<r-100}\psi_{1}||_{L_{t}^{\infty}L_{x}^{\infty}}\\
&||P_{k_{3},\kappa'}Q^{\pm}_{r-100\leq .< r-k_{1}+
k_{3}-100}\psi_{3}||_{L_{t}^{2}L_{x}^{2}}\\
&\leq C\sup_{\pm,\pm}\sum_{\kappa,\kappa'\in
K_{\frac{r-k_{1}}{2}-50},\mbox{dist}(\kappa,\kappa') = C
2^{\frac{r-k_{1}}{2}}}C2^{\frac{k_{1}}{2}}||P_{k_{1},\kappa}Q^{\pm}_{<r-100}\psi_{1}||_{L_{t}^{\infty}L_{x}^{2}}\\
&||P_{k_{3},\kappa'}Q^{\pm}_{r-100\leq .<r-k_{1}+
k_{3}-100}\psi_{3}||_{L_{t}^{2}L_{x}^{2}}\\
&\leq C 2^{-\frac{r}{2}}||P_{k_{1}}\psi_{1}||_{S[k_{1}]}||P_{k_{3}}\psi_{3}||_{S[k_{3}]}\\
\end{split}\nonumber\end{equation}

Hence inserting the preceding two results into \eqref{almost done}
and utilizing lemma~\ref{bilinear1}, we deduce

\begin{equation}\begin{split}
&||\sum_{r<k_{1}+O(1)}\Box
Q_{r}(\nabla^{-1}P_{k_{1}}Q_{<r-100}\psi_{1}P_{k_{3}}Q_{<r-100}\psi_{3})\nabla^{-1}P_{k_{2}}\psi_{2}||_{N[0]}\\
&\leq
\sum_{r<k_{1}+O(1)}C2^{\delta(r-k_{2})}\prod_{l=1}^{3}||P_{k_{l}}\psi_{k_{l}}||_{S[k_{l}]}\\
&\leq C2^{\delta(k_{1}-k_{2})}\prod_{l=1}^{3}||P_{k_{l}}\psi_{k_{l}}||_{S[k_{l}]}\\
\end{split}\end{equation}

2.A.2): At least one input of $(,)$ has modulation $\geq
2^{r-100}$:

\begin{equation}\begin{split}
&||Q_{r}(\nabla^{-1}P_{k_{1}}Q_{k_{1}-100>. \geq
r-100}\psi_{1}P_{k_{3}}Q_{<k_{1}-100}\psi_{3})||_{L_{t}^{2}L_{x}^{2}}\\
&\leq C ||\nabla^{-1}P_{k_{1}}Q_{k_{1}-100>.\geq
r-100}\psi_{1}||_{L_{t}^{2}L_{x}^{\infty}}||P_{k_{3}}Q_{<k_{1}-100}\psi_{3}||_{L_{t}^{\infty}L_{x}^{2}}\\
&\leq
C2^{-\frac{r}{2}}||P_{k_{1}}\psi_{1}||_{S[k_{1}]}||P_{k_{3}}\psi_{3}||_{S[k_{3}]}\\
\end{split}\nonumber\end{equation}

\begin{equation}\begin{split}
&||Q_{r}(\nabla^{-1}P_{k_{1}}Q_{<r-100}\psi_{1}P_{k_{3}}Q_{k_{1}-100>.
. \geq
r-100} \psi_{3})||_{L_{t}^{2}L_{x}^{2}}\\
&\leq
C||\nabla^{-1}P_{k_{1}}Q_{<r-100}\psi_{1}||_{L_{t}^{\infty}L_{x}^{\infty}}||P_{k_{3}}Q_{k_{1}-100>.\geq
r-100}\psi_{3}||_{L_{t}^{2}L_{x}^{2}}\\
&\leq 2^{-\frac{r}{2}}||P_{k_{1}}\psi_{1}||_{S[k_{1}]}||P_{k_{3}}\psi_{3}||_{S[k_{3}]}\\
\end{split}\nonumber\end{equation}

Now proceed as before to conclude

\begin{equation}\begin{split}
&||\sum_{r<k_{1}+O(1)}\Box
Q_{r}(\nabla^{-1}P_{k_{1}}Q_{<k_{1}-100}\psi_{1}P_{k_{3}}Q_{<k_{1}-100}\psi_{3})\nabla^{-1}P_{k_{2}}\psi_{2}||_{N[0]}\\
&\leq
C2^{\delta(k_{1}-k_{2})}\prod_{l=1}^{3}||P_{k_{l}}\psi_{l}||_{S[k_{l}]}\\
\end{split}\end{equation}

2.B):

\begin{equation}
\Box
(\nabla^{-1}P_{k_{1}}Q_{<k_{1}+100}\psi_{1})P_{k_{3}}Q_{<k_{1}-100}\psi_{3}\nabla^{-1}P_{k_{2}}\psi_{2}
\end{equation}

We can estimate this directly:

\begin{equation}\begin{split}
&||P_{0}[\Box
(\nabla^{-1}P_{k_{1}}Q_{<k_{1}+100}\psi_{1})P_{k_{3}}Q_{<k_{1}-100}\psi_{3}\nabla^{-1}P_{k_{2}}\psi_{2}]||_{N[0]}\\
&\leq||P_{0}[\Box
(\nabla^{-1}P_{k_{1}}Q_{<k_{1}+100}\psi_{1})P_{k_{3}}Q_{\geq
k_{1}-100}\psi_{3}\nabla^{-1}P_{k_{2}}\psi_{2}]||_{L_{t}^{1}L_{x}^{2}}\\
&+||P_{0}[\Box
(\nabla^{-1}P_{k_{1}}Q_{<k_{1}+100}\psi_{1})P_{k_{3}}\psi_{3}\nabla^{-1}P_{k_{2}}\psi_{2}]||_{L_{t}^{1}L_{x}^{2}}\\
&\leq C ||\Box
(\nabla^{-1}P_{k_{1}}Q_{<k_{1}+100}\psi_{1})||_{L_{t}^{2}L_{x}^{\infty}}||P_{k_{3}}Q_{\geq
k_{1}-100}\psi_{3}||_{L_{t}^{2}L_{x}^{2}}||\nabla^{-1}P_{k_{2}}\psi_{2}||_{L_{t}^{\infty}L_{x}^{\infty}}\\
&+C||\Box
(\nabla^{-1}P_{k_{1}}Q_{<k_{1}+100}\psi_{1})||_{L_{t}^{2}L_{x}^{\infty}}||P_{k_{3}}\psi_{3}||_{L_{t}^{4}L_{x}^{4}}||\nabla^{-1}P_{k_{2}}\psi_{2}||_{L_{t}^{4}L_{x}^{4}}\\
&\leq
C2^{\frac{k_{1}}{2}}\prod_{l=1}^{3}||P_{k_{l}}\psi_{l}||_{S[k_{l}]}\\
\end{split}\end{equation}

which is of course acceptable because of $k_{1}<k_{2}-10<O(1)$.
\\

2.C):

\begin{equation}
\nabla^{-1}P_{k_{1}}Q_{<k_{1}+100}\psi_{1}\Box(P_{k_{3}}Q_{<k_{1}-100}\psi_{3})\nabla^{-1}P_{k_{2}}\psi_{2}
\end{equation}

This is again straightforward because

\begin{equation}\begin{split}
&||P_{0}[\nabla^{-1}P_{k_{1}}Q_{<k_{1}+100}\psi_{1}\Box(P_{k_{3}}Q_{<k_{1}-100}\psi_{3})\nabla^{-1}P_{k_{2}}\psi_{2}]||_{L_{t}^{1}L_{x}^{2}}\\
&\leq
C||\nabla^{-1}P_{k_{1}}Q_{<k_{1}+100}\psi_{1}||_{L_{t}^{4}L_{x}^{\infty}}||\nabla^{-1}P_{k_{2}}\psi_{2}||_{L_{t}^{4}L_{x}^{\infty}}||\Box(P_{k_{3}}Q_{<k_{1}-100}\psi_{3})||_{L_{t}^{2}L_{x}^{2}}\\
&\leq C2^{\frac{k_{1}-k_{2}}{4}}\prod_{l=1}^{3}||P_{k_{l}}\psi_{l}||_{S[k_{l}]}\\
\end{split}\end{equation}

This finishes the proof of the lemma, and thereby the proof of
Case \begin{bf}1.A)\end{bf}.
\end{proof}

\newpage
\begin{bf}1.C): Estimation of C. \end{bf}
\\

This is the case corresponding to high-high interactions in $[,]$,
and can be rewritten as

\begin{equation}
\sum_{k\leq
\min\{k_{1}+O(1),-10\}}P_{0}(\sum_{j=1}^{3}\triangle^{-1}\partial_{j}P_{k}[R_{\nu}P_{k_{1}}\psi_{1}R_{j}P_{k_{2}}\psi_{2}-R_{j}P_{k_{1}}\psi_{1}R_{\nu}P_{k_{2}}\psi_{2}]P_{k_{3}}\partial^{\nu}\psi_{3})
\nonumber\end{equation}

We want to proceed in analogy to the case 1.A), by first reducing
output and input $P_{k_{3}}\psi_{3}$ to small modulation, in this
case modulation $<2^{k-100}$, where $k\leq \min\{k_{1}+O(1),-10\}$
is held fixed. Since we are eventually summing over $k$, we want
to obtain an exponential gain in the difference $k-k_{1}$. Keep in
mind that $k_{1}=k_{2}+O(1)$ for this case. We will use the
"imbedded $Q_{\nu j}$ null-form".
\\

\begin{bf}Output has large modulation\end{bf}
\\

1.1): Output has modulation $2^{l}$ with $k-100\leq l\leq
k_{1}+100$, $P_{k_{3}}\psi_{3}$ reduced to modulation
$<2^{k-100}$: We use here the simple identity $R_{\nu}f R_{j}g -
R_{j}f R_{\nu}g = \partial_{\nu}(\nabla^{-1}f
R_{j}g)-\partial_{j}(\nabla^{-1}f R_{\nu}g)$ in order to pull out
a derivative of $[,]$. This will allow us to play the modulation
of the output against the larger frequencies of the inputs of
$[,]$:

\begin{equation}\begin{split}
&||\sum_{k-100\leq l\leq
k_{1}+100}P_{0}Q_{l}(\sum_{j=1}^{3}P_{k}\triangle^{-1}\partial_{j}[R_{\nu}P_{k_{1}}\psi_{1}R_{j}P_{k_{2}}\psi_{2}-R_{j}P_{k_{1}}\psi_{1}R_{\nu}P_{k_{2}}\psi_{2}]\\&P_{k_{3}}Q_{<k-100}\psi_{3})||_{N[0]}\\
&\leq\sum_{k-100\leq l\leq k_{1}+100}C||P_{0}Q_{l}(P_{k}Q_{\leq
l+O(1)}\triangle^{-1}\nabla_{x,t}\partial_{j}[\nabla_{x,t}\nabla^{-1}P_{k_{1}}\psi_{1}\nabla^{-1}P_{k_{2}}\psi_{2}]\\&P_{k_{3}}Q_{<k-100}\partial^{\nu}\psi_{3})||_{\dot{X}_{0}^{-\frac{1}{2},-\frac{1}{2},1}}\\
&\leq\sum_{k-100\leq l\leq k_{1}+100}C
2^{-\frac{l}{2}}2^{\frac{k}{2}}2^{l}||\nabla_{x,t}\nabla^{-1}P_{k_{1}}\psi_{1}\nabla^{-1}P_{k_{2}}\psi_{2}||_{L_{t}^{2}L_{x}^{2}}||P_{k_{3}}Q_{<k-100}\partial^{\nu}\psi_{3}||_{L_{t}^{\infty}L_{x}^{2}}\\
&\leq C\sum_{k-100\leq l\leq
k_{1}+100}2^{\frac{k+l}{2}}2^{-k_{1}}\prod_{a=1}^{3}||P_{k_{a}}\psi_{a}||_{S[k_{a}]}\\
&\leq C2^{\frac{k-k_{1}}{2}}\prod_{a=1}^{3}||P_{k_{a}}\psi_{a}||_{S[k_{a}]}\\
\end{split}\nonumber\end{equation}

1.2): Output has modulation $2^{l}$ with $k-100\leq l\leq
k_{1}+100$, $P_{k_{3}}\psi_{3}$ reduced to modulation $\geq
2^{k-100}$:

\begin{equation}\begin{split}
&||\sum_{k-100\leq l\leq
k_{1}+100}P_{0}Q_{l}(\sum_{j=1}^{3}P_{k}\triangle^{-1}\partial_{j}[R_{\nu}P_{k_{1}}\psi_{1}R_{j}P_{k_{2}}\psi_{2}-R_{j}P_{k_{1}}\psi_{1}R_{\nu}P_{k_{2}}\psi_{2}]\\&P_{k_{3}}Q_{\geq k-100}\psi_{3})||_{\dot{X}_{0}^{-\frac{1}{2},-\frac{1}{2},1}}\\
&\leq\sum_{k-100\leq l\leq k_{1}+100}C
2^{-\frac{l}{2}}||P_{k}\triangle^{-1}\partial_{j}[R_{\nu}P_{k_{1}}\psi_{1}R_{j}P_{k_{2}}\psi_{2}-R_{j}P_{k_{1}}\psi_{1}R_{\nu}P_{k_{2}}\psi_{2}]||_{L_{t}^{\infty}L_{x}^{\infty}}\\&||P_{k_{3}}Q_{\geq
k-100}\psi_{3}||_{L_{t}^{2}L_{x}^{2}}\\
&\leq\sum_{k-100\leq l\leq
k_{1}+100}C2^{2k-\frac{l}{2}}||R_{\nu}P_{k_{1}}\psi_{1}R_{j}P_{k_{2}}\psi_{2}-R_{j}P_{k_{1}}\psi_{1}R_{\nu}P_{k_{2}}\psi_{2}||_{L_{t}^{\infty}L_{x}^{1}}\\&
||P_{k_{3}}Q_{\geq
k-100}\psi_{3}||_{L_{t}^{2}L_{x}^{2}}\\
&\leq \sum_{k-100\leq l\leq
k_{1}+100}C2^{\frac{3k-l}{2}-k_{1}}\prod_{a=1}^{3}||P_{k_{a}}\psi_{a}||_{S[k_{a}]}\\
&\leq C2^{k-k_{1}}\prod_{a=1}^{3}||P_{k_{a}}\psi_{a}||_{S[k_{a}]}\\
\end{split}\end{equation}

1.3): Output has modulation $2^{l}$ with $l>k_{1}+100$,
$P_{k_{3}}\psi_{3}$ reduced to modulation $<2^{l-100}$: this is
the case corresponding to very large modulation (by comparison
with the occuring frequencies) of the output. This condition then
entails that at least one input has at least comparable
modulation. Thus we can write

\begin{equation}\label{mess}\begin{split}
&\sum_{l>k_{1}+100}\sum_{j=1}^{3}P_{0}Q_{l}(\triangle^{-1}\partial_{j}P_{k}[R_{\nu}P_{k_{1}}\psi_{1}R_{j}P_{k_{2}}\psi_{2}-R_{j}P_{k_{1}}\psi_{1}R_{\nu}P_{k_{2}}\psi_{2}]P_{k_{3}}Q_{<l-100}\partial^{\nu}\psi_{3})\\
&=\sum_{l>k_{1}+100}\sum_{j=1}^{3}P_{0}Q_{l}(\triangle^{-1}\partial_{j}P_{k}[R_{\nu}P_{k_{1}}Q_{\geq
l-100} \psi_{1}R_{j}P_{k_{2}}\psi_{2}-R_{j}P_{k_{1}}Q_{\geq l-100}
\psi_{1}R_{\nu}P_{k_{2}}\psi_{2}]\\&P_{k_{3}}Q_{<l-100}\partial^{\nu}\psi_{3})\\
&+\sum_{l>k_{1}+100}\sum_{j=1}^{3}P_{0}Q_{l}(\triangle^{-1}\partial_{j}P_{k}[R_{\nu}P_{k_{1}}Q_{<
l-100} \psi_{1}R_{j}P_{k_{2}}Q_{\geq l-100}
\psi_{2}-R_{j}P_{k_{1}}Q_{<l-100}
\psi_{1}\\&R_{\nu}P_{k_{2}}Q_{\geq l-100}\psi_{2}]P_{k_{3}}Q_{<l-100}\partial^{\nu}\psi_{3})\\
\end{split}\end{equation}

For example, we can estimate

\begin{equation}\begin{split}
&\sum_{l>k_{1}+100}\sum_{j=1}^{3}||P_{0}Q_{l}[\triangle^{-1}\partial_{j}P_{k}(R_{\nu}P_{k_{1}}Q_{\geq
l-100}\psi_{1}R_{j}P_{k_{2}}\psi_{2})P_{k_{3}}Q_{<l-100}\partial^{\nu}\psi_{3}]_{N[0]}\\
&\leq\sum_{l>k_{1}+100}\sum_{j=1}^{3}C2^{-\frac{l}{2}}||\triangle^{-1}\partial_{j}P_{k}(R_{\nu}P_{k_{1}}Q_{\geq
l-100}\psi_{1}R_{j}P_{k_{2}}\psi_{2})||_{L_{t}^{2}L_{x}^{\infty}}\\&||P_{k_{3}}Q_{<l-100}\partial^{\nu}\psi_{3}||_{L_{t}^{\infty}L_{x}^{2}}\\
&\sum_{l>k_{1}+100}\sum_{j=1}^{3}2^{-\frac{l}{2}}2^{2k}||R_{\nu}P_{k_{1}}Q_{\geq
l-100}\psi_{1}||_{L_{t}^{2}L_{x}^{2}}||R_{j}P_{k_{2}}\psi_{2}||_{L_{t}^{\infty}L_{x}^{2}}\\&||P_{k_{3}}Q_{<l-100}\partial^{\nu}\psi_{3}||_{L_{t}^{\infty}L_{x}^{2}}\\
&\leq\sum_{l>k_{1}+100}C2^{2k-k_{1}-l}\prod_{a=1}^{3}||P_{k_{a}}\psi_{a}||_{S[k_{a}]}\\
&\leq C2^{2(k-k_{1})}\prod_{a=1}^{3}||P_{k_{a}}\psi_{a}||_{S[k_{a}]}\\
\end{split}\end{equation}

The other terms in \eqref{mess} are estimated similarly and
therefore left out.
\\

1.4): Output has modulation $2^{l}$ with $l>k_{1}+100$,
$P_{k_{3}}\psi_{3}$ reduced to modulation $\geq 2^{l-100}$:

\begin{equation}\begin{split}
&\sum_{l>k_{1}+100}||P_{0}Q_{l}[P_{k}\triangle^{-1}\partial_{j}(R_{\nu}P_{k_{1}}\psi_{1}R_{j}\psi_{2}-R_{j}P_{k_{1}}\psi_{1}R_{\nu}\psi_{2})P_{k_{3}}Q_{\geq
l-100}\partial^{\nu}\psi_{3}]||_{L_{t}^{1}L_{x}^{2}}\\
&\leq\sum_{l>k_{1}+100}||P_{k}\triangle^{-1}\partial_{j}(R_{\nu}P_{k_{1}}\psi_{1}R_{j}\psi_{2}-R_{j}P_{k_{1}}\psi_{1}R_{\nu}\psi_{2})||_{L_{t}^{2}L_{x}^{\infty}}||P_{k_{3}}Q_{\geq
l-100}\partial^{\nu}\psi_{3}||_{L_{t}^{2}L_{x}^{2}}\\
&\leq\sum_{l>k_{1}+100}C2^{\frac{k-l}{2}}\prod_{a=1}^{3}||P_{k_{a}}\psi_{a}||_{S[k_{a}]}\\
&\leq C2^{\frac{k-k_{1}}{2}}\prod_{a=1}^{3}||P_{k_{a}}\psi_{a}||_{S[k_{a}]}\\
\end{split}\end{equation}

\begin{bf}Output has small modulation\end{bf}
\\

2.1): Output has modulation $<2^{k-100}$, $P_{k_{3}}\psi_{3}$
reduced to modulation $\geq 2^{k+100}$: notice the identity

\begin{equation}\begin{split}
&P_{0}Q_{<k-100}[P_{k}\triangle^{-1}\partial_{j}(R_{\nu}P_{k_{1}}\psi_{1}R_{j}\psi_{2}-R_{j}P_{k_{1}}\psi_{1}R_{\nu}P_{k_{2}}\psi_{2})P_{k_{3}}Q_{\geq
k+100}\partial^{\nu}\psi_{3}]\\
&=\sum_{l\geq
k+100}P_{0}Q_{<k-100}[P_{k}Q_{l+O(1)}\triangle^{-1}\partial_{j}(R_{\nu}P_{k_{1}}\psi_{1}R_{j}\psi_{2}-R_{j}P_{k_{1}}\psi_{1}R_{\nu}P_{k_{2}}\psi_{2})\\&P_{k_{3}}Q_{l}\partial^{\nu}\psi_{3}]\\
\end{split}\nonumber\end{equation}

Next, observe that

\begin{equation}\begin{split}
&P_{k}Q_{l+O(1)}(P_{k_{1}}f P_{k_{2}}g) =\\
&P_{k}Q_{l+O(1)}(P_{k_{1}}P_{l+O(1)}Q_{<l-100}f)
P_{k_{2}}Q_{<l-100}g)\\&+P_{k}Q_{l+O(1)}(P_{k_{1}}Q_{\geq l-100}f
P_{k_{2}}g)\\
&+P_{k}Q_{l+O(1)}(P_{k_{1}}Q_{<l-100}f P_{k_{2}}Q_{\geq l-100}g)\\
\end{split}\nonumber\end{equation}

Thus we conclude that this case can be estimated by

\begin{equation}\begin{split}
&\sum_{l\geq k+100}C
2^{\frac{k}{2}}||R_{\nu}P_{l+O(1)}P_{k_{1}}Q_{<l-100}\psi_{1}||_{L_{t}^{4}L_{x}^{4}}||R_{j}P_{k_{2}}Q_{<l-100}\psi_{2}||_{L_{t}^{4}L_{x}^{4}}
||P_{k_{3}}Q_{l}\partial^{\nu}\psi_{3}||_{L_{t}^{2}L_{x}^{2}}\\
&+\sum_{l\geq k+100}C2^{2k}||R_{\nu}P_{k_{1}}Q_{\geq
l-100}\psi_{1}||_{L_{t}^{2}L_{x}^{2}}||R_{j}P_{k_{2}}\psi_{2}||_{L_{t}^{\infty}L_{x}^{2}}||P_{k_{3}}Q_{l}\psi_{3}||_{L_{t}^{2}L_{x}^{2}}\\
&+\sum_{l\geq
k+100}2^{2k}||R_{\nu}P_{k_{1}}Q_{<l-100}\psi_{1}||_{L_{t}^{\infty}L_{x}^{2}}||R_{j}P_{k_{2}}Q_{\geq
l-100}\psi_{2}||_{L_{t}^{2}L_{x}^{2}}||P_{k_{3}}Q_{l}\partial^{\nu}\psi_{3}||_{L_{t}^{2}L_{x}^{2}}\\
&\leq \sum_{l\geq
k+100,l=k_{1}+O(1)}C2^{\frac{k-l}{2}}\prod_{a=1}^{3}||P_{k_{a}}\psi_{a}||_{S[k_{a}]}+\sum_{l\geq
k+100}C2^{2k-l-k_{1}}\prod_{a=1}^{3}||P_{k_{a}}\psi_{a}||_{S[k_{a}]}\\
&\leq C2^{\frac{k-k_{1}}{2}}\prod_{a=1}^{3}||P_{k_{a}}\psi_{a}||_{S[k_{a}]}\\
\end{split}\nonumber\end{equation}

Having reduced output and $P_{k_{3}}\psi_{3}$ to small modulation,
we can now invoke the identity \eqref{special identity} and
proceed in exact analogy with case 1.A). Since this does not
entail any additional difficulties, it is left out. This then
finishes the low-high case \begin{bf}1)\end{bf}.
\\

\begin{bf}2): High-High Interactions: $[,]$ at frequency
$2^{k_{3}+O(1)}\geq 2^{-10}$, $k_{3}\geq -10$, i.e. $k_{3}\geq
O(1)$.\end{bf}
\\

We will again utilize the algebraic structure of the null-form,
but in a somewhat different fashion than before. The first step
will involve reduction to comparably low modulation of output and
third input. Having achieved this, we will try to pull out the
$\partial_{\nu}$ from the input of $[,]$ with larger frequency:
Assume w.l.o.g. in the sequel that the 1st input of $[,]$ has
larger frequency, i.e. $k_{2}\leq k_{1}$.

2.1): $P_{k_{3}}\partial^{\nu}\psi_{3}$ has modulation $\geq
2^{k_{3}}$:

\begin{equation}\begin{split}
&||P_{0}(\sum_{j=1}^{3}\triangle^{-1}\partial_{j}[R_{\nu}P_{k_{1}}\psi_{1}R_{j}P_{k_{2}}\psi_{2}-R_{j}P_{k_{1}}\psi_{1}R_{\nu}P_{k_{2}}\psi_{2}]\partial^{\nu}P_{k_{3}}Q_{\geq k_{3}}\psi_{3})||_{L_{t}^{1}L_{x}^{2}}\\
&\leq
C\sum_{j=1}^{3}2^{-k_{3}}||R_{\nu}P_{k_{1}}\psi_{1}R_{j}P_{k_{2}}\psi_{2}-R_{j}P_{k_{1}}\psi_{1}R_{\nu}P_{k_{2}}\psi_{2}||_{L_{t}^{2}L_{x}^{2}}||\partial^{\nu}P_{k_{3}}Q_{\geq
k_{3}}\psi_{3}||_{L_{t}^{2}L_{x}^{2}}\\
&\leq C2^{-k_{3}}2^{\frac{k_{2}-k_{1}}{6}}\prod_{a=1}^{3}||P_{k_{a}}\psi_{a}||_{S[k_{a}]}\\
\end{split}\end{equation}

2.2): Output at modulation $\geq 2^{0}$, $P_{k_{3}}\psi_{3}$ at
modulation $<2^{k_{3}}$:

\begin{equation}\begin{split}
&||P_{0}Q_{\geq 0}(\sum_{j=1}^{3}\triangle^{-1}\partial_{j}[R_{\nu}P_{k_{1}}\psi_{1}R_{j}P_{k_{2}}\psi_{2}-R_{j}P_{k_{1}}\psi_{1}R_{\nu}P_{k_{2}}\psi_{2}]\partial^{\nu}P_{k_{3}}Q_{<k_{3}}\psi_{3})||_{\dot{X}_{0}^{-\frac{1}{2},-\frac{1}{2},1}}\\
&\leq
C2^{-k_{3}}||R_{\nu}P_{k_{1}}\psi_{1}R_{j}P_{k_{2}}\psi_{2}-R_{j}P_{k_{1}}\psi_{1}R_{\nu}P_{k_{2}}\psi_{2}||_{L_{t}^{2}L_{x}^{2}}||\partial^{\nu}P_{k_{3}}Q_{<
k_{3}}\psi_{3}||_{L_{t}^{\infty}L_{x}^{2}}\\
&\leq
C2^{-\frac{k_{3}}{2}}2^{\frac{k_{2}-k_{1}}{6}}\prod_{a=1}^{3}||P_{k_{a}}\psi_{a}||_{S[k_{a}]}\\
\end{split}\end{equation}

2.3): Output at modulation $<2^{0}$, $P_{k_{3}}\psi_{3}$ at
modulation $< 2^{k_{3}}$: A simple algebraic manipulation reduces
this case to the sum

\begin{equation}\label{High-High}\begin{split}
&||P_{0}Q_{<0}(\sum_{j=1}^{3}\triangle^{-1}\partial_{j}\partial_{\nu}[R_{j}P_{k_{2}}\psi_{2}\nabla^{-1}P_{k_{1}}\psi_{1}]\partial^{\nu}P_{k_{3}}Q_{<k_{3}}\psi_{3})||_{N[0]}\\
&+||P_{0}Q_{<0}(\sum_{j=1}^{3}R_{\nu}P_{k_{2}}\psi_{2}\nabla^{-1}P_{k_{1}}\psi_{1}P_{k_{3}}Q_{<k_{3}}\partial^{\nu}\psi_{3})||_{N[0]}\\
\end{split}\end{equation}

For the first of these, we want to remove $\partial_{\nu}$ from
$[,]$, as this term might have very large modulation. We are in a
favorable situation since letting $\partial_{\nu}$ fall on the
output is harmless on account of our assumptions, while letting it
fall on the third term $P_{k_{3}}Q_{<k_{3}}\partial^{\nu}\psi_{3}$
is quite useful, as it produces a $\Box$-operator. Thus we
majorize the first term of \eqref{High-High} by

\begin{equation}\begin{split}
&\sum_{j=1}^{3}||P_{0}Q_{<0}\partial_{\nu}(\triangle^{-1}\partial_{j}[R_{j}P_{k_{2}}\psi_{2}\nabla^{-1}P_{k_{1}}\psi_{1}]P_{k_{3}}Q_{<k_{3}}\partial^{\nu}\psi_{3})||_{N[0]}\\
&+\sum_{j=1}^{3}||P_{0}Q_{<0}(\triangle^{-1}\partial_{j}[R_{j}P_{k_{2}}\psi_{2}\nabla^{-1}P_{k_{1}}\psi_{1}]\Box
P_{k_{3}}Q_{<k_{3}}\psi_{3})||_{N[0]}\\
\end{split}\nonumber\end{equation}

For the first term in the immediately preceding, we have

\begin{equation}\begin{split}
&\sum_{j=1}^{3}||P_{0}Q_{<0}\partial_{\nu}(\triangle^{-1}\partial_{j}[R_{j}P_{k_{2}}\psi_{2}\nabla^{-1}P_{k_{1}}\psi_{1}]P_{k_{3}}Q_{<k_{3}}\partial^{\nu}\psi_{3})||_{L_{t}^{1}L_{x}^{2}}\\
&\leq
C2^{-k_{3}-k_{1}}||P_{k_{2}}\psi_{2}||_{L_{t}^{3}L_{x}^{6}}||P_{k_{1}}\psi_{1}||_{L_{t}^{3}L_{x}^{6}}||P_{k_{3}}Q_{<k_{3}}\partial^{\nu}\psi_{3}||_{L_{t}^{3}L_{x}^{6}}\\
&\leq C 2^{-k_{1}}2^{\frac{k_{1}+k_{2}+k_{3}}{6}}\prod_{a=1}^{3}||P_{k_{a}}\psi_{a}||_{S[k_{a}]}\\
&\leq C2^{-\frac{k_{3}}{2}}2^{\frac{k_{2}-k_{1}}{6}}\prod_{a=1}^{3}||P_{k_{a}}\psi_{a}||_{S[k_{a}]}\\
\end{split}\end{equation}

We have used here that in the present situation $k_{1}\geq
k_{3}+O(1)$.
\\

For the 2nd term, we have

\begin{equation}\begin{split}
&||P_{0}Q_{<0}(\triangle^{-1}\partial_{j}[R_{j}P_{k_{2}}\psi_{2}\nabla^{-1}P_{k_{1}}\psi_{1}]\Box
P_{k_{3}}Q_{<k_{3}}\psi_{3})||_{L_{t}^{1}L_{x}^{2}}\\
&\leq
C2^{-k_{1}-k_{3}}||R_{j}P_{k_{2}}\psi_{2}||_{L_{t}^{3}L_{x}^{6}}||P_{k_{1}}\psi_{1}||_{L_{t}^{6}L_{x}^{3}}||\Box
P_{k_{3}}Q_{<k_{3}}\psi_{3}||_{L_{t}^{2}L_{x}^{2}}\\
&\leq C2^{-k_{1}}2^{\frac{k_{2}-k_{1}}{6}}\prod_{a=1}^{3}||P_{k_{a}}\psi_{a}||_{S[k_{a}]}\\
&\leq C 2^{-k_{3}}2^{\frac{k_{2}-k_{1}}{6}}\prod_{a=1}^{3}||P_{k_{a}}\psi_{a}||_{S[k_{a}]}\\
\end{split}\end{equation}

Hence in order to finish case \begin{bf}2)\end{bf}, we need to
deal with the 2nd term in \eqref{High-High}. But this is
immediate, referring to lemma~\ref{almost endpoint} in the
preceding section, as well as the disposability of $P_{0}Q_{<0}$,
$P_{k_{3}}Q_{<k_{3}}$.
\\

\begin{bf}3): High-Low Interactions: $[,]$ at frequency $\geq 2^{-10}$,
$k_{3}<-10$.\end{bf}
\\

This case is the most elementary on account of the fact that a low
frequency term is hit by a derivative. Moreover, it can be dealt
with by the same methods as in the immediately preceding case, so
we shall discuss it only briefly. First, one reduces the output to
modulation $<2^{0}$ and the third input $P_{k_{3}}\psi_{3}$ to
modulation $<2^{k_{3}}$. Then, as we only have to deal with the
case when the inputs $\psi_{1},\psi_{2}$ are at frequencies
$2^{k_{1}},\,2^{k_{2}}$ with $k_{1} = k_{2}+O(1) >>O(1)$
(otherwise, evaluate all inputs in $L_{t}^{3}L_{x}^{6}$), we use
as before the identity

\begin{equation}\begin{split}
&\sum_{j=1}^{3}P_{0}Q_{<0}(\triangle^{-1}\partial_{j}[R_{\nu}P_{k_{1}}\psi_{1}R_{j}P_{k_{2}}\psi_{2}-R_{j}P_{k_{1}}\psi_{1}R_{\nu}P_{k_{2}}\psi_{2}]P_{k_{3}}Q_{<k_{3}}\partial^{\nu}\psi_{3})\\
& =
\sum_{j=1}^{3}P_{0}Q_{<0}\partial_{\nu}(\triangle^{-1}\partial_{j}[\nabla^{-1}P_{k_{1}}\psi_{1}R_{j}P_{k_{2}}\psi_{2}]\partial^{\nu}P_{k_{3}}Q_{<k_{3}}\psi_{3})\\
&-\sum_{j=1}^{3}P_{0}Q_{<0}(\triangle^{-1}\partial_{j}[\nabla^{-1}P_{k_{1}}\psi_{1}R_{j}P_{k_{2}}\psi_{2}]\Box
P_{k_{3}}Q_{<k_{3}}\psi_{3})\\
&-P_{0}Q_{<0}(\nabla^{-1}P_{k_{1}}\psi_{1}R_{\nu}P_{k_{2}}\psi_{2})\partial^{\nu}P_{k_{3}}Q_{<k_{3}}\psi_{3}\\
\end{split}\end{equation}

Each of these can be easily estimated:

\begin{equation}\begin{split}
&||\sum_{j=1}^{3}P_{0}Q_{<0}\partial_{\nu}(\triangle^{-1}\partial_{j}[\nabla^{-1}P_{k_{1}}\psi_{1}R_{j}P_{k_{2}}\psi_{2}]\partial^{\nu}P_{k_{3}}Q_{<k_{3}}\psi_{3})||_{N[0]}\\
&\leq\sum_{j=1}^{3}||P_{0}Q_{<0}\partial_{\nu}(\triangle^{-1}\partial_{j}[\nabla^{-1}P_{k_{1}}\psi_{1}R_{j}P_{k_{2}}\psi_{2}]\partial^{\nu}P_{k_{3}}Q_{<k_{3}}\psi_{3})||_{L_{t}^{1}L_{x}^{2}}\\
&\leq
\sum_{j=1}^{3}C||\nabla^{-1}P_{k_{1}}\psi_{1}||_{L_{t}^{3}L_{x}^{6}}||R_{j}P_{k_{2}}\psi_{2}||_{L_{t}^{3}L_{x}^{6}}||\partial^{\nu}P_{k_{3}}Q_{<k_{3}}\psi_{3}||_{L_{t}^{3}L_{x}^{6}}\\
&\leq C
2^{-\frac{2k_{1}}{3}}2^{-\frac{7|k_{3}|}{6}}\prod_{a=1}^{3}||P_{k_{a}}\psi_{a}||_{S[k_{a}]}\\
\end{split}\end{equation}

\begin{equation}\begin{split}
&||\sum_{j=1}^{3}P_{0}Q_{<0}(\triangle^{-1}\partial_{j}[\nabla^{-1}P_{k_{1}}\psi_{1}R_{j}P_{k_{2}}\psi_{2}]\Box
P_{k_{3}}Q_{<k_{3}}\psi_{3})||_{N[0]}\\
&\leq\sum_{j=1}^{3}||P_{0}Q_{<0}(\triangle^{-1}\partial_{j}[\nabla^{-1}P_{k_{1}}\psi_{1}R_{j}P_{k_{2}}\psi_{2}]\Box
P_{k_{3}}Q_{<k_{3}}\psi_{3})||_{L_{t}^{1}L_{x}^{2}}\\
&\leq\sum_{j=1}^{3}C||\nabla^{-1}P_{k_{1}}\psi_{1}R_{j}P_{k_{2}}\psi_{2}||_{L_{t}^{2}L_{x}^{2}}||\Box
P_{k_{3}}Q_{<k_{3}}\psi_{3}||_{L_{t}^{2}L_{x}^{\infty}}\\
&\leq C
2^{-k_{1}}2^{-\frac{5|k_{3}|}{2}}\prod_{a=1}^{3}||P_{k_{a}}\psi_{a}||_{S[k_{a}]}\\
\end{split}\end{equation}

\begin{equation}\begin{split}
&||P_{0}Q_{<0}(\nabla^{-1}P_{k_{1}}\psi_{1}R_{\nu}P_{k_{2}}\psi_{2})\partial^{\nu}P_{k_{3}}Q_{<k_{3}}\psi_{3}||_{L_{t}^{1}L_{x}^{2}}\\
&\leq C2^{-\frac{2k_{1}}{3}}2^{-\frac{7|k_{3}|}{6}}\prod_{a=1}^{3}||P_{k_{a}}\psi_{a}||_{S[k_{a}]}\\
\end{split}\end{equation}

Of course, in order to verify the statement of the Proposition, we
can even discard the explicit gain in $k_{1}$. This finishes the
proof of the Proposition.
\end{proof}

\subsection{The Gauge Change Estimate}

\begin{proposition}
Let $f(x)$ be a smooth function with all derivatives up to and
including fourth order bounded. Then, provided
$\tilde{\phi}_{j},\psi$ are Schwartz functions with
$\max\{||P_{k}\tilde{\phi}_{j}||_{S[k]},
||P_{k}\psi||_{S[k]}\}\leq c_{k}$ for a frequency envelope $c_{k}$
as in the previous section, we have

\begin{equation}
||P_{k}(f(\triangle^{-1}\sum_{j=1}^{3}\partial_{j}\tilde{\phi}_{j})\psi)||_{S[k]}\leq
Cc_{k}
\end{equation}
\end{proposition}

\begin{proof}
The proof of this assertion will consist in the careful analysis
of many different cases. In particular, the following observation
will be used many times: Let $u\in C^{\infty}(\mathbf{R}^{n})$,
and $F\in C^{\infty}(\mathbf{R})$.

\begin{lemma}\label{basic idea}
Let $X$ be a translation invariant norm defined on all measurable
functions, and $F$, $u$ as before. Then

\begin{equation}
||P_{0}(F(u))||_{X}\leq C ||P_{0}(\nabla u F'(u))||_{X}
\end{equation}
\end{lemma}

The proof of this is immediate: We have

\begin{equation}\begin{split}
& P_{0}(F(u))= P_{0}\tilde{P}_{0}(F(u))
=P_{0}(\sum_{j=1}^{n}\triangle^{-1}\partial_{j}^{2}\tilde{P_{0}}(F(u)))\\&=
P_{0}\sum_{j}\int_{\mathbf{R}^{n}}\triangle^{-1}\partial_{j}\hat{m_{0}}(x-y)\partial_{j}(F(u))(y)dy\\
&=P_{0}\sum_{j}\int_{\mathbf{R}^{n}}a_{j}(x-y)\partial_{j}(F(u))(y)dy\\
\end{split}\end{equation}

where $a_{j}(y) = \triangle^{-1}\partial_{j}(\hat{m_{0}})(y)$,
$m_{0}$ is the symbol of $\tilde{P}_{0}$, and $\tilde{P}_{0}$ is a
Fourier multiplier like $P_{0}$ whose symbol equals $1$ on the
Fourier support of $P_{0}$. Since $a(y)$ has finite $L^{1}$-norm,
the claim follows.
\\

The reason why the previous lemma might be useful is that by
hitting the function $f$ in the Proposition with a derivative, we
gain an extra factor which is morally equivalent to $\phi$. The
more such factors are present, the more freedom is gained in
proving the necessary estimates. In particular, we must and can
avoid to expand the function into a power series(which would
require real analyticity anyways), for then the crucial
$L^{\infty}$ bound would be lost.
\\

The following lemma will also be used many times in the sequel:

\begin{lemma}\label{two-two}
Assume that $||P_{k}\psi_{j}||_{S[k]}\leq c_{k}$, with $c_{k}$ a
frequency envelope as before,
$\psi_{j}\in\mathcal{S}(\mathbf{R}^{3+1})$, $j\in{1,2,3}$.
Provided that $k\leq j+O(1)$, we have

\begin{equation}
||P_{k}Q_{j}(f(\sum_{j=1}^{3}\triangle^{-1}\partial_{j}\psi_{j}))||_{L_{t}^{2}L_{x}^{2}}\leq
C 2^{-\frac{3j}{2}-\frac{k}{2}}
\end{equation}
\end{lemma}

\begin{proof}
Replace $\sum_{j=1}^{3}\triangle^{-1}\partial_{j}\psi$ formally
for simplicity's sake by $\nabla^{-1}\psi$. For the sake of
simplicity, assume that $|O(1)|$ in the formulation of the theorem
and throughout the rest of the paper is $\leq 50$. We split the
expression as follows:
\begin{equation}\label{26}\begin{split}
&P_{k}Q_{j}(f(\nabla^{-1}\psi))= \Box^{-1}P_{k}Q_{j}(\Box\nabla^{-1}\psi f'(\nabla^{-1}\psi))\\&+\Box^{-1}P_{k}Q_{j}(\partial_{\nu}\nabla^{-1}\psi\partial^{\nu}\nabla^{-1}\psi f''(\nabla^{-1}\psi))\\
&=\Box^{-1}P_{k}Q_{j}(P_{k+O(1)}\Box\nabla^{-1}\psi P_{<k-100}(f'(\nabla^{-1}\psi)))\\
&+\Box^{-1}P_{k}Q_{j}(P_{<k-100}\Box\nabla^{-1}\psi P_{k+O(1)}(f'(\nabla^{-1}\psi))\\
&+\sum_{l_{1}=l_{2}+O(1)\geq
k+O(1)}\Box^{-1}P_{k}Q_{j}(P_{l_{1}}\Box\nabla^{-1}\psi P_{l_{2}}(f'(\nabla^{-1}\psi))\\
&+\Box^{-1}P_{k}Q_{j}(\partial_{\nu}\nabla^{-1}\psi\partial^{\nu}\nabla^{-1}\psi f''(\nabla^{-1}\psi))\\
\end{split}\end{equation}
where $\Box^{-1}P_{k}Q_{j}$ is a Fourier multiplier with symbol
$\frac{m_{k}(|\xi|)m_{j}(||\tau|-|\xi||)}{|\tau|^{2}-|\xi|^{2}}$.
This operator is actually disposable on account of $j\geq k+O(1)$,
but its boundedness on $L_{t}^{2}L_{x}^{2}$ is entirely
elementary. Indeed, its operator norm is dominated by $C2^{-2j}$.
The first 3 terms above represent the usual trichotomy into
high-low, low-high,
high-high interactions. We treat these first, the last term being most elementary.  \\

1): High-Low: Split this as follows:
\begin{equation}\label{27}\begin{split}
&\Box^{-1}P_{k}Q_{j}(P_{k+O(1)}\Box\nabla^{-1}\psi P_{<k-100}(f'(\nabla^{-1}\psi)))\\
=&\Box^{-1}P_{k}Q_{j}(P_{k+O(1)}Q_{<j+100}\Box\nabla^{-1}\psi
P_{<k-100}(f'(\nabla^{-1}\psi)))\\&+\Box^{-1}P_{k}Q_{j}(P_{k+O(1)}Q_{\geq
j+100}
\Box\nabla^{-1}\psi P_{<k-100}(f'(\nabla^{-1}\psi)))\\
\end{split}\end{equation}
The $L_{t}^{2}L_{x}^{2}$-norm of the first term can be estimated
by $C 2^{\frac{j-k}{2}}2^{-2j}c_{k}\leq C c_{k}
2^{-\frac{3}{2}j-\frac{k}{2}}$ by placing
$P_{k+O(1)}Q_{<j-100}\Box\nabla^{-1}\psi$ into
$L_{t}^{2}L_{x}^{2}$. \\
As to the 2nd term, write it as
\begin{equation}
\Box^{-1}P_{k}Q_{j}\sum_{l\geq
j+100}P_{k+O(1)}Q_{l}\Box\nabla^{-1}\psi
P_{<k-100}Q_{l+O(1)}(f'(\nabla^{-1}\psi))
\end{equation}
where we have used the fact that if $(\tau_{1},\xi_{1})$ denotes a
point in the Fourier support of $\Box\nabla^{-1}\psi$ and
$(\tau_{2},\xi_{2})$ denotes a point in the Fourier support of
$P_{<k-100}f'(\nabla^{-1}\psi)$, while $(\tau,\xi)$ denotes as
point in the Fourier support of the output, we have the condition
\begin{equation}
2^{j}\sim ||\tau|-|\xi||=|\tau_{1}\pm |\xi_{1}|+\tau_{2}\pm
|\xi_{2}|\mp |\xi_{1}|\mp|\xi_{2}|\pm |\xi_{1}+\xi_{2}||
\end{equation}
where the signs have been chosen in such a way that $|\tau_{a}\pm
|\xi_{a}||= ||\tau_{a}|-|\xi_{a}||,|\tau_{1}+\tau_{2}\pm
|\xi_{1}+\xi_{2}||=|\tau_{1}+\tau_{2}|-|\xi_{1}+\xi_{2}||,
a\in{1,2}$. This forces the factor $P_{<k-100}f'(\nabla^{-1}\psi)$
to be microlocalized at the indicated modulation. Now use the
$L_{t}^{4}L_{x}^{4} $-boundedness of the operator
$\nabla_{t}^{-1}P_{<k}Q_{l+O(1)}$ with symbol $\frac{m_{\leq
k}(|\xi|)m_{l+O(1)}(||\tau|-|\xi||)}{\tau}$; indeed, this operator
is given by convolution with a kernel whose $L^{1}$-norm is $\leq
C2^{-l}$. This allows us to estimate this term by
\begin{equation}\begin{split}
&||\sum_{l\geq
j+100}P_{k+O(1)}Q_{l}\Box\nabla^{-1}\psi P_{<k-100}Q_{l+O(1)}(f'(\nabla^{-1}\psi))||_{L_{t}^{2}L_{x}^{2}}\\
&\leq ||\sum_{l\geq
j+100}P_{k+O(1)}Q_{l}\Box\nabla^{-1}\psi\nabla_{t}^{-1}P_{<k-100}Q_{l+O(1)}(\partial_{t}\nabla^{-1}\psi f''(\nabla^{-1}\psi))||_{L_{t}^{2}L_{x}^{2}}\\
&\leq\sum_{l\geq
j+100}||P_{k+O(1)}Q_{l}\Box\nabla^{-1}\psi||_{L_{t}^{4}L_{x}^{4}}||\nabla_{t}^{-1}P_{<k-100}Q_{l+O(1)}(\nabla^{-1}\partial_{t}\psi f''(\nabla^{-1}\psi))||_{L_{t}^{4}L_{x}^{4}}\\
&\leq\sum_{l\geq j+100}C2^{l}c_{k}2^{-l}2^{\frac{k-l}{4}}\leq C c_{k}2^{\frac{k-j}{4}}\\
\end{split}\nonumber\end{equation}
Hence the $L_{t}^{2}L_{x}^{2}$-norm of the 2nd term in \eqref{27}
can be estimated by an expression
of the desired form.\\

2): Low-High: Split this term as follows:
\begin{equation}\begin{split}
&\Box^{-1}P_{k}Q_{j}(P_{<k-100}\Box\nabla^{-1}\psi P_{k+O(1)}(f'(\nabla^{-1}\psi))\\
&=\Box^{-1}P_{k}Q_{j}(P_{<k-100}Q_{<j+100}\Box\nabla^{-1}\psi
P_{k+O(1)}(f'(\nabla^{-1}\psi))\\&+\Box^{-1}P_{k}Q_{j}(P_{<k-100}Q_{\geq
j+100}\Box\nabla^{-1}\psi P_{k+O(1)}(f'(\nabla^{-1}\psi)))\\
\end{split}\nonumber\end{equation}
We have
\begin{equation}\begin{split}
&\Box^{-1}P_{k}Q_{j}(P_{<k-100}Q_{<j+100}\Box\nabla^{-1}\psi P_{k+O(1)}(f'(\nabla^{-1}\psi)))\\
&=\Box^{-1}P_{k}Q_{j}(P_{<k-100}Q_{<j+100}\Box\nabla^{-1}\psi\sum_{a=1}^{3}P_{k+O(1)}\triangle^{-1}\partial_{a}(\partial_{a}\nabla^{-1}\psi f''(\nabla^{-1}\psi)))\\
\end{split}\nonumber\end{equation}

whence
\begin{equation}\begin{split}
&||\Box^{-1}P_{k}Q_{j}(P_{<k-100}Q_{<j+100}\Box\nabla^{-1}\psi P_{k+O(1)}(f'(\nabla^{-1}\psi)))||_{L_{t}^{2}L_{x}^{2}}\\
&\leq
2^{-2j}||P_{<k-100}Q_{<j+100}\Box\nabla^{-1}\psi||_{L_{t}^{2}L_{x}^{6}}||\sum_{a=1}^{3}P_{k+O(1)}\triangle^{-1}\partial_{a}(\partial_{a}\nabla^{-1}\psi f''(\nabla^{-1}\psi))||_{L_{t}^{\infty}L_{x}^{3}}\\
& \leq 2^{-2j}C 2^{\frac{k}{2}}2^{\frac{j}{2}}2^{-k}\leq C
2^{-\frac{3j}{2}-\frac{k}{2}}\\
\end{split}\nonumber\end{equation}
Also,
\begin{equation}\begin{split}
&\sum_{l\geq j+100}\Box^{-1}P_{k}Q_{j}(P_{<k-100}Q_{l}\Box\nabla^{-1}\psi P_{k+O(1)}Q_{l+O(1)}(f'(\nabla^{-1}\psi)))\\
&=\sum_{l\geq j+100}\Box^{-1}P_{k}Q_{j}(P_{<k-100}Q_{l}\Box\nabla^{-1}\psi\nabla_{t}^{-1}P_{k+O(1)}Q_{l+O(1)}(\partial_{t}\nabla^{-1}\psi f''(\nabla^{-1}\psi)))\\
\end{split}\nonumber\end{equation}
Therefore
\begin{equation}\begin{split}
&\sum_{l\geq j+100}||\Box^{-1}P_{k}Q_{j}(P_{<k-100}Q_{l}\Box\nabla^{-1}\psi P_{k+O(1)}Q_{l+O(1)}(f'(\nabla^{-1}\psi)))||_{L_{t}^{2}L_{x}^{2}}\\
&\leq 2^{-2j}\sum_{l\geq
j+100}\sum_{r<k-100}||P_{r}Q_{l}\Box\nabla^{-1}\psi||_{L_{t}^{4}L_{x}^{4}}2^{-l}||\partial_{t}\nabla^{-1}\psi f''(\nabla^{-1}\psi)||_{L_{t}^{4}L_{x}^{4}}\\
&\leq 2^{-2j}\sum_{l\geq j+100}\sum_{r<k-100}C
2^{\frac{r-l}{4}} c_{r}\leq 2^{-2j}C\\
\end{split}\nonumber\end{equation}

3): High-High interactions: We have
\begin{equation}\label{32}\begin{split}
&\sum_{l_{1}=l_{2}+O(1)\geq
k+O(1)}\Box^{-1}P_{k}Q_{j}(P_{l_{1}}\Box\nabla^{-1}\psi P_{l_{2}}(f'(\nabla^{-1}\psi))\\
&=\sum_{l_{1}=l_{2}+O(1)\geq
k+O(1)}\Box^{-1}P_{k}Q_{j}(P_{l_{1}}Q_{\leq l_{1}+10}(\Box\nabla^{-1}\psi)\sum_{a=1}^{3}P_{l_{2}}\triangle^{-1}\partial_{a}(\partial_{a}\nabla^{-1}\psi f''(\nabla^{-1}\psi)))\\
&+\sum_{l_{1}=l_{2}+O(1)\geq
k+O(1)}\Box^{-1}P_{k}Q_{j}(P_{l_{1}}Q_{>l_{1}+10}(\Box\nabla^{-1}\psi)P_{l_{2}}(f'(\nabla^{-1}\psi)))\\
\end{split}\end{equation}
But $||P_{l_{1}}Q_{\leq
l_{1}+10}(\Box\nabla^{-1}\psi)||_{L_{t}^{2}L_{x}^{2}}\leq C\mu
c_{l_{1}}$, and $||\partial_{a}\nabla^{-1}\psi
f''(\nabla^{-1}\psi)||_{L_{t}^{\infty}L_{x}^{3}}\leq C$ by the
fundamental theorem of Littlewood-Paley theory; also,
$P_{l_{2}}\triangle^{-1}\partial_{a}$ is given by convolution with
a kernel whose $L^{1}$-norm is bounded by $C2^{-l_{2}}$, so we can
estimate
\begin{equation}\begin{split}
&||\sum_{l_{1}=l_{2}+O(1)\geq
k+O(1)}\Box^{-1}P_{k}Q_{j}(P_{l_{1}}Q_{\leq l_{1}+10}(\Box\nabla^{-1}\psi)\sum_{a=1}^{3}P_{l_{2}}\triangle^{-1}\partial_{a}(\partial_{a}\nabla^{-1}\psi f'(\nabla^{-1}\psi)))||_{L_{t}^{2}L_{x}^{2}}\\
&\leq 2^{-2j}\sum_{l_{1}=l_{2}+O(1)\geq k+O(1)}2^{k}C
c_{l_{1}}2^{-l_{2}}\leq 2^{-2j}C
\end{split}\nonumber\end{equation}
where in the last step we have used Bernstein's inequality.\\
As to the 2nd term in \eqref{32}, we have
\begin{equation}\label{horror}\begin{split}
&\sum_{l_{1}=l_{2}+O(1)\geq
k+O(1)}\Box^{-1}P_{k}Q_{j}(P_{l_{1}}Q_{>l_{1}+10}(\Box\nabla^{-1}\psi)P_{l_{2}}(f'(\nabla^{-1}\psi)))\\
&=\sum_{l_{1}=l_{2}+O(1)\geq
k+O(1)}\Box^{-1}P_{k}Q_{j}(P_{l_{1}}Q_{j+100\geq.>l_{1}+10}(\Box\nabla^{-1}\psi)P_{l_{2}}(f'(\nabla^{-1}\psi)))\\
&+\sum_{l_{1}=l_{2}+O(1)\geq
k+O(1)}\Box^{-1}P_{k}Q_{j}(P_{l_{1}}Q_{>\max\{j+100,l_{1}+10\}}(\Box\nabla^{-1}\psi)P_{l_{2}}(f'(\nabla^{-1}\psi)))\\
\end{split}\end{equation}
The 2nd term in the immediately preceding is estimated as follows:
\begin{equation}\begin{split}
&||\sum_{l_{1}=l_{2}+O(1)\geq
k+O(1)}\Box^{-1}P_{k}Q_{j}(P_{l_{1}}Q_{j+100\geq.>l_{1}+10}(\Box\nabla^{-1}\psi)P_{l_{2}}(f'(\nabla^{-1}\psi)))||_{L_{t}^{2}L_{x}^{2}}\\
&\leq \sum_{j+100\geq r>l_{1}+10}\sum_{l_{1}=l_{2}+O(1)\geq
k+O(1)}||\Box^{-1}P_{k}Q_{j}(P_{l_{1}}Q_{r}(\Box\nabla^{-1}\psi)P_{l_{2}}(f'(\nabla^{-1}\psi)))||_{L_{t}^{2}L_{x}^{2}}\\
&\leq
C2^{-2j}||P_{l_{1}}Q_{r}(\Box\nabla^{-1}\psi)||_{L_{t}^{2}L_{x}^{2}}||P_{l_{2}}(f'(\nabla^{-1}\psi))||_{L_{t}^{\infty}L_{x}^{\infty}}\\
&\leq C2^{-2j}\sum_{j+100\geq
r>l_{1}+10}\sum_{l_{1}=l_{2}+O(1)\geq
k+O(1)}2^{\frac{r-l_{1}}{2}}c_{l_{1}}\leq
2^{-\frac{3j}{2}-\frac{k}{2}}C
\end{split}\nonumber\end{equation}
As to the 3rd term in \eqref{horror}, we have
\begin{equation}\begin{split}
&\sum_{r>\max\{j+100,l_{1}+10\}}\Box^{-1}P_{k}Q_{j}(Q_{r}\Box\nabla^{-1}\psi P_{l_{2}}(f'(\nabla^{-1}\psi)))\\
&=\sum_{r>\max\{j+100,l_{1}+10\}}\Box^{-1}P_{k}Q_{j}(Q_{r}\Box\nabla^{-1}\psi P_{l_{2}}Q_{r+O(1)}(f'(\nabla^{-1}\psi)))\\
\end{split}\nonumber\end{equation}
Arguing as in the corresponding high-modulation-inputs case for
the high-low or low-high interactions, we have
\begin{equation}\begin{split}
&\sum_{r>\max\{j+100,l_{1}+10\}}\Box^{-1}P_{k}Q_{j}(P_{l_{1}}Q_{r}\Box\nabla^{-1}\psi P_{l_{2}}Q_{r+O(1)}(f'(\nabla^{-1}\psi)))\\
&=\sum_{r>\max\{j+100,l_{1}+10\}}\Box^{-1}P_{k}Q_{j}(P_{l_{1}}Q_{r}\Box\nabla^{-1}\psi\nabla_{t}^{-1}P_{l_{2}}Q_{r+O(1)}(\partial_{t}\nabla^{-1}\psi f''(\nabla^{-1}\psi)))\\
&=\sum_{r>\max\{j+100,l_{1}+10\}}\\&\Box^{-1}P_{k}Q_{j}(P_{l_{1}}Q_{r}\Box\nabla^{-1}\psi\nabla_{t}^{-1}P_{l_{2}}Q_{r+O(1)}(P_{l_{2}+O(1)}\partial_{t}\nabla^{-1}\psi P_{<l_{2}-100}f''(\nabla^{-1}\psi)))\\
&+\sum_{r>\max\{j+100,l_{1}+10\}}\\&\Box^{-1}P_{k}Q_{j}(P_{l_{1}}Q_{r}\Box\nabla^{-1}\psi\nabla_{t}^{-1}P_{l_{2}}Q_{r+O(1)}(\partial_{t}\nabla^{-1}\psi\sum_{a=1}^{3}P_{\geq l_{2}-100}\triangle^{-1}\partial_{a}(\partial_{a}\nabla^{-1}\psi f'''(\nabla^{-1}\psi)))\\
\end{split}\nonumber\end{equation}
Now
\begin{equation}\begin{split}
&\sum_{l_{1}=l_{2}+O(1)}\sum_{r>\max\{j+100,l_{1}+10\}}||\Box^{-1}P_{k}Q_{j}(P_{l_{1}}Q_{r}\Box\nabla^{-1}\psi\nabla_{t}^{-1}P_{l_{2}}Q_{r+O(1)}(P_{l_{2}+O(1)}\partial_{t}\nabla^{-1}\psi\\&P_{<l_{2}-100}f''(\nabla^{-1}\psi)))||_{L_{t}^{2}L_{x}^{2}}\\
&\leq 2^{-2j}\sum_{l_{1}=l_{2}+O(1)}\sum_{r>\max\{j+100,l_{1}+10\}}||P_{l_{1}}Q_{r}\Box\nabla^{-1}\psi||_{L_{t}^{4}L_{x}^{4}}2^{-r}||P_{l_{2}+O(1)}\partial_{t}\nabla^{-1}\psi\\&P_{<l_{2}-100}f''(\nabla^{-1}\psi)||_{L_{t}^{4}L_{x}^{4}}\\
&\leq
2^{-2j}\sum_{l_{1}=l_{2}+O(1)}\sum_{r>\max\{j+100,l_{1}+10\}}2^{r}2^{\frac{l_{1}-r}{4}}c_{l_{1}}2^{-r}c_{l_{2}}\leq
2^{-2j}C\\
\end{split}\nonumber\end{equation}
Moreover
\begin{equation}\begin{split}
&\sum_{l_{1}=l_{2}+O(1)}\sum_{r>\max\{j+100,l_{1}+10\}}||\Box^{-1}P_{k}Q_{j}(P_{l_{1}}Q_{r}\Box\nabla^{-1}\psi\nabla_{t}^{-1}P_{l_{2}}Q_{r+O(1)}(\partial_{t}\nabla^{-1}\psi\\&\sum_{a=1}^{3}P_{\geq l_{2}-100}\triangle^{-1}\partial_{a}(\partial_{a}\nabla^{-1}\psi f'''(\nabla^{-1}\psi)))||_{L_{t}^{2}L_{x}^{2}}\\
&\leq 2^{-2j}\sum_{l_{1}=l_{2}+O(1)}\sum_{r>\max\{j+100,l_{1}+10\}}2^{k}||P_{l_{1}}Q_{r}\Box\nabla^{-1}\psi||_{L_{t}^{4}L_{x}^{4}}2^{-r}||\partial_{t}\nabla^{-1}\psi||_{L_{t}^{4}L_{x}^{4}}C2^{-l_{2}}\\&||\partial_{a}\nabla^{-1}\psi||_{L_{t}^{\infty}L_{x}^{3}}\\
&\leq
2^{-2j}\sum_{l_{1}=l_{2}+O(1)}\sum_{r>\max\{j+100,l_{1}+10\}}C
2^{k-l_{2}}2^{\frac{l_{1}-r}{4}}\leq
2^{-2j}C\\
\end{split}\nonumber\end{equation}
where we have used Bernstein's inequality for the last step.\\

The last term in \eqref{26} can be estimated easily by placing
$\partial^{\nu}\nabla^{-1}\psi, \partial^{\nu}\nabla^{-1}\psi $
into $L_{t}^{4}L_{x}^{4}$.
\end{proof}

Continuing with the proof of the Proposition, we split into
High-Low, Low-High and High-High interactions:

\begin{equation}\begin{split}
&P_{0}[f(\sum_{j=1}^{3}\triangle^{-1}\partial_{j}\tilde{\phi}_{j})\psi]
=
P_{-10,10}f(\sum_{j=1}^{3}\triangle^{-1}\partial_{j}\tilde{\phi}_{j})P_{<15}\psi
+P_{<-10}f(\sum_{j=1}^{3}\triangle^{-1}\partial_{j}\tilde{\phi}_{j})P_{-10,10}\psi\\
&+\sum_{k_{1}>10,
k_{2}=k_{1}+O(1)}P_{k_{1}}f(\sum_{j=1}^{3}\triangle^{-1}\partial_{j}\tilde{\phi}_{j})P_{k_{2}}\psi\\
\end{split}\nonumber\end{equation}

Introduce the notation
$\sum_{j=1}^{3}\triangle^{-1}\partial_{j}\tilde{\phi}_{j}:=\Phi$.

\newpage

\begin{bf}High-Low Interactions.\end{bf}
\\

This case is the most elementary, as one can immediately employ
lemma~\ref{basic idea} in order to introduce a bilinear structure,
without incurring any losses.
\\

1): $L_{t}^{\infty}L_{x}^{2}$-norm of the output. We have

\begin{equation}\begin{split}
&||P_{0}\partial_{t}[P_{-10,10}f(\Phi)P_{<15}\psi]||_{L_{t}^{\infty}L_{x}^{2}}\\
&\leq||P_{0}[P_{-10,10}(\partial_{t}\Phi
f'(\Phi))P_{<15}\psi]||_{L_{t}^{\infty}L_{x}^{2}}\\
&+\sum_{a=1}^{3}||P_{0}[\triangle^{-1}\partial_{a}P_{-10,10}(\partial_{a}\Phi
f'(\Phi))P_{<15}\psi]||_{L_{t}^{\infty}L_{x}^{2}}\\
&\leq C\sum_{k_{1}\leq 0,
k_{2}<15}||P_{k_{1}}\nabla_{x,t}\Phi||_{L_{t}^{\infty}L_{x}^{4}}||P_{k_{2}}\psi||_{L_{t}^{\infty}L_{x}^{4}}\\
&+C\sum_{k_{1}>0,k_{2}<15}||P_{k_{1}}\nabla_{x,t}\Phi||_{L_{t}^{\infty}L_{x}^{2}}||P_{k_{2}}\psi||_{L_{t}^{\infty}L_{x}^{\infty}}\\
&\leq C\sum_{k_{1}\leq
0,k_{2}<15}2^{\frac{k_{1}+k_{2}}{4}}c_{k_{1}}c_{k_{2}} +
\sum_{k_{1}>0,
k_{2}<15}2^{-\frac{k_{1}}{2}}2^{k_{2}}c_{k_{1}}c_{k_{2}}\\
&\leq Cc_{0}^{2}\\
\end{split}\end{equation}

2): $\dot{X}_{0}^{\frac{1}{2},\frac{1}{2},\infty}$-norm of the
output.
\\

2.1) Output at modulation $2^{j}$ with $j\leq 100$:

\begin{equation}\begin{split}
&||P_{0}Q_{j}[P_{-10,10}f(\Phi)P_{<15}\psi]||_{\dot{X}_{0}^{\frac{1}{2},\frac{1}{2},\infty}}\\
&\leq\sum_{a=1}^{3}||P_{0}Q_{j}[\triangle^{-1}\partial_{a}P_{-10,10}(\partial_{a}\Phi
f'(\Phi))P_{<15}\psi]||_{\dot{X}_{0}^{\frac{1}{2},\frac{1}{2},\infty}}\\
&\leq C2^{\frac{j}{2}}\sum_{a=1}^{3}||\sum_{k_{1}\leq
-100}\triangle^{-1}\partial_{a}P_{-10,10}(\partial_{a}P_{k_{1}}\Phi
f'(\Phi))P_{<15}\psi||_{L_{t}^{2}L_{x}^{2}}\\
&+C2^{\frac{j}{2}}\sum_{a=1}^{3}||\sum_{k_{1}>-100}\triangle^{-1}\partial_{a}P_{-10,10}(\partial_{a}P_{k_{1}}\Phi
f'(\Phi))P_{<15}\psi||_{L_{t}^{2}L_{x}^{2}}\\
&\leq C2^{\frac{j}{2}}\sum_{a,b=1}^{3}||\sum_{k_{1}\leq
-100}\triangle^{-1}\partial_{a}P_{-10,10}(P_{k_{1}}\partial_{a}\Phi\triangle^{-1}\partial_{b}P_{-20,20}(\partial_{b}\Phi
f''(\Phi))P_{<15}\psi||_{L_{t}^{2}L_{x}^{2}}\\
&+C2^{\frac{j}{2}}\sum_{a=1}^{3}||\sum_{k_{1}>-100}\triangle^{-1}\partial_{a}P_{-10,10}(\partial_{a}P_{k_{1}}\Phi
f'(\Phi))P_{<15}\psi||_{L_{t}^{2}L_{x}^{2}}\\
\end{split}\end{equation}

The first of the two immediately preceding expressions can be
estimated by

\begin{equation}\begin{split}
& C\sum_{k_{1}\leq
-10}||P_{<-100}\nabla\Phi||_{L_{t}^{4}L_{x}^{4}}||\nabla\Phi||_{L_{t}^{4}L_{x}^{4}}||P_{<15}\psi||_{L_{t}^{\infty}L_{x}^{\infty}}\\
&\leq C\sum_{k<15}2^{k}c_{k}\leq Cc_{0}\\
\end{split}\end{equation}

by definition of a frequency envelope.\\
The 2nd term can be estimated by

\begin{equation}\begin{split}
&C\sum_{k_{1}>-100,
k_{2}<15}||P_{k_{1}}\nabla\Phi||_{L_{t}^{6}L_{x}^{3}}||P_{k_{2}}\psi||_{L_{t}^{3}L_{x}^{6}}\\
&\leq
C\sum_{k_{1}>-100,k_{2}<15}2^{\frac{k_{2}-k_{1}}{6}}c_{k_{1}}c_{k_{2}}\leq
C c_{0}^{2}\\
\end{split}\end{equation}

2.2) Output has modulation $2^{j}$ with $j>100$:

A): $P_{<15 }\psi$ has modulation $\geq 2^{j-100}$:

\begin{equation}\begin{split}
&||P_{0}Q_{j}[P_{-10,10}f(\Phi)P_{<15}Q_{\geq
j-100}\psi]||_{\dot{X}_{0}^{\frac{1}{2},\frac{1}{2},\infty}}\\
&\leq C2^{\frac{3j}{2}}||P_{<15}Q_{\geq j-100}\psi||_{L_{t}^{2}L_{x}^{2}}\\
&\leq C\sum_{k<15}C2^{\frac{k}{2}}c_{k}\leq Cc_{0}\\
\end{split}\end{equation}

B): $P_{<15}\psi$ has modulation $<2^{j-100}$; this implies that
$P_{-10,10}f(\Phi)$ has to be at modulation $\sim 2^{j}$:

\begin{equation}\begin{split}
&2^{\frac{j}{2}}||P_{0}Q_{j}\partial_{t}[P_{-10,10}Q_{j+O(1)}(f(\Phi))P_{<15}Q_{<j-100}\psi]||_{L_{t}^{2}L_{x}^{2}}\\
&\leq 2^{\frac{j}{2}}||P_{0}Q_{j}[P_{-10,10}Q_{j+O(1)}(f(\Phi))P_{<15}Q_{<j-100}\partial_{t}\psi]||_{L_{t}^{2}L_{x}^{2}}\\
&+2^{\frac{j}{2}}||P_{0}Q_{j}[P_{-10,10}Q_{j+O(1)}(\partial_{t}f(\Phi))P_{<15}Q_{<j-100}\psi]||_{L_{t}^{2}L_{x}^{2}}\\
\end{split}\end{equation}

The first summand is estimated by means of lemma~\ref{two-two}:

\begin{equation}\begin{split}
&2^{\frac{j}{2}}||P_{0}Q_{j}[P_{-10,10}Q_{j+O(1)}(f(\Phi))P_{<15}Q_{<j-100}\partial_{t}\psi]||_{L_{t}^{2}L_{x}^{2}}\\
&\leq
C2^{\frac{j}{2}}||P_{-10,10}Q_{j+O(1)}(f(\Phi))||_{L_{t}^{2}L_{x}^{2}}||P_{<15}Q_{<j-100}\partial_{t}\psi||_{L_{t}^{\infty}L_{x}^{\infty}}\\
&\leq Cc_{0}\\
\end{split}\end{equation}

The 2nd summand is estimated similarly.
\\

Note that in the preceding we have established the boundedness of
the $X^{\frac{1}{2},\frac{1}{2},1}$-norm of the part of
$P_{0}[P_{-10,10}(f(\Phi))P_{<15}\psi]$ with small modulation. In
particular, we control the third component of the $S[0]$-norm, and
are done with high-low interactions.

\newpage
\begin{bf} High-High Interactions.\end{bf}
\\

This case is more difficult. In particular, we have to cope with
the situation that two high-frequency terms whose Fourier supports
live very close to the light cone result in an output of small
frequency but very far away from the light cone, which renders the
$X^{\frac{1}{2},\frac{1}{2},\infty}$-component more difficult to
control. This case would be impossible to handle in $2$ space
dimensions.
\\

1): $L_{t}^{\infty}L_{x}^{2}$-norm of the output:

\begin{equation}\begin{split}
&||\partial_{t}(\sum_{k_{1}>10,
k_{1}=k_{2}+O(1)}P_{0}[P_{k_{1}}(f(\Phi))P_{k_{2}}\psi])||_{L_{t}^{\infty}L_{x}^{2}}\\
&\leq C\sum_{k_{1}>10,
k_{1}=k_{2}+O(1)}||P_{k_{1}}(f(\Phi))\partial_{t}P_{k_{2}}\psi||_{L_{t}^{\infty}L_{x}^{2}}\\
&+C\sum_{k_{1}>10, k_{1}=k_{2}+O(1)}||P_{k_{1}}(\partial_{t}\Phi
f'(\Phi))P_{k_{2}}\psi||_{L_{t}^{\infty}L_{x}^{2}}\\
&\leq C\sum_{k_{1}>10,
k_{1}=k_{2}+O(1)}||\sum_{a=1}^{3}\triangle^{-1}P_{k_{1}}\partial_{a}(\partial_{a}\Phi
f'(\Phi))\partial_{t}P_{k_{2}}\psi||_{L_{t}^{\infty}L_{x}^{2}}\\
&+C\sum_{k_{1}>10, k_{1}=k_{2}+O(1)}||P_{k_{1}}(\partial_{t}\Phi
f'(\Phi))P_{k_{2}}\psi||_{L_{t}^{\infty}L_{x}^{2}}\\
&\leq\sum_{k_{1}>10,
k_{1}=k_{2}+O(1)}C||\nabla_{x,t}\Phi||_{L_{t}^{\infty}L_{x}^{3}}(2^{-k_{1}}||\partial_{t}P_{k_{2}}\psi||_{L_{t}^{\infty}L_{x}^{2}}+||P_{k_{2}}\psi||_{L_{t}^{\infty}L_{x}^{2}})\\
&\leq C\sum_{k_{1}>10}2^{-\frac{k_{1}}{2}}c_{k_{1}}\leq Cc_{0}\\
\end{split}\end{equation}

2): $\dot{X}_{0}^{\frac{1}{2},\frac{1}{2},\infty}$-norm of the
output: We can easily control the
$\dot{X}_{0}^{\frac{1}{2},\frac{1}{2},1}$-norm of the output
restricted to small modulations, which in particular controls the
third component of $S[0]$ for the output. Hence consider now the
case when the modulation $2^{j}$ of the output is very large, i.e.
$j>>1$. Split the output as follows:

\begin{equation}\label{tedious}\begin{split}
&\sum_{k_{1}=k_{2}+O(1)>10}P_{0}Q_{j}(P_{k_{1}}f(\Phi)P_{k_{2}}\psi)\\
&=\sum_{k_{1}=k_{2}+O(1),j+10\geq k_{1}>10}P_{0}Q_{j}(P_{k_{1}}Q_{<j-10}f(\Phi)P_{k_{2}}Q_{<j-10}\psi)\\
&+\sum_{k_{1}=k_{2}+O(1),j+10\geq k_{1}>10}P_{0}Q_{j}(P_{k_{1}}Q_{\geq j-10}f(\Phi)P_{k_{2}}\psi)\\
&+\sum_{k_{1}=k_{2}+O(1),j+10\geq k_{1}>10}P_{0}Q_{j}(P_{k_{1}}Q_{<j-10}f(\Phi)P_{k_{2}}Q_{\geq j-10}\psi)\\
&+\sum_{k_{1}=k_{2}+O(1)>j+10}P_{0}Q_{j}(P_{k_{1}}f(\Phi)P_{k_{2}}Q_{<j-10}\psi)\\
&+\sum_{k_{1}=k_{2}+O(1)>j+10}(P_{k_{1}}f(\Phi)P_{k_{2}}Q_{\geq j-10}\psi)\\
\end{split}\end{equation}

The 2nd term in the immediately preceding is the most difficult,
as we cannot employ the $X^{\frac{1}{2},\frac{1}{2},\infty}$-norm
of the inputs. Instead, we will have to resort to angular
localization of the inputs, and exploit the third component of
$S[k]$:
\\

A): $\sum_{k_{1}=k_{2}+O(1),j+10\geq
k_{1}>10}P_{0}Q_{j}(P_{k_{1}}Q_{<j-10}f(\Phi)P_{k_{2}}Q_{<j-10}\psi)$:
First, observe that if $(\tau_{1},\xi_{1}),\,(\tau_{2},\xi_{2})$
are points in the Fourier supports of the inputs
$P_{k_{1}}Q_{<j-10}f(\Phi)$, $P_{k_{2}}Q_{<j-10}\psi$ respectively
that contribute toward the output, we need to have

\begin{equation}
2^{j}\sim ||\tau|-|\xi||=|\tau_{1}\pm |\xi_{1}|+\tau_{2}\pm
|\xi_{2}|\mp|\xi_{1}|\mp|\xi_{2}|\pm |\xi_{1}+\xi_{2}||
\nonumber\end{equation}

where signs have been chosen in such fashion that
$|\tau_{i}|-|\xi_{i}|| = |\tau_{i}\pm |\xi_{i}||,
||\tau_{1}+\tau_{2}|-|\xi_{1}+\xi_{2}|| = |\tau_{1}+\tau_{2}\pm
|\xi_{1}+\xi_{2}||,\, i=1,2$. Now $|\xi_{1}+\xi_{2}|\sim 2^{0}$,
whence by our assumptions on the Fourier supports of the inputs we
conclude  that the signs of $\mp|\xi_{1}|,\mp|\xi_{2}|$ must be
identical; moreover $k_{1}=k_{2}+O(1) = j+O(1)$. Observe in the
sequel that in this case $P_{k_{1}}Q_{<j-10}$ is a disposable
operator, viz. the definition in section 3.
\\

We need to estimate

\begin{equation}\begin{split}
& 2^{\frac{j}{2}}||\sum_{k_{1}=k_{2}+O(1),j+10\geq
k_{1}>10}P_{0}Q_{j}\partial_{t}(P_{k_{1}}Q_{<j-10}f(\Phi)P_{k_{2}}Q_{<j-10}\psi)||_{L_{t}^{2}L_{x}^{2}}\\
&\leq 2^{\frac{j}{2}}||\sum_{k_{1}=k_{2}+O(1),j+10\geq
k_{1}>10}P_{0}Q_{j}(P_{k_{1}}Q_{<j-10}f(\Phi)P_{k_{2}}Q_{<j-10}\partial_{t}\psi)||_{L_{t}^{2}L_{x}^{2}}\\
&+2^{\frac{j}{2}}||\sum_{k_{1}=k_{2}+O(1),j+10\geq
k_{1}>10}P_{0}Q_{j}(P_{k_{1}}Q_{<j-10}(\partial_{t}\Phi f'(\Phi))P_{k_{2}}Q_{<j-10}\psi)||_{L_{t}^{2}L_{x}^{2}}\\
&\leq 2^{\frac{j}{2}}||\sum_{k_{1}=k_{2}+O(1),j+10\geq
k_{1}>10}P_{0}Q_{j}(\sum_{a=1}^{3}\triangle^{-1}\partial_{a}P_{k_{1}}Q_{<j-10}(\partial_{a}\Phi f'(\Phi))\\&P_{k_{2}}Q_{<j-10}\partial_{t}\psi)||_{L_{t}^{2}L_{x}^{2}}\\
&+2^{\frac{j}{2}}||\sum_{k_{1}=k_{2}+O(1),j+10\geq
k_{1}>10}P_{0}Q_{j}(P_{k_{1}}Q_{<j-10}(\partial_{t}\Phi f'(\Phi))P_{k_{2}}Q_{<j-10}\psi)||_{L_{t}^{2}L_{x}^{2}}\\
\end{split}\end{equation}

We estimate the last term here, the last term but one being dealt
with in an identical manner.
\\

A.1): $\partial_{t}\Phi$ at modulation $<2^{j-100}$, $f'(\Phi)$ at
modulation $<2^{j-100}$ and frequency $<2^{j-100}$:\\
denote points in the Fourier supports of
$P_{<j-100}Q_{<j-100}f'(\Phi)$, $Q_{<j-100}\partial_{t}\Phi$  that
contribute toward the output
$P_{k_{1}}Q_{<j-10}(Q_{<j-100}\partial_{t}\Phi
P_{<j-100}Q_{<j-100}f'(\Phi))$ as \\ $(\tau_{b},\xi_{b}),\,
(\tau_{a},\xi_{a})$ respectively. First, note that $|\xi_{a}|\sim
2^{k_{1}}$. Next, $|\tau_{b}|\leq 2^{k_{1}-80}$, whence
$|\tau_{a}|\sim 2^{k_{1}}$. Also, $\tau_{a}$ has the same sign as
the $\tau_{1}$ of points $(\tau_{1},\xi_{1})$ in the Fourier
support of $P_{k_{1}}Q_{<j-10}(Q_{<j-100}\partial_{t}\Phi
P_{<j-100}Q_{<j-100}f'(\Phi))$. Thus recalling the comments of the
preceding paragraph, we can microlocalize
$Q_{<j-100}\partial_{t}\Phi$, $P_{k_{2}}Q_{<j-10}\psi$ to the same
half-space $\tau><0$.
\\

A.1.1): $P_{<j-100}Q_{<j-100}f'(\Phi)$ at frequency $\geq 2^{0}$:
fixing this frequency to be $2^{l},\,0\leq l\leq j-100$ for now,
we microlocalize the Fourier support of
$P_{k_{1}+O(1)}Q_{<j-100}\partial_{t}\Phi$ to a cap of size
$2^{\frac{3}{4}(l-k_{1})}$ when projected onto $S^{2}$. This then
implies that $P_{k_{2}}Q_{<j-10}\psi$ can be restricted to have
Fourier support in an approximately (up to $O(1)$ choices)
opposite cap. Now utilize the important bilinear inequality

\begin{equation}
||\phi\psi||_{L_{t}^{2}L_{x}^{2}}\leq
2^{\frac{k'-k}{2}}\frac{|\kappa'|}{dist(\kappa,\kappa')}||\phi||_{S[k,\kappa]}||\psi||_{S[k',\kappa']}
\nonumber\end{equation}

Hence

\begin{equation}\begin{split}
&||P_{0}Q_{j}[P_{k_{1}}Q_{<j-10}(\partial_{t}Q_{<j-100}\Phi
P_{j-100\geq . \geq
0}Q_{<j-100}(f'(\Phi)))P_{k_{2}}Q_{<j-10}\psi]||_{L_{t}^{2}L_{x}^{2}}\\
&\leq\sum_{\pm}\sum_{0\leq l\leq j-100}\sum_{\kappa,\kappa'\in
K_{\frac{3}{4}(l-k_{1})},\,dist(\kappa,-\kappa')\leq
C2^{\frac{3}{4}(l-k_{1})}}||P_{k_{1}}Q_{<j-10}(P_{k_{1}+O(1),\kappa}\partial_{t}Q^{\pm}_{<j-100}\Phi\\
&P_{l}Q_{<j-100}(f'(\Phi)))P_{k_{2},\kappa'}Q^{\pm}_{<j-10}\psi||_{L_{t}^{2}L_{x}^{2}}\\
\end{split}\end{equation}

a): Both $P_{k_{1}+O(1),\kappa}\partial_{t}Q^{\pm}_{<j-100}\Phi,
P_{k_{2},\kappa'}Q^{\pm}_{<j-10}\psi$ restricted to modulation
$<2^{\frac{3l-k_{1}}{2}}$: we can estimate
\begin{equation}\begin{split}
&\sum_{0\leq l\leq j-100}\sum_{\kappa,\kappa'\in
K_{\frac{3}{4}(l-k_{1})},\,dist(\kappa,-\kappa')\leq
C2^{\frac{3}{4}(l-k_{1})}}||P_{0}Q_{j}[P_{k_{1}}Q_{<j-10}(P_{k_{1}+O(1),\kappa}\partial_{t}Q^{\pm}_{<\frac{3l-k_{1}}{2}}\Phi\\
&P_{l}Q_{<j-100}(f'(\Phi)))P_{k_{2},\kappa'}Q^{\pm}_{<\frac{3l-k_{1}}{2}}\psi]||_{L_{t}^{2}L_{x}^{2}}\\
&\leq 4\sup_{\pm}\sum_{0\leq l\leq
j-100}2^{\frac{3}{4}(l-k_{1})}(\sum_{\kappa\in
K_{\frac{3}{4}(l-k_{1})}}||P_{k_{1}+O(1),\kappa}\partial_{t}Q^{\pm}_{<\frac{3l-k_{1}}{2}}\Phi||_{S[k_{1}+O(1),\pm\kappa]}^{2})^{\frac{1}{2}}\\
&(\sum_{\kappa\in
K_{\frac{3}{4}(l-k_{1})}}||P_{k_{2},\kappa}Q^{\pm}_{<\frac{3l-k_{1}}{2}}\Phi||_{S[k_{2},\pm\kappa]}^{2})^{\frac{1}{2}}\\
&||P_{l}Q_{<j-100}(f'(\Phi))||_{L_{t}^{\infty}L_{x}^{3}}\\
&\leq C\sup_{\pm}\sum_{0\leq l\leq
j-100}2^{\frac{3}{4}(l-k_{1})}[(\sum_{|a|\leq O(1)}\sum_{\kappa\in
K_{\frac{3}{4}(l-k_{1})}}||P_{k_{1}+a,\pm\kappa}\partial_{t}Q^{\pm}_{<\frac{3l-k_{1}}{2}+a}\Phi||_{S[k_{1}+a,\kappa]}^{2})^{\frac{1}{2}}\\
&+(\sum_{\kappa\in
K_{\frac{3}{4}(l-k_{1})}}||P_{k_{1}+a,\pm\kappa}\partial_{t}Q^{\pm}_{\frac{3l-k_{1}}{2}<.<\frac{3l-k_{1}}{2}+a}\Phi||_{\dot{X}_{k_{1}+a}^{\frac{1}{2},\frac{1}{2},1}}^{2})^{\frac{1}{2}}]\\
&[(\sum_{\kappa\in
K_{\frac{3}{4}(l-k_{1})}}||P_{k_{2},\pm\kappa}Q^{\pm}_{<\frac{3l-k_{2}}{2}}\psi||_{S[k_{2},\kappa]}^{2})^{\frac{1}{2}}+(P_{k_{2},\pm\kappa}Q^{\pm}_{\frac{3l-k_{2}}{2}<
.<
\frac{3l-k_{1}}{2}}\psi||_{\dot{X}_{k_{2}}^{\frac{1}{2},\frac{1}{2},1}}^{2})^{\frac{1}{2}}]\\
&||P_{l}Q_{<j-100}(f'(\Phi))||_{L_{t}^{\infty}L_{x}^{3}}\\
&\leq C\sum_{0\leq l\leq
j-100}2^{\frac{3}{4}(l-k_{1})}2^{-l}||P_{k_{1}+O(1)}\Phi||_{S[k_{1}+O(1)]}||P_{k_{2}}\psi||_{S[k_{2}]}\\
&\leq C\sum_{0\leq l\leq
j-100}2^{\frac{3}{4}(l-k_{1})}2^{-l}c_{k_{1}}c_{k_{2}}\leq C
2^{-\frac{3k_{1}}{4}}c_{k_{1}}^{2}\leq 2^{-\frac{3j}{4}}Cc_{j}^{2}\\
\end{split}\nonumber\end{equation}

because of $j=k_{1}+O(1)$.
\\

b): At least one of
$P_{k_{1}+O(1),\kappa}\partial_{t}Q^{\pm}_{<j-100}\Phi,
P_{k_{2},\kappa'}Q^{\pm}_{<j-10}\psi$ has modulation $\geq
2^{\frac{3l-k_{1}}{2}}$: this is handled by placing the input with
large modulation into $L_{t}^{2}L_{x}^{2}$ and the other input
into $L_{t}^{\infty}L_{x}^{2}$. One thereby obtains the upper
bound

\begin{equation}
\sum_{l\geq
0}C2^{-\frac{3l}{4}-\frac{3k_{1}}{4}}c_{k_{1}}c_{k_{2}}\leq
C2^{-\frac{3j}{4}}c_{j}^{2} 
\nonumber\end{equation}

Again, this suffices to establish the Proposition.
\\

A.1.2): $P_{<j-100}Q_{<j-100}f'(\Phi)$ at frequency $<2^{0}$: This
is handled as in the preceding case, except that now
$P_{k_{1}+O(1)}Q_{<j-100}\partial_{t}\Phi$, and
$P_{k_{2}}Q_{<j-10}\psi$ can be microlocalized in approximately
opposite caps of size $2^{-\frac{3k_{1}}{4}}$. Also,
$P_{<0}Q_{<j-100}(f'(\Phi))$ is estimated in the
$L_{t}^{\infty}L_{x}^{\infty}$-norm. Otherwise, the argument is
identical to the immediately preceding.
\\

A.2): $f'(\Phi)$ at frequency $\geq 2^{j-100}$:

\begin{equation}\begin{split}
&2^{\frac{j}{2}}||\sum_{k_{1}=k_{2}+O(1),k_{1}>10}P_{0}Q_{j}[P_{k_{1}}Q_{<j-10}(\partial_{t}\Phi
P_{\geq
j-100}(f'(\Phi)))P_{k_{2}}Q_{<j-10}\psi]||_{L_{t}^{2}L_{x}^{2}}\\
&\leq
2^{\frac{j}{2}}||\sum_{k_{1}=k_{2}+O(1),k_{1}>10}P_{0}Q_{j}[P_{k_{1}}Q_{<j-10}(\partial_{t}\Phi\sum_{a=1}^{3}\triangle^{-1}\partial_{a}P_{\geq
j-100}(\partial_{a}\Phi
f''(\Phi)))\\&P_{k_{2}}Q_{<j-10}\psi]||_{L_{t}^{2}L_{x}^{2}}\\
&\leq
C\sum_{k_{1}=j+O(1)=k_{2}+O(1)}2^{-\frac{j}{2}}||\partial_{t}\Phi||_{L_{t}^{\infty}L_{x}^{3}}||\partial_{a}\Phi||_{L_{t}^{4}L_{x}^{4}}||P_{k_{2}}Q_{<j-10}\psi||_{L_{t}^{4}L_{x}^{4}}\\
&\leq C\sum_{k_{1}=j+O(1)=k_{2}+O(1)}2^{-\frac{j}{2}}c_{k_{1}}\leq
Cc_{0}\\
\end{split}\end{equation}

A.3): $\partial_{t}\Phi$ at modulation $\geq 2^{j-100}$:

\begin{equation}\begin{split}
&2^{\frac{j}{2}}||\sum_{k_{1}=k_{2}+O(1),k_{1}>10}P_{0}Q_{j}[P_{k_{1}}Q_{\geq
j-100}\partial_{t}\Phi
P_{<k_{1}-100}Q_{<j-100}(f'(\Phi))P_{k_{2}}Q_{<j-10}\psi]||_{L_{t}^{2}L_{x}^{2}}\\
&\leq C\sum_{k_{1}=j+O(1)}2^{\frac{j}{2}}||P_{k_{1}}Q_{\geq
j-100}\partial_{t}\Phi||_{L_{t}^{2}L_{x}^{2}}||P_{k_{2}}Q_{<j-10}\psi||_{L_{t}^{\infty}L_{x}^{2}}\\
&\leq
C\sum_{k_{1}=j+O(1)}2^{\frac{j}{2}}2^{-\frac{j}{2}}2^{-k_{1}}c_{k_{1}}\leq
Cc_{0}\\
\end{split}\end{equation}

A.4): $f'(\Phi)$ at modulation $\geq 2^{j-100}$: use
lemma~\ref{two-two} to conclude:

\begin{equation}\begin{split}
&2^{\frac{j}{2}}||\sum_{k_{1}=k_{2}+O(1),k_{1}>10}P_{0}Q_{j}[P_{k_{1}}Q_{<j-10}(\partial_{t}\Phi
P_{<j-100}Q_{\geq
j-100}(f'(\Phi)))P_{k_{2}}Q_{<j-10}\psi||_{L_{t}^{2}L_{x}^{2}}\\
&\leq
2^{\frac{j}{2}}C\sum_{k_{1}=k_{2}+O(1)=j+O(1)}||P_{<j-100}Q_{\geq
j-100}(f'(\Phi))||_{L_{t}^{2}L_{x}^{\infty}}||P_{k_{2}}Q_{<j-10}\psi||_{L_{t}^{\infty}L_{x}^{2}}\\&||P_{k_{1}+O(1)}\partial_{t}\Phi||_{L_{t}^{\infty}L_{x}^{2}}\\
&\leq
C\sum_{k_{1}=k_{2}+O(1)=j+O(1)}2^{\frac{j}{2}}2^{j}2^{-\frac{3j}{2}}2^{-k_{2}}c_{k_{2}}\leq
Cc_{0}\\
\end{split}\end{equation}

Returning to the remaining terms of \eqref{tedious}:
\\

B): $\sum_{k_{1}=k_{2}+O(1),j+10\geq
k_{1}>10}P_{0}Q_{j}(P_{k_{1}}Q_{\geq j-10}f(\Phi)P_{k_{2}}\psi)$:

\begin{equation}\begin{split}
&2^{\frac{3j}{2}}||\sum_{k_{1}=k_{2}+O(1),j+10\geq
k_{1}>10}P_{0}Q_{j}[P_{k_{1}}Q_{\geq
j-10}(f(\Phi))P_{k_{2}}\psi]||_{L_{t}^{2}L_{x}^{2}}\\
&\leq C2^{\frac{3j}{2}}\sum_{k_{1}=k_{2}+O(1)}||P_{k_{1}}Q_{\geq
j-10}(f(\Phi))||_{L_{t}^{2}L_{x}^{2}}||P_{k_{2}}\psi||_{L_{t}^{\infty}L_{x}^{2}}\\
&\leq\sum_{10<k_{1}=k_{2}+O(1)\leq j+10}C2^{-k_{1}}c_{k_{1}}\leq
Cc_{0}\\
\end{split}\end{equation}

C): $\sum_{k_{1}=k_{2}+O(1),j+10\geq
k_{1}>10}P_{0}Q_{j}(P_{k_{1}}Q_{<j-10}f(\Phi)P_{k_{2}}Q_{\geq
j-10}\psi)$:

\begin{equation}\begin{split}
&2^{\frac{3j}{2}}||\sum_{k_{1}=k_{2}+O(1), j+10\geq
k_{1}>10}P_{0}Q_{j}[P_{k_{1}}Q_{<j-10}(f(\Phi))P_{k_{2}}Q_{\geq
j-10}\psi]||_{L_{t}^{2}L_{x}^{2}}\\
&\leq\sum_{k_{1}=k_{2}+O(1), j+10\geq k_{1}>10}\sum_{l\geq
j-10}C2^{\frac{3j}{2}}||P_{k_{1}}Q_{<j-10}(\sum_{a=1}^{3}\triangle^{-1}\partial_{a}(\partial_{a}\Phi
f'(\Phi))||_{L_{t}^{\infty}L_{x}^{3}}\\&||P_{k_{2}}Q_{l}\psi||_{L_{t}^{2}L_{x}^{2}}\\
&\leq\sum_{k_{1}=k_{2}+O(1), j+10\geq k_{1}>10}\sum_{l\geq
j-10}C2^{\frac{3j}{2}}2^{-k_{1}}2^{\frac{k_{2}}{2}}2^{-\frac{3l}{2}}c_{k_{2}}\\
&\leq Cc_{0}\\
\end{split}\end{equation}

D):
$\sum_{k_{1}=k_{2}+O(1)>j+10}P_{0}Q_{j}(P_{k_{1}}f(\Phi)P_{k_{2}}Q_{<j-10}\psi)$.

The argument is a little more convoluted here, since
$P_{k_{1}}Q_{j}$ isn't disposable anymore. Note that $f(\Phi)$ can
be restricted to modulation $\geq 2^{j+O(1)}$.

\begin{equation}\begin{split}
&2^{\frac{3j}{2}}\sum_{k_{1}=k_{2}+O(1)>j+10}||P_{0}Q_{j}[P_{k_{1}}f(\Phi)P_{k_{2}}Q_{<j-10}\psi]||_{L_{t}^{2}L_{x}^{2}}\\
&\leq
C2^{\frac{3j}{2}}\sum_{k_{1}=k_{2}+O(1)>j+10}||P_{k_{1}}Q_{\geq
j+O(1)}(f(\Phi))||_{L_{t}^{2}L_{x}^{2}}||P_{k_{2}}Q_{<j-10}\psi||_{L_{t}^{\infty}L_{x}^{2}}\\
&\leq
C2^{\frac{3j}{2}}\sum_{k_{1}=k_{2}+O(1)>j+10}||\sum_{a=1}^{3}\triangle^{-1}\partial_{a}P_{k_{1}}Q_{\geq
j+O(1)}[\partial_{a}P_{k_{1}+O(1)}Q_{\geq j+O(1)}\Phi
\\&P_{<j-100}Q_{<j-100}(f'(\Phi))]||_{L_{t}^{2}L_{x}^{2}}||P_{k_{2}}Q_{<j-10}\psi||_{L_{t}^{\infty}L_{x}^{2}}\\
&+C2^{\frac{3j}{2}}\sum_{k_{1}=k_{2}+O(1)>j+10}||\sum_{a=1}^{3}\triangle^{-1}\partial_{a}P_{k_{1}}Q_{\geq
j+O(1)}[\partial_{a}P_{k_{1}+O(1)}\Phi
\\&P_{<j-100}Q_{\geq j-100}(f'(\Phi))||_{L_{t}^{2}L_{x}^{2}}||P_{k_{2}}Q_{<j-10}\psi||_{L_{t}^{\infty}L_{x}^{2}}\\
&+C2^{\frac{3j}{2}}\sum_{k_{1}=k_{2}+O(1)>j+10}||\sum_{a=1}^{3}\triangle^{-1}\partial_{a}P_{k_{1}}Q_{\geq
j+O(1)}[\partial_{a}\Phi
P_{\geq j-100}(f'(\Phi))]||_{L_{t}^{2}L_{x}^{2}}\\
&||P_{k_{2}}Q_{<j-10}\psi||_{L_{t}^{\infty}L_{x}^{2}}\\
\end{split}\nonumber\end{equation}

Each of these terms is straightforward to estimate. We have

\begin{equation}\begin{split}
&2^{\frac{3j}{2}}\sum_{k_{1}=k_{2}+O(1)>j+10}||\sum_{a=1}^{3}\triangle^{-1}\partial_{a}P_{k_{1}}Q_{\geq
j+O(1)}[\partial_{a}P_{k_{1}+O(1)}Q_{\geq j+O(1)}\Phi
\\&P_{<j-100}Q_{<j-100}(f'(\Phi))]||_{L_{t}^{2}L_{x}^{2}}||P_{k_{2}}Q_{<j-10}\psi||_{L_{t}^{\infty}L_{x}^{2}}\\
&\leq
C2^{\frac{3j}{2}}\sum_{k_{1}=k_{2}+O(1)>j+10}2^{-k_{1}}||P_{k_{1}+O(1)}Q_{\geq
j+O(1)}
\nabla\Phi||_{L_{t}^{2}L_{x}^{2}}||P_{<j-100}Q_{<j-100}(f'(\Phi))||_{L_{t}^{\infty}L_{x}^{\infty}}\\
&||P_{k_{2}}Q_{<j-10}\psi||_{L_{t}^{\infty}L_{x}^{2}}\\
&\leq
C\sum_{k_{1}=k_{2}+O(1)>j+10}2^{\frac{3j}{2}}2^{-k_{1}}2^{-\frac{k_{2}}{2}}2^{-\frac{k_{1}}{2}}2^{-\frac{j}{2}}c_{k_{1}}\\
&\leq C\sum_{k_{1}=k_{2}+O(1)>j+O(1)}2^{j-2k_{1}}c_{k_{1}}\leq
Cc_{0}\\
\end{split}\nonumber\end{equation}

\begin{equation}\begin{split}
&2^{\frac{3j}{2}}\sum_{k_{1}=k_{2}+O(1)>j+10}||\sum_{a=1}^{3}\triangle^{-1}\partial_{a}P_{k_{1}}Q_{\geq
j+O(1)}[\partial_{a}P_{k_{1}+O(1)}\Phi
P_{<j-100}Q_{\geq j-100}(f'(\Phi))]||_{L_{t}^{2}L_{x}^{2}}\\
&||P_{k_{2}}Q_{<j-10}\psi||_{L_{t}^{\infty}L_{x}^{2}}\\
&\leq
C\sum_{k_{1}=k_{2}+O(1)>j+10}2^{\frac{3j}{2}}2^{-k_{1}}||P_{k_{1}+O(1)}\nabla\Phi||_{L_{t}^{\infty}L_{x}^{2}}||P_{<j-100}Q_{\geq
j-100}(f'(\Phi))||_{L_{t}^{2}L_{x}^{\infty}}\\&||P_{k_{2}}Q_{<j-10}\psi||_{L_{t}^{\infty}L_{x}^{2}}\\
&\leq
C\sum_{k_{1}=k_{2}+O(1)>j+10}2^{\frac{3j}{2}}2^{-k_{1}}c_{k_{1}}2^{j}2^{-\frac{3j}{2}}2^{-\frac{k_{2}}{2}}\\
&\leq C\sum_{k_{1}=k_{2}+O(1)>j+10}2^{-\frac{k_{1}}{2}}c_{k_{1}}\leq Cc_{0}\\
\end{split}\nonumber\end{equation}

Here we have used lemma~\ref{two-two}. \\
Finally

\begin{equation}\begin{split}
&2^{\frac{3j}{2}}\sum_{k_{1}=k_{2}+O(1)>j+10}||\sum_{a=1}^{3}\triangle^{-1}\partial_{a}P_{k_{1}}Q_{\geq
j+O(1)}[\partial_{a}\Phi
P_{\geq j-100}(f'(\Phi))]||_{L_{t}^{2}L_{x}^{2}}||P_{k_{2}}Q_{<j-10}\psi||_{L_{t}^{\infty}L_{x}^{2}}\\
&\leq\sum_{k_{1}=k_{2}+O(1)>j+10}C2^{\frac{3j}{2}}\sum_{a,b=1}^{3}2^{-k_{1}}||\partial_{a}\Phi||_{L_{t}^{4}L_{x}^{4}}||P_{\geq
j-100}\triangle^{-1}\partial_{b}(\partial_{b}\Phi
f''(\Phi))||_{L_{t}^{4}L_{x}^{4}}\\&||P_{k_{2}}Q_{<j-10}\psi||_{L_{t}^{\infty}L_{x}^{2}}\\
&\leq
\sum_{k_{1}=k_{2}+O(1)>j+10}C2^{\frac{3j}{2}-k_{1}-j}c_{k_{1}}\leq
Cc_{0}\\
\end{split}\nonumber\end{equation}

\newpage
\begin{bf} Low-High Interactions \end{bf}
\\

This case is hard as controlling the
$\dot{X}_{0}^{\frac{1}{2},\frac{1}{2},\infty}$-norm of the output
for small modulations forces us to utilize lemma~\ref{basic idea}
when $f(\Phi)$ is at small frequency, hence incurring an
exponential loss which we can only make up for by invoking another
case of angular localization and employing the third component of
$S[k]$.
\\

\begin{bf}1)\end{bf}: $L_{t}^{\infty}L_{x}^{2}$-norm of the output:

\begin{equation}\begin{split}
&||\partial_{t}[P_{<-10}(f(\Phi))P_{-10,10}\psi]||_{L_{t}^{\infty}L_{x}^{2}}\\
&\leq
C||P_{-10,10}\partial_{t}\psi||_{L_{t}^{\infty}L_{x}^{2}}+||P_{<-10}(\partial_{t}\Phi
f'(\Phi))||_{L_{t}^{\infty}L_{x}^{3}}||P_{-10,10}\psi||_{L_{t}^{\infty}L_{x}^{2}}\\
&\leq Cc_{0}\\
\end{split}\end{equation}

\begin{bf}2)\end{bf}: $\dot{X}_{0}^{\frac{1}{2},\frac{1}{2},\infty} $-norm of the output:
For high modulations of the output, this can be done exactly as in
the high-low case. Now assume that the modulation of the output is
$2^{j}\leq 2^{-10}$.

\begin{equation}\begin{split}
&P_{0}Q_{j}[P_{<-10}(f(\Phi))P_{-10,10}\psi]\\
&=P_{0}Q_{j}[P_{>j-100}(f(\Phi))P_{-10,10}\psi]\\
&+P_{0}Q_{j}[P_{\leq j-100}Q_{\geq j-100}(f(\Phi))P_{-10,10}\psi]\\
&+P_{0}Q_{j}[P_{\leq j-100}Q_{<j-100}f(\Phi)P_{-10,10}Q_{\geq
j-100}\psi]\\
\end{split}\end{equation}

1): $P_{0}Q_{j}[P_{>j-100}(f(\Phi))P_{-10,10}\psi]$: reformulate
this term as

\begin{equation}
\sum_{\tilde{j}>j-100}P_{0}Q_{j}(\sum_{a=1}^{3}\triangle^{-1}\partial_{a}P_{\tilde{j}}[\partial_{a}\Phi
f'(\Phi)]P_{-10,10}\psi)
\end{equation}

1.1): $f'(\Phi)$ at frequency $<2^{j-100}$ and modulation
$<2^{j-100}$, $\partial_{a}\Phi$ at modulation $<2^{j-100}$; also,
$\psi$ at modulation $<2^{j-100}$. This is the worst possible
scenario since we cannot introduce higher linearity by iterating
lemma~\ref{basic idea}. We will have to resort to angular
localization:

\begin{equation}\begin{split}
&\sum_{\tilde{j}>j-100}P_{0}Q_{j}(\sum_{a=1}^{3}\triangle^{-1}\partial_{a}P_{\tilde{j}}[\partial_{a}Q_{<j-100}\Phi
P_{<j-100}Q_{<j-100}f'(\Phi)]P_{-10,10}Q_{<j-100}\psi)\\
&=\sum_{\tilde{j}>j-100}P_{0}Q_{j}(\sum_{a=1}^{3}\triangle^{-1}\partial_{a}P_{\tilde{j}}[\partial_{a}P_{\tilde{j}+O(1)}Q_{<j-100}\Phi
P_{<j-100}Q_{<j-100}f'(\Phi)]\\&P_{-10,10}Q_{<j-100}\psi)\\
&=P_{0}Q_{j}[(\int_{\mathbf{R}^{3}}a(y)P_{\tilde{j}+O(1)}Q_{<j-100}\partial_{a}\Phi(x-y)P_{<j-100}Q_{<j-100}f'(\Phi)(x-y)dy)\\
&P_{-10,10}Q_{<j-100}\psi]\\
&=\int_{\mathbf{R}^{3}}a(y)P_{0}Q_{j}(P_{-15,15}Q_{j+O(1)}[P_{\tilde{j}+O(1)}Q_{<j-100}\partial_{a}\Phi(x-y)
P_{-10,10}Q_{<j-100}\psi(x)]\\&P_{<j-100}Q_{<j-100}f'(\Phi)(x-y))dy\\
\end{split}\nonumber\end{equation}

Here $a(x)$ denotes the kernel associated with the operator
$\triangle^{-1}\partial_{a}P_{\tilde{j}}$.
\\

The reason for reformulating the expression as above is that
provided the inputs
$P_{\tilde{j}+O(1)}Q_{<j-100}\partial_{a}\Phi(x-y)$,
$P_{-10,10}Q_{<j-100}\psi(x)$ of $[,]$ in the last term are
microlocalized to a half-space $\tau><0$, and such that the
projection of their Fourier supports to $S^{2}$ are supported on
caps $\kappa_{1},\kappa_{2}$ respectively, of radius
$2^{\frac{j-\tilde{j}}{2}-10}$, then $\pm\kappa_{1}$,
$\pm\kappa_{2}$ are at distance $\sim 2^{\frac{j-\tilde{j}}{2}}$.
The signs are chosen to be $+1$ for the upper half-space and $-1$
for the lower half-space. For this simple geometric fact see
lemma(13.2) in \cite{Tao 2}. Hence we have the identity:

\begin{equation}\label{muehsal}\begin{split}
&P_{-15,15}Q_{j+O(1)}[P_{\tilde{j}+O(1)}Q_{<j-100}\partial_{a}\Phi
P_{-10,10}Q_{<j-100}\psi(x)]\\
&=\sum_{\pm}\sum_{\kappa,\kappa'\in
K_{\frac{j-\tilde{j}}{2}},dist(\pm\kappa,\pm\kappa')\sim
2^{\frac{j-\tilde{j}}{2}}}P_{-15,15}Q_{j+O(1)}[P_{\tilde{j}+O(1),\kappa}Q^{\pm}_{<j-100}\partial_{a}\Phi
\\&P_{-10,10,\kappa'}Q^{\pm}_{<j-100}\psi(x)]\\
&=\sum_{\pm}\sum_{\kappa,\kappa'\in
K_{\frac{j-\tilde{j}}{2}},dist(\pm\kappa,\pm\kappa')\sim
2^{\frac{j-\tilde{j}}{2}}}P_{-15,15}Q_{j+O(1)}[P_{\tilde{j}+O(1),\kappa}Q^{\pm}_{<j-100}\partial_{a}\Phi
\\&P_{-10,10,\kappa'}Q^{\pm}_{\leq j-\tilde{j}}\psi(x)]\\
&-\sum_{\pm}\sum_{\kappa,\kappa'\in
K_{\frac{j-\tilde{j}}{2}},dist(\pm\kappa,\pm\kappa')\sim
2^{\frac{j-\tilde{j}}{2}}}P_{-15,15}Q_{j+O(1)}[P_{\tilde{j}+O(1),\kappa}Q^{\pm}_{<j-100}\partial_{a}\Phi
\\&P_{-10,10,\kappa'}Q^{\pm}_{j-\tilde{j}\geq .\geq j-100}\psi(x)]\\
\end{split}\end{equation}

Using as usual the bilinear inequality \eqref{crux}, we now
deduce:

\begin{equation}\begin{split}
&||\sum_{\pm}\sum_{\kappa,\kappa'\in K_{\frac{j-\tilde{j}}{2}-10},
dist(\pm\kappa,\pm\kappa')\sim
2^{\frac{j-\tilde{j}}{2}}}P_{-15,15}Q_{j+O(1)}[P_{\tilde{j}+O(1),\kappa}Q^{\pm}_{<j-100}\partial_{a}\Phi
\\&P_{-10,10,\kappa'}Q^{\pm}_{\leq j-\tilde{j}}\psi(x)]||_{L_{t}^{2}L_{x}^{2}}\\
&\leq
C\frac{2^{\frac{\tilde{j}}{2}}2^{\frac{j-\tilde{j}}{2}}}{2^{\frac{j-\tilde{j}}{2}}}\sum_{\kappa,\kappa'\in
K_{\frac{j-\tilde{j}}{2}-10, dist(\kappa,\kappa')\sim
2^{\frac{j-\tilde{j}}{2}}}}||P_{\tilde{j}+O(1),\pm\kappa}Q^{\pm}_{<j-100}\partial_{a}\Phi||_{S[\tilde{j},\kappa]}\\
&||P_{-10,10,\pm\kappa'}Q^{\pm}_{\leq
j-\tilde{j}}\psi(x)||_{S[0,\kappa']}\\
\end{split}\end{equation}

\begin{equation}\begin{split}
&\leq C2^{\frac{\tilde{j}}{2}}(\sum_{\kappa\in
K_{\frac{j-\tilde{j}}{2}-10}}||P_{\tilde{j}+O(1),\pm\kappa}Q^{\pm}_{<j-100}\partial_{a}\Phi||_{S[\tilde{j},\kappa]}^{2})^{\frac{1}{2}}\\
&(\sum_{\kappa\in
K_{\frac{j-\tilde{j}}{2}-10}}||P_{-10,10,\pm\kappa'}Q_{\leq
j-\tilde{j}}\psi(x)||_{S[0,\kappa']}^{2})^{\frac{1}{2}}\\
&\leq C2^{\frac{\tilde{j}}{2}}\sum_{|a-\tilde{j}|\leq O(1)}
\sum_{|b|\leq O(1)}[(\sum_{\kappa\in
K_{\frac{j-\tilde{j}}{2}-10}}||P_{a,\pm\kappa}Q^{\pm}_{<a+j-\tilde{j}-20}\partial_{a}\Phi||_{S[a,\kappa]}^{2})^{\frac{1}{2}}\\
&+(\sum_{\kappa\in
K_{\frac{j-\tilde{j}}{2}-10}}||P_{a,\pm\kappa}Q^{\pm}_{j-100<.<a+j-\tilde{j}-20}\partial_{a}\Phi||_{\dot{X}_{a}^{\frac{1}{2},\frac{1}{2},1}}^{2})^{\frac{1}{2}}]\\
&\times [(\sum_{\kappa\in
K_{\frac{j-\tilde{j}}{2}-10}}||P_{b,\pm\kappa}Q^{\pm}_{<j-\tilde{j}-20}\psi||_{S[b,\kappa]}^{2})^{\frac{1}{2}}\\
&+(\sum_{\kappa\in
K_{\frac{j-\tilde{j}}{2}-10}}||P_{b,\pm\kappa}Q^{\pm}_{j-\tilde{j}-20<.<j-\tilde{j}}\psi||_{\dot{X}_{b}^{\frac{1}{2},\frac{1}{2},1}}^{2})^{\frac{1}{2}}]\\
&\leq C2^{\frac{\tilde{j}}{2}}c_{0}c_{\tilde{j}}\\
\end{split}\end{equation}

Now as to the 2nd term in \eqref{muehsal}, note that

\begin{equation}\begin{split}
&||\sum_{\kappa,\kappa'\in
K_{\frac{j-\tilde{j}}{2}},dist(\pm\kappa,\pm\kappa')\sim
2^{\frac{j-\tilde{j}}{2}}}P_{-15,15}Q_{j+O(1)}[P_{\tilde{j}+O(1),\kappa}Q^{\pm}_{<j-100}\partial_{a}\Phi
\\&P_{-10,10,\kappa'}Q^{\pm}_{j-\tilde{j}\geq .\geq j-100}\psi(x)]||_{L_{t}^{2}L_{x}^{2}}\\
&\leq (\sum_{\kappa\in
K_{\frac{j-\tilde{j}}{2}-10}}||P_{\tilde{j}+O(1),\kappa}Q^{\pm}_{<j-100}\partial_{a}\Phi||_{L_{t}^{\infty}L_{x}^{\infty}}^{2})^{\frac{1}{2}}\\
&(\sum_{\kappa\in
K_{\frac{j-\tilde{j}}{2}-10}}||P_{-10,10,\kappa}Q^{\pm}_{j-\tilde{j}\geq
.\geq j-100}\psi(x)||_{L_{t}^{2}L_{x}^{2}}^{2})^{\frac{1}{2}}\\
\end{split}\nonumber\end{equation}

By means of Bernstein's and Strichartz' inequality, we know that

\begin{equation}\label{nice inequality}
||P_{k}Q_{j}\phi||_{L_{t}^{\infty}L_{x}^{\infty}}\leq
C2^{\frac{j}{4}}2^{\frac{3k}{4}}||\phi||_{L_{t}^{4}L_{x}^{4}}\leq
C2^{\frac{3j}{4}}2^{\frac{5k}{4}}||\phi||_{L_{t}^{2}L_{x}^{2}}
\end{equation}

Therefore

\begin{equation}\begin{split}
&||\sum_{\kappa,\kappa'\in
K_{\frac{j-\tilde{j}}{2}},dist(\pm\kappa,\pm\kappa')\sim
2^{\frac{j-\tilde{j}}{2}}}P_{-15,15}Q_{j+O(1)}[P_{\tilde{j}+O(1),\kappa}Q^{\pm}_{<j-100}\partial_{a}\Phi
\\&P_{-10,10,\kappa'}Q^{\pm}_{j-\tilde{j}\geq .\geq j-100}\psi(x)]||_{L_{t}^{2}L_{x}^{2}}\\
&\leq
C2^{-\frac{j}{4}+\frac{3\tilde{j}}{4}}||P_{\tilde{j}+O(1)}Q_{<j-100}\Phi||_{\dot{X}_{\tilde{j}+O(1)}^{\frac{1}{2},\frac{1}{2},\infty}}||P_{-10,10}Q_{j-\tilde{j}\geq
.\geq j-100}\psi||_{\dot{X}_{-10,10}^{\frac{1}{2},\frac{1}{2},\infty}}\\
&\leq C2^{-\frac{j}{4}+\frac{3\tilde{j}}{4}}c_{0}c_{\tilde{j}}\\
\end{split}\end{equation}

We can conclude that

\begin{equation}\begin{split}
&2^{\frac{j}{2}}||\sum_{\tilde{j}>j-100}P_{0}Q_{j}(\sum_{a=1}^{3}\triangle^{-1}\partial_{a}P_{\tilde{j}}[\partial_{a}\Phi
f'(\Phi)]P_{-10,10}\psi)||_{L_{t}^{2}L_{x}^{2}}\\
&\leq\sum_{\tilde{j}>j-100}C(2^{\frac{j-\tilde{j}}{2}}+2^{\frac{j-\tilde{j}}{4}})c_{\tilde{j}}c_{0}\leq
Cc_{0}\\
\end{split}\end{equation}

1.2): $f'(\Phi)$ at frequency $<2^{j-100}$ and modulation
$<2^{j-100}$, $\partial_{a}\Phi$ at modulation $\geq 2^{j-100}$,
$\psi$ at modulation $<2^{j-100}$:

\begin{equation}\begin{split}
&\sum_{\tilde{j}>j-100}P_{0}Q_{j}(\sum_{a=1}^{3}\triangle^{-1}\partial_{a}P_{\tilde{j}}[\partial_{a}P_{\tilde{j}+O(1)}Q_{\geq
j-100}\Phi P_{<j-100}Q_{<j-100}f'(\Phi)]P_{-10,10}Q_{<j-100}\psi)\\
&=\sum_{\tilde{j}>j-100}\sum_{l\geq
j-100}P_{0}Q_{j}(\sum_{a=1}^{3}\triangle^{-1}\partial_{a}P_{\tilde{j}}[P_{\tilde{j}+O(1)}Q_{l}\Phi
P_{<j-100}Q_{<j-100}f'(\Phi)]P_{-10,10}Q_{<j-100}\psi)\\
\end{split}\nonumber\end{equation}

The $\dot{X}_{0}^{\frac{1}{2},\frac{1}{2},\infty}$-norm of this
can be majorized by

\begin{equation}
C\sum_{\tilde{j}\geq l\geq
j-100}2^{\frac{j}{2}}2^{-\tilde{j}}c_{\tilde{j}}2^{\tilde{j}}2^{-\frac{l}{2}}2^{\frac{l-\tilde{j}}{4}}c_{0}
+\sum_{l>\tilde{j}}2^{\frac{j}{2}}2^{-\tilde{j}}c_{\tilde{j}}2^{\tilde{j}}2^{-\frac{l}{2}}c_{0}\leq
Cc_{0} \nonumber\end{equation}

We have used here the improved Bernstein's inequality
\eqref{improved Bernstein}.
\\

1.3): $f'(\Phi)$ has frequency $\geq 2^{\tilde{j}-10}$:

\begin{equation}\begin{split}
&2^{\frac{j}{2}}||\triangle^{-1}\partial_{a}P_{\tilde{j}}[\partial_{a}\Phi
P_{\geq
\tilde{j}-10}(f'(\Phi))]P_{-10,10}Q_{<j-100}\psi||_{L_{t}^{2}L_{x}^{2}}\\
&\leq
\sum_{b=1}^{3}2^{\frac{j}{2}}||\triangle^{-1}\partial_{a}P_{\tilde{j}}[\partial_{a}\Phi
P_{\geq \tilde{j}-10}\triangle^{-1}\partial_{b}(\partial_{b}\Phi
f''(\Phi))]P_{-10,10}Q_{<j-100}\psi||_{L_{t}^{2}L_{x}^{2}}\\
&\leq
\sum_{b=1}^{3}2^{\frac{j}{2}}||\triangle^{-1}\partial_{a}P_{\tilde{j}}[\partial_{a}\Phi
P_{\geq \tilde{j}-10}\triangle^{-1}\partial_{b}(\partial_{b}\Phi
f''(\Phi))]||_{L_{t}^{2}L_{x}^{\infty}}||P_{-10,10}Q_{<j-100}\psi||_{L_{t}^{\infty}L_{x}^{2}}\\
&\leq C2^{\frac{j-\tilde{j}}{2}}c_{0}\\
\end{split}\nonumber\end{equation}

by placing $\partial_{a}\Phi$, $\partial_{b}\Phi$ into
$L_{t}^{4}L_{x}^{4}$, and $P_{-10,10}Q_{<j-100}\psi$ into
$L_{t}^{\infty}L_{x}^{2}$. This can be summed over
$\tilde{j}>j-100$.
\\

1.4): $f'(\Phi)$ has frequency between $2^{j-100}$ and
$2^{\tilde{j}-10}$:

\begin{equation}\begin{split}
&2^{\frac{j}{2}}||\triangle^{-1}\partial_{a}P_{\tilde{j}}[\partial_{a}\Phi
P_{j-100<.<\tilde{j}-10}(f'(\Phi))]P_{-10,10}Q_{<j-10}\psi||_{L_{t}^{2}L_{x}^{2}}\\
&\leq
\sum_{b=1}^{3}2^{\frac{j}{2}}||\triangle^{-1}\partial_{a}P_{\tilde{j}}[\partial_{a}\Phi
P_{j-100<.<\tilde{j}-10}\triangle^{-1}\partial_{b}(\partial_{b}\Phi
f''(\Phi))]\\&P_{-10,10}Q_{<j-10}\psi||_{L_{t}^{2}L_{x}^{2}}\\
&\leq
\sum_{b=1}^{3}2^{\frac{j}{2}}||\triangle^{-1}\partial_{a}P_{\tilde{j}}[\partial_{a}P_{\tilde{j}}\Phi
P_{j-100<.<\tilde{j}-10}\triangle^{-1}\partial_{b}(\partial_{b}\Phi
f''(\Phi))]||_{L_{t}^{2}L_{x}^{\infty}}\\&||P_{-10,10}Q_{<j-10}\psi||_{L_{t}^{\infty}L_{x}^{2}}\\
&\leq
\sum_{b=1}^{3}C2^{\frac{j}{2}}2^{-\tilde{j}}||P_{\tilde{j}}\partial_{a}\Phi||_{L_{t}^{4}L_{x}^{\infty}}
||P_{j-100<.<\tilde{j}-10}\triangle^{-1}\partial_{b}(\partial_{b}\Phi
f''(\Phi))||_{L_{t}^{4}L_{x}^{\infty}}\\&||P_{-10,10}Q_{<j-10}\psi||_{L_{t}^{\infty}L_{x}^{2}}\\
&\leq C2^{\frac{j-\tilde{j}}{4}}c_{0}\\
\end{split}\nonumber\end{equation}

Again this can be summed over $\tilde{j}>j-100$ to yield a bound
of the required form.
\\

1.5): $f'(\Phi)$ has modulation $\geq 2^{j-100}$ and frequency
$<2^{j-100}$. Split into the cases $P_{\tilde{j}+O(1)}Q_{\geq
j-100}\partial_{a}\Phi$ and
$P_{\tilde{j}+O(1)}Q_{<j-100}\partial_{a}\Phi$. The first is dealt
with as in 1.2), taking into account the disposability of the
operator $P_{<j-100}Q_{\geq j-100}$. Now

\begin{equation}\begin{split}
&2^{\frac{j}{2}}||\triangle^{-1}\partial_{a}P_{\tilde{j}}[P_{\tilde{j}+O(1)}Q_{<j-100}(\partial_{a}\Phi)P_{<j-100}Q_{\geq
j-100}(f'(\Phi))]P_{-10,10}Q_{<j-100}\psi||_{L_{t}^{2}L_{x}^{2}}\\
&\leq
C2^{-\tilde{j}}2^{\frac{j}{2}}||P_{\tilde{j}+O(1)}Q_{<j-100}\partial_{a}\Phi||_{L_{t}^{\infty}L_{x}^{\infty}}||P_{<j-100}Q_{\geq
j-100}(f'(\Phi))||_{L_{t}^{2}L_{x}^{\infty}}\\&
||P_{<j-100}Q_{<j-100}\psi||_{L_{t}^{\infty}L_{x}^{2}}\\
&\leq C2^{\frac{j-\tilde{j}}{4}}c_{\tilde{j}}c_{0}\\
\end{split}\nonumber\end{equation}

where we have used \eqref{nice inequality}. Summing over
$\tilde{j}\geq j-100$, this case is dealt with as well.
\\

2): $P_{0}Q_{j}[P_{\leq j-100}Q_{\geq
j-100}(f(\Phi))P_{-10,10}\psi]$: This is handled using
lemma~\ref{two-two}:

\begin{equation}\begin{split}
&2^{\frac{j}{2}}||P_{0}Q_{j}[P_{\leq j-100}Q_{\geq
j-100}(f(\Phi))P_{-10,10}\psi]||_{L_{t}^{2}L_{x}^{2}}\\
&\leq C2^{\frac{j}{2}}||P_{\leq j-100}Q_{\geq
j-100}(f(\Phi))||_{L_{t}^{2}L_{x}^{\infty}}\\
&||P_{-10,10}\psi||_{L_{t}^{\infty}L_{x}^{2}}\\
&\leq Cc_{0}\\
\end{split}\end{equation}

3): $P_{0}Q_{j}[P_{\leq j-100}Q_{<j-100}f(\Phi)P_{-10,10}Q_{\geq
j-100}\psi]$:

\begin{equation}\begin{split}
&2^{\frac{j}{2}}||P_{0}Q_{j}[P_{\leq
j-100}Q_{<j-100}f(\Phi)P_{-10,10}Q_{\geq
j-100}\psi]||_{L_{t}^{2}L_{x}^{2}}\\
&\leq C2^{\frac{j}{2}}||P_{-10,10}Q_{\geq
j-100}\psi||_{L_{t}^{2}L_{x}^{2}}\\
&\leq Cc_{0}\\
\end{split}\end{equation}

We have simply placed $P_{\leq j-100}Q_{<j-100}f(\Phi)$ into
$L_{t}^{\infty}L_{x}^{\infty}$.
\\

\begin{bf}3)\end{bf}
\\

Finally, we have to deal with the 3rd component of $S[0]$ in the
low-high interaction case, i.e. we need to estimate

\begin{equation}
\sup_{\pm}\sup_{l\leq -10}(\sum_{\kappa\in
K_{l}}||P_{0,\pm\kappa}Q^{\pm}_{<2l}(f(\Phi)\psi)||_{S[0,\kappa]}^{2})^{\frac{1}{2}}
\end{equation}

We split into 4 cases:
\\
1.a): $\psi$ has modulation $<2^{l}$, $f(\Phi)$ has frequency
$<2^{l-10}$ and modulation $<2^{-10}$, where we have fixed $l\leq
-10$:

\begin{equation}\begin{split}
&\sup_{\pm}(\sum_{\kappa\in
K_{l}}||P_{0,\pm\kappa}Q^{\pm}_{<2l}(P_{<l-10}Q_{<-10}[f(\Phi)]Q_{<2l}\psi)||_{S[0,\kappa]}^{2})^{\frac{1}{2}}\\
&\leq \sup_{\pm}\sum_{|a|<O(1)}(\sum_{\kappa\in
K_{l}}||\sum_{\kappa'\in
K_{l-10},\kappa'\subset\kappa}P_{0,\pm\kappa}Q^{\pm}_{<2l}[(P_{<l-10}Q_{<-10}f(\Phi))\\&P_{a,\pm\kappa'}Q^{\pm}_{<a+2l-20}\psi]||_{S[0,\kappa]}^{2})^{\frac{1}{2}}\\
&+\sup_{\pm}\sum_{|a|<O(1)}(\sum_{\kappa\in
K_{l}}||\sum_{\kappa'\in
K_{l-10},\kappa'\subset\kappa}P_{0,\pm\kappa}Q^{\pm}_{<2l}[(P_{<l-10}Q_{<-10}f(\Phi))\\&P_{a,\pm\kappa'}Q^{\pm}_{2l>.\geq a+2l-20}\psi]||_{S[0,\kappa]}^{2})^{\frac{1}{2}}\\
&\leq C\sup_{\pm}\sum_{|a|\leq O(1)}(\sum_{\kappa'\in
K_{l-10}}||P_{a,\pm\kappa'}Q^{\pm}_{<a+2l-20}\psi||_{S[a,\kappa']}^{2})^{\frac{1}{2}}\\
&+C\sup_{\pm}\sum_{|a|\leq O(1)}(\sum_{\kappa'\in
K_{l-10}}||P_{a,\pm\kappa'}Q^{\pm}_{2l>.\geq
a+2l-20}\psi||_{\dot{X}_{a}^{\frac{1}{2},\frac{1}{2},1}})^{\frac{1}{2}}\\
&\leq Cc_{0}\\
\end{split}\end{equation}

We have used here the facts that $||fg||_{S[k,\kappa]}\leq
C||f||_{L^{\infty}}||g||_{S[k,\kappa]}$,
$||\phi||_{S[k,\kappa]}\leq C||\phi||_{S[k'\kappa']}$ provided
$\kappa'\subset\kappa, k=k'+O(1)$, as well as the disposability of
$P_{0,\pm\kappa}Q_{<2l}$ and the inequality
$||P_{k}\phi||_{S[k,\kappa]}\leq
C||P_{k}\phi||_{\dot{X}_{k}^{\frac{1}{2},\frac{1}{2},1}}$. For
these facts, see \cite{Tao 2}. Also, we have used here that
$P_{k,\kappa}$ microlocalizes to a concentric cap inside $\kappa$,
of half its size.
\\

1.b): $\psi$ has modulation $<2^{l}$, $f(\Phi)$ has frequency
$<2^{l-10}$ and modulation $\geq 2^{-10}$, where we have fixed
$l\leq -10$: This is much more elementary using
lemma~\ref{two-two}, hence left out.

2): $\psi$ has modulation $\geq 2^{2l}$ and $f(\Phi)$ has
frequency $<2^{l-10}$: we have

\begin{equation}\begin{split}
&(\sum_{\kappa\in
K_{l}}||P_{0\pm\kappa}Q^{\pm}_{<2l}(P_{<l-10}[f(\Phi)]Q_{\geq
2l}\psi)||_{S[0,\kappa]}^{2})^{\frac{1}{2}}\\
&\leq C\sum_{r<2l}(\sum_{\kappa\in
K_{l}}||P_{0,\pm\kappa}Q^{\pm}_{r}(P_{<l-10}[f(\Phi)]Q_{\geq
2l}\psi)||_{\dot{X}_{0}^{\frac{1}{2},\frac{1}{2},1}}^{2})^{\frac{1}{2}}\\
&\leq
C\sum_{r<2l}2^{\frac{r}{2}}||P_{<l-10}[f(\Phi)]P_{-5,5}Q_{\geq
2l}\psi||_{L_{t}^{2}L_{x}^{2}}\\
&\leq C\sum_{r<2l}2^{\frac{r}{2}-l}c_{0}\leq Cc_{0}\\
\end{split}\end{equation}

3): $f(\Phi)$ has frequency $\geq 2^{l-10}$:

\begin{equation}\begin{split}
&(\sum_{\kappa\in K_{l}}||P_{0,\pm\kappa}Q^{\pm}_{<2l}(P_{\geq
l-10}[f(\Phi)]P_{-10,10}\psi)||_{S[0,\kappa]}^{2})^{\frac{1}{2}}\\
&\leq C\sum_{r<2l}||P_{0}Q^{\pm}_{r}(P_{\geq
l-10}[f(\Phi)]P_{-10,10}\psi)||_{\dot{X}_{0}^{\frac{1}{2},\frac{1}{2},1}}\\
&\leq
C\sum_{r<2l}2^{\frac{r}{2}}||\sum_{a=1}^{3}\triangle^{-1}\partial_{a}P_{\geq
l-10}(\partial_{a}\Phi
f'(\Phi))P_{-10,10}\psi||_{L_{t}^{2}L_{x}^{2}}\\
&\leq
\sum_{r<2l}C2^{\frac{r}{2}-l}||\sum_{a=1}^{3}\partial_{a}\Phi
f'(\Phi)||_{L_{t}^{4}L_{x}^{4}}||P_{-10,10}\psi||_{L_{t}^{4}L_{x}^{4}}\\
&\leq C\sum_{r<2l}2^{\frac{r}{2}-l}c_{0}\leq Cc_{0}\\
\end{split}\end{equation}

This finishes the low-high case, and the Proposition is
established.
\end{proof}

\begin{bf} Remark \end{bf}: Note that we have proved more than the
Proposition states: provided that $||P_{k}\psi||_{S[k]}\leq \mu
c_{k}$, $||P_{k}\tilde{\phi}||_{S[k]}\leq M c_{k}$, and
$\sqrt{\sum_{k} c_{k}^{2}}\leq \epsilon$, we have that

\begin{equation}
||P_{k}[f(\triangle^{-1}\partial_{j}\tilde{\phi}_{j})\psi]||_{S[k]}\leq
C(Q(M\epsilon) +1)\mu c_{k}
\end{equation}

for some polynomial $Q(x)$. This is the statement needed to close
the bootstrapping argument.

\newpage
\begin{lemma}(T.Tao)\nonumber:
Let $j\leq\min(k_{1},k_{2})+O(1)$. Then
\begin{equation}
||P_{k}(F\psi)||_{N[k]}\leq
C2^{-\delta_{1}|k-\max(k_{1},k_{2})|}2^{-\delta_{2}|j-\min(k_{1},k_{2})|}||F||_{\dot{X}_{k_{1}}^{\frac{1}{2},-\frac{1}{2},\infty}}
||\psi||_{S[k_{2}]}
\end{equation}
for all Schwartz functions $F$ with Fourier support at frequency
$2^{k_{1}}$ and modulation $2^{j}$ while $\psi$ is at frequency
$2^{k_{2}}$, $\delta_{1},\delta_{2}>0$.
\end{lemma}

\begin{proof}

We need only prove the high-high interaction case. Proceeding as
in Tao's paper, we split into the following cases:\\
\\
Rescale to $k_{1}=k_{2}+O(1)=0$, whence $k\leq O(1)$. Also, let
$C>>1$:
\\

\begin{bf}1)\end{bf}: The estimate for $P_{k}(FQ_{\geq j-C}\psi)$:
\\

\begin{bf}1a)\end{bf}: $j<100k$:
\begin{equation}\begin{split}
&2^{-\frac{k}{2}}||P_{k}(FQ_{\geq j-C}\psi)||_{L_{t}^{1}L_{x}^{2}}
\leq 2^{-\frac{k}{2}}||F||_{L_{t}^{2}L_{x}^{\infty}}||Q_{\geq
j-C}\psi||_{L_{t}^{2}L_{x}^{2}}\\&\leq 2^{-\frac{j}{2}}2^{\delta
j}2^{-\frac{k}{2}}
||F||_{L_{t}^{2}L_{x}^{2}}||\psi||_{S[k_{2}]}\\
\end{split}\end{equation}
where $\delta$ can be chosen to be $\frac{1}{4}$, by the improved
Bernstein's  inequality \eqref{improved Bernstein}. Clearly the
gain $2^{\frac{j}{4}}$ makes more than up for the
$2^{-\frac{k}{2}}$-loss.
\\

\begin{bf}1b)\end{bf}: $j\geq 100k$: use Bernstein's inequality to
get
\begin{equation}\begin{split}
&2^{-\frac{k}{2}}||P_{k}(P_{k_{1}}FP_{k_{2}}Q_{\geq
j-C}\psi||_{L_{t}^{1}L_{x}^{2}}\leq C
2^{-\frac{k}{2}}2^{\frac{3k}{2}}||P_{k}(P_{k_{1}}F
P_{k_{2}}Q_{\geq
j-C}\psi)||_{L_{t}^{1}L_{x}^{1}}\\
&\leq 2^{k} ||F||_{L_{t}^{2}L_{x}^{2}}||Q_{\geq
j-C}\psi||_{L_{t}^{2}L_{x}^{2}}\leq
2^{k}2^{-\frac{j}{2}}||F||_{L_{t}^{2}L_{x}^{2}}||\psi||_{S[k_{2}]}\\
\end{split}\end{equation}
This is acceptable.
\\

\begin{bf}2)\end{bf}: The estimate for $P_{k}Q_{\geq
j-C}(FQ_{<j-C}\psi)$:

\begin{bf}2a)\end{bf}: $j<100k$: \\
\begin{equation}\begin{split}
&2^{-\frac{j}{2}}2^{-\frac{k}{2}}||P_{k}Q_{\geq
j-C}(FQ_{<j-C}\psi)||_{L_{t}^{2}L_{x}^{2}}\leq
2^{-\frac{j}{2}}2^{-\frac{k}{2}}||F||_{L_{t}^{2}L_{x}^{\infty}}||Q_{<j-C}\psi||_{L_{t}^{\infty}L_{x}^{2}}\\
&\leq C 2^{-\frac{k}{2}}2^{-\frac{j}{2}}2^{\delta j}
||F||_{L_{t}^{2}L_{x}^{2}}||\psi||_{S[k_{2}]}\\
\end{split}\end{equation}
Once again, we use that $\delta =\frac{1}{4}$, in order to
conclude this case.
\\

\begin{bf}2b)\end{bf}: $100k\leq j<O(1)$:\\
Use Bernstein's inequality:
\begin{equation}\begin{split}
&2^{-\frac{j}{2}}2^{-\frac{k}{2}}||P_{k}Q_{\geq
j-C}(FQ_{<j-C}\psi)||_{L_{t}^{2}L_{x}^{2}}\leq
C2^{-\frac{j}{2}}2^{k}2^{-\frac{k}{2}}||P_{k}Q_{\geq
j-C}(FQ_{<j-C}\psi)||_{L_{t}^{2}L_{x}^{1}}\\
&\leq 2^{-\frac{j}{2}}2^{k}||F||_{L_{t}^{2}L_{x}^{2}}||\psi||_{S[k_{2}]}\\
\end{split}\end{equation}
This is again acceptable.\\

\begin{bf}3)\end{bf}: The estimate for
$P_{k}Q_{<j-C}(FQ_{<j-C}\psi)$: this is the only case where we
have to work a little more carefully, and deviate from Tao's
proof: note that $j\leq k+O(1)$ here. The idea in this case is to
use the angular separation of the inputs $F, Q_{<j-C}\psi$, which
is obtained from Tao's lemma(13.2) for the "imbalanced case": this
allows us to conclude that upon localizing $F, Q_{<j-C}$ in  caps
of size $2^{-10}2^{k}2^{\frac{j-k}{2}}$, their angular separation
in Tao's signed sense has to be $\sim 2^{\frac{j+k}{2}}$. In
particular, fixing such a cap for $F$, there can be at most
$O(1)$-possible caps for $Q_{<j-C}\psi$. Moreover, we can also
conclude from that lemma that the angular separation between the
output and $F$ is $\sim 2^{\frac{j-k}{2}}$. Now utilize the
fundamental inequality
\begin{equation}
||\phi\psi||_{NFA[\kappa]}\leq
C\frac{|\kappa'|^{\frac{1}{2}}2^{\frac{k'}{2}}}{2^{\frac{k}{2}}dist(\kappa,\kappa')}||\phi||_{L_{t}^{2}L_{x}^{2}}||\psi||_{S[k',\kappa']}
\end{equation}
to estimate the previous by (letting $l_{1}=\frac{k+j}{2}-10,
l_{2}=\frac{j-k}{2}-10$)
\begin{equation}\begin{split}
&||P_{k}Q_{<j-C}(FQ_{<j-C}\psi)||_{N[k]}\\&\leq
2^{-\frac{k}{2}}\sum_{\pm,\pm,\pm}(\sum_{\kappa''\in
K_{l_{2}}}||\sum_{\kappa\in
K_{l_{1}},\,dist(\pm\kappa'',\pm\kappa)\sim 2^{l_{2}}}
\sum_{\kappa'\in K_{l_{1}},\,dist(\pm\kappa,\pm\kappa')\sim
2^{l_{1}}}\\&
P_{k,\kappa''}Q^{\pm}_{<j-C}(P_{k_{1},\kappa}Q^{\pm}_{j}FP_{k_{2},\kappa'}Q^{\pm}_{<j-C}\psi)||_{NFA[\kappa'']}^{2})^{\frac{1}{2}}\\
&\leq
C2^{-\frac{k}{2}}\frac{2^{\frac{k+j}{2}}}{2^{\frac{j-k}{2}}}(\sum_{\kappa''\in
K_{l_{2}}}[\sum_{\kappa\in
K_{l_{1}},\,dist(\pm\kappa'',\pm\kappa)\sim
2^{l_{2}}}||P_{k_{1},\kappa}Q^{\pm}_{j}F||_{L_{t}^{2}L_{x}^{2}}^{2}]\\&[\sum_{\kappa'\in
K_{l_{1}},\,dist(\pm\kappa',\pm\kappa'')\sim
2^{l_{2}}}||P_{k_{2},\kappa'}Q^{\pm}_{<j-C}\psi||_{S[k_{1},\pm\kappa']}^{2}])^{\frac{1}{2}}\\
&\leq C
2^{\frac{k}{2}}2^{\frac{j}{2}}[2^{-\frac{j}{2}}||Q_{j}F||_{L_{t}^{2}L_{x}^{2}}](\sum_{\kappa'\in
K_{l_{1}}}||P_{k_{2},\pm\kappa'}Q^{\pm}_{<j+k+k_{2}-20}\psi||_{S[k_{2},\kappa']}^{2})^{\frac{1}{2}}\\
&+2^{\frac{k}{2}}|k|2^{\frac{j}{2}}[2^{-\frac{j}{2}}||Q_{j}F||_{L_{t}^{2}L_{x}^{2}}]
||Q_{j-C>.\geq
j+k+k_{2}-20}P_{k_{2}}\psi||_{\dot{X}_{k_{2}}^{\frac{1}{2},\frac{1}{2},\infty}}
\end{split}\end{equation}
where we have discarded the summation over $\kappa''$ on account
of the distance condition $dist(\kappa',\kappa'')\sim 2^{l_{2}}$.
Of course, the above is acceptable.
\end{proof}

\begin{lemma}\nonumber(T.Tao)
Let $\phi$,$\psi$ be Schwarz functions on $\mathbf{R}^{3+1}$. Then
letting $k_{1}=k_{2}+O(1)\geq O(1)$, we have
\begin{equation}
||P_{0}[R_{\nu}P_{k_{1}}\phi\partial^{\nu}P_{k_{2}}\psi]||_{N[0]}\leq
C2^{-\delta
k_{1}}||P_{k_{1}}\phi||_{S[k_{1}]}||P_{k_{2}}\psi||_{S[k_{2}]}
\end{equation}
for some $\delta>0$.
\end{lemma}

\begin{proof}
We only have to prove this high-high interaction case. As usual,
we split into the cases corresponding to large modulations of the
inputs, when the
$\dot{X}_{k}^{\frac{1}{2},\frac{1}{2},\infty}$-component becomes
effective, and small modulations of the inputs, when the
$Q_{0}$-structure and lemma~\ref{bilinear1} kick in.
\\

1): $P_{0}[R_{\nu}P_{k_{1}}Q_{\geq k_{1}+100}
\phi\partial^{\nu}P_{k_{2}}Q_{<k_{2}-100}\psi]$:

\begin{equation}\begin{split}
&||P_{0}[R_{\nu}P_{k_{1}}Q_{\geq k_{1}+100}
\phi\partial^{\nu}P_{k_{2}}Q_{<k_{2}-100}\psi]||_{\dot{X}_{0}^{\frac{1}{2},-\frac{1}{2},1}}\\
&\leq\sum_{l\geq k_{1}+100}||P_{0}Q_{l+O(1)}[R_{\nu}P_{k_{1}}Q_{l}
\phi\partial^{\nu}P_{k_{2}}Q_{<k_{2}-100}\psi]||_{\dot{X}_{0}^{\frac{1}{2},-\frac{1}{2},1}}\\
&\leq\sum_{l\geq
k_{1}+100}C2^{-\frac{l}{2}}||R_{\nu}P_{k_{1}}Q_{l}
\phi||_{L_{t}^{2}L_{x}^{2}}||\partial^{\nu}P_{k_{2}}Q_{<k_{2}-100}\psi||_{L_{t}^{\infty}L_{x}^{2}}\\
&\leq\sum_{l\geq
k_{1}+100}C2^{-\frac{l}{2}}2^{-\frac{l}{2}}2^{-\frac{k_{1}}{2}}||P_{k_{1}}\phi||_{S[k_{1}]}2^{\frac{k_{2}}{2}}||P_{k_{2}}\psi||_{S[k_{2}]}\\
&\leq C2^{-k_{1}}||P_{k_{1}}\phi||_{S[k_{1}]}||P_{k_{2}}\psi||_{S[k_{2}]}\\
\end{split}\end{equation}

2): $P_{0}[R_{\nu}P_{k_{1}}Q_{\geq k_{1}+100}
\phi\partial^{\nu}P_{k_{2}}Q_{\geq k_{2}-100}\psi]$:

\begin{equation}\begin{split}
&||P_{0}[R_{\nu}P_{k_{1}}Q_{\geq k_{1}+100}
\phi\partial^{\nu}P_{k_{2}}Q_{\geq
k_{2}-100}\psi]||_{L_{t}^{1}L_{x}^{2}}\\
&\leq C||R_{\nu}P_{k_{1}}Q_{\geq k_{1}+100}
\phi||_{L_{t}^{2}L_{x}^{2}}||\partial^{\nu}P_{k_{2}}Q_{\geq
k_{2}-100}\psi||_{L_{t}^{2}L_{x}^{2}}\\
&\leq
C2^{-k_{1}}||P_{k_{1}}\phi||_{S[k_{1}]}||P_{k_{2}}\psi||_{S[k_{2}]}\\
\end{split}\end{equation}

3): $P_{0}[R_{\nu}P_{k_{1}}Q_{<k_{1}+100}
\phi\partial^{\nu}P_{k_{2}}Q_{\geq k_{2}+200}\psi]$: This is like
case 1).
\\

4):  $P_{0}[R_{\nu}P_{k_{1}}Q_{<k_{1}+100}
\phi\partial^{\nu}P_{k_{2}}Q_{<k_{2}+200}\psi]$: Use the
$Q_{0}$-structure to reduce this to the estimation of the
following terms:
\\

4.1): $\Box P _{0}[P_{k_{1}}Q_{<k_{1}+100} \frac{\phi}{\nabla}
P_{k_{2}}Q_{<k_{2}+200}\psi]$:

\begin{equation}\begin{split}
&||\Box P_{0}[P_{k_{1}}Q_{<k_{1}+100}
\frac{\phi}{\nabla}P_{k_{2}}Q_{<k_{2}+200}\psi]||_{N[0]}\\
&\leq ||\Box P_{0}Q_{<k_{1}+200}[P_{k_{1}}Q_{<k_{1}+100}
\frac{\phi}{\nabla}P_{k_{2}}Q_{<k_{2}+200}\psi]||_{\dot{X}_{0}^{\frac{1}{2},-\frac{1}{2},1}}\\
&\leq
C2^{-\frac{k_{1}}{2}}||P_{k_{1}}\phi||_{L_{t}^{4}L_{x}^{4}}||P_{k_{2}}\psi||_{L_{t}^{4}L_{x}^{4}}\\
&\leq
C2^{-\frac{k_{1}}{2}}||P_{k_{1}}\phi||_{S[k_{1}]}||P_{k_{2}}\psi||_{S[k_{2}]}\\
\end{split}\end{equation}

4.2): $P _{0}[P_{k_{1}}Q_{<k_{1}+100} \frac{\Box\phi}{\nabla}
P_{k_{2}}Q_{<k_{2}+200}\psi]$: use lemma~\ref{bilinear1} to
conclude that

\begin{equation}\begin{split}
&||P _{0}[P_{k_{1}}Q_{<k_{1}+100} \frac{\Box\phi}{\nabla}
P_{k_{2}}Q_{<k_{2}+200}\psi]||_{N[0]}\\
&\leq C2^{-\delta k_{1}}||P_{k_{1}}Q_{<k_{1}+100}
\frac{\Box\phi}{\nabla}||_{\dot{X}_{k_{1}}^{\frac{1}{2},-\frac{1}{2},\infty}}||P_{k_{2}}Q_{<k_{2}+200}\psi||_{S[k_{2}]}\\
&\leq C2^{-\delta
k_{1}}||P_{k_{1}}\phi||_{S[k_{1}]}||P_{k_{2}}\psi||_{S[k_{2}]}\\
\end{split}\end{equation}

4.3): $P _{0}[P_{k_{1}}Q_{<k_{1}+100} \frac{\phi}{\nabla}
P_{k_{2}}Q_{<k_{2}+200}\Box\psi]$: this is similar to the
preceding case.
\end{proof}

\end{document}